\documentclass{amsart}
\usepackage{amssymb,amscd}
\usepackage[french]{babel}
\usepackage[latin1]{inputenc}
\usepackage[T1]{fontenc}
\usepackage[11pt,curve,matrix,arrow,frame,graph]{xy}
\newcounter{intro}

\newtheorem{theo}[intro]{Théorème}

\newtheorem{thm}{Théorème}[section]
\newtheorem{lem}[thm]{Lemme}
\newtheorem{prop}[thm]{Proposition}
\newtheorem{cor}[thm]{Corollaire}

\theoremstyle{remark}
\newtheorem{rems}[thm]{Remarques}
\newtheorem{rem}[thm]{Remarque}
\newtheorem{defi}[thm]{Définition}
\newtheorem*{merci}{Remerciements}

\numberwithin{equation}{section}
\newcounter{counteroman}
\newenvironment{enumeroman}{\begin{list}{\roman{counteroman})}{\usecounter{counteroman}}}{\end{list}}

\newcommand{\cref}[1]{Corollary~\ref{#1}}

\newcommand{\R}{\mathbb{R}}
\newcommand{\C}{\mathbb{C}}
\renewcommand{\H}{\mathbb{H}}

\newcommand{\N}{\mathbb{N}}\newcommand{\bQ}{\mathbb{Q}}
\newcommand{\Z}{\mathbb{Z}}
\newcommand{\bS}{\mathbb{S}}
\newcommand{\bB}{\mathbb{B}}

\newcommand{\cH}{\mathcal{H}}

\newcommand{\cD}{\mathcal{D}}
\newcommand{\cC}{\mathcal{C}}
\newcommand{\cO}{\mathcal{O}}

\def\et{{\rm\ et\ }}
\def\si{{\rm\ si\ }}

\let\eps=\varepsilon

\DeclareMathOperator{\Sp}{Sp}
\DeclareMathOperator{\Su}{SU}

\DeclareMathOperator{\Id}{Id}
\DeclareMathOperator{\ima}{Im}

\DeclareMathOperator{\vol}{vol}
\DeclareMathOperator{\age}{age}
\DeclareMathOperator{\card}{card}
\DeclareMathOperator{\inte}{int}
\def\cro#1#2{\mathrel{\langle {#1},{#2}\rangle}}

\def\Lk{\Lambda^kT^*}
\def\Lkm{{\Lambda^{k-1}T^*}}
\def\Lkp{{\Lambda^{k+1}T^*}}
\def\loc{{\text loc}}

\begin{document}

\title[Cohomologie $L^2$ des variétés QALE]
{Cohomologie $L^2$ des variétés QALE}

\author{Gilles Carron}
\address{Laboratoire de Math\'ematiques Jean Leray (UMR 6629), 
Universit\'e de Nantes,
2, rue de la Houssini\`ere, B.P.~92208, 44322 Nantes Cedex~3, France}
\email{Gilles.Carron@math.univ-nantes.fr}


\date{18 janvier 2005}
\begin{abstract}
Nous donnons une interprétation topologique des espaces de formes harmoniques 
$L^2$ de certaines des variétés QALE introduites par D. Joyce. 
Nous introduisons pour cela un critère analytique qui nous permet d'utiliser
des bouts de suites exactes de Mayer-Vietoris.
\end{abstract}
\maketitle
\section{Introduction}

Sur une variété riemannienne compacte, le célèbre théorème de Hodge et de Rham identifie les espaces de
formes harmoniques $L^2$ aux groupes de cohomologie réels de la variété. Lorsque $(X,g)$ est une variété
complète et non compacte, ses espaces de formes harmoniques $L^2$
$$\cH^k(X)=\left\{\alpha\in L^2(\Lk X), d\alpha=d^*\alpha=0\right\}$$
peuvent être identifiés aux espaces de $L^2$ cohomologie réduite :
$$\cH^k(X)\simeq \H^k(X)=\left\{\alpha\in L^2(\Lk X), d\alpha=0\right\}/
\overline{ \left\{d\alpha,\alpha\in L^2(\Lkm X), d\alpha\in L^2 \right\} }$$
où l'adhérence est prise pour la topologie $L^2$. Cependant cette cohomologie $L^2$ réduite n'est
pas  en général pratique à calculer, car elle ne vérifie pas les suites exactes de Mayer-Vietoris. 
La cohomologie $L^2$ non réduite définie par
$${}^{nr}\H^k(X)=\left\{\alpha\in L^2(\Lk X),\ d\alpha=0\right\}/
\left\{d\alpha,\alpha\in L^2(\Lkm X),\ d\alpha\in L^2 \right\}$$
peut se calculer à l'aide de suites exactes de Mayer Vietoris, cependant ces espaces ne coïncident
que lorsque l'image de $d$ est fermée et dans le cas contraire les espaces de cohomologie $L^2$ non réduite 
sont de dimension infinie : ils ne permettent pas de récupérer la cohomologie $L^2$ réduite.

De nombreux travaux donnent une interprétation topologique de ces espaces de cohomologie $L^2$ réduite
 en fonction de la topologie et de la géométrie à l'infini. Par exemple, dans \cite{APS}, M. Atiyah,
  V. Patodi et I. Singer calculent la cohomologie $L^2$ réduite des variétés à bouts cylindriques, les travaux
  de A. Borel, B. Casselman, E. Looijenga, L. Saper, M. Stern et S. Zucker ont permis d'identifier la cohomologie $L^2$ des variétés
 hermitiennes localement symétriques
 \footnote{dans ce cas l'image de $d$ est fermée et cohomologie $L^2$
  réduite et non réduite coïncident.} à la cohomologie d'intersection de leurs compactifications de
  Borel-Serre-Satake (\cite{B1,BC,Loo,SS,Z1,Z2}), les variétés à courbure négative ou asymptotiquement nulle ont aussi fait l'objet de
  nombreux travaux (\cite{Lo,M,MP,Y,Cargafa}).
  
Récemment, A. Sen a utilisé une dualité issue de la M-théorie pour prédire l'existence
de formes harmoniques $L^2$ sur l'espace réduit des monopoles de charges $k$ sur $\R^3$ (\cite{S}). Grâce à
des arguments également issus de la physique théorique, il existe
aujourd'hui de nombreuses prédictions concernant les espaces de formes harmoniques $L^2$ sur certaines
variétés riemanniennes complètes à holonomie exceptionnelle.
Dans un papier remarquable, N. Hitchin montre qu'une variété hyperkählerienne complète 
de dimension réelle $4d$ obtenue par réduction hyperkählerienne de $\R^{4n}$ ne peut porter de formes harmoniques $L^2$ 
qu'en degré $2d$ \cite{H} ; pour cela il utilise une amélioration due à J.Jost et K.Zuo (\cite{JZ}) d'un
résultat de M. Gromov (\cite{Gro}
cf. aussi \cite{CX,Mc}). N. Hitchin détermine aussi les formes harmoniques $L^2$ de certaines de ces variétés dont
  l'espace réduit des monopoles de charge $2$ et il vérifie dans ce cas les prédictions de A. Sen. En fait
  les travaux de G. Segal et A. Selby montrent les prédictions de A. Sen équivalent à démontrer que pour ces variétés 
  la cohomologie $L^2$ réduite est isomorphe à l'image de la cohomologie à support compact dans la
  cohomologie absolue (\cite{SSe}).
 Un travail fondamental de T. Hausel, E. Hunsicker, R. Mazzeo interprète les espaces de cohomologie $L^2$
 réduite des variétés dont la géométrie à l'infini est de type "fibred boundary" (à bord fibré) ou "
 fibred cusp (à pointe fibré)  en terme de cohomologie d'intersection de la compactification obtenue en
 écrasant les fibres de la fibration à l'infini (\cite{HHM}). Ces géométries introduites par R. Mazzeo et
  R. Melrose modélisent les variétés localement symétrique de $\bQ$-rang $1$ et certaines géométries
   (ALE, ALF, ALG) des variétés à holonomie exceptionnelle (\cite{MM}) et ce travail leurs permet de confirmer ou d'infirmer certaines des
 prédictions issues de  physique théorique (cf. la septième partie de \cite{HHM} pour plus de précisions).
 
 L'objet de  ce travail est la classe de variétés à holonomie $\Su(n)$ ou $\Sp(n/2)$ construite par D. Joyce
 (théorème 9.3.3 de \cite{joyce}):
 
 \begin{theo} Soit $G\subset \Su (n)$ un sous groupe fini et $X\stackrel{\pi}{\longrightarrow} \C^n/G$ une
 résolution crépante de $\C^n/G$, si $X$ porte une métrique kählerienne QALE asymptote à $\C^n/G$ alors elle
 porte aussi une métrique Kähler-Einstein
 plate ; de plus si $G\subset \Sp (n/2)$ alors cette métrique est hyperkählerienne.
 \end{theo}

Nous ne donnons pas la définition de résolution crépante (cf. \cite{joyce} page 126 par exemple), nous
signalons uniquement qu'une telle résolution n'existe pas toujours. Nous dirons plus loin deux mots sur
la géométrie des variétés QALE (Quasi-Asymptotiquement Localement Euclidienne). Certaines de ces résolutions comme
les schémas de Hilbert sur $\C^2$ sont aussi connues sous le nom d'espace des instantons non commutatifs
(\cite{NS, Nakajima}). Notre résultat principal est le suivant :
  
\begin{theo}\label{cohocrepante} Soit $G\subset \Su (n)$ un sous groupe fini, on suppose que $\bS^{2n-1}/G
$ le quotient de la sphère par $G$ soit à singularités isolées. Si $X\stackrel{\pi}{\longrightarrow} \C^n/G$ une
 résolution crépante de $\C^n/G$ équipée d'une métrique QALE asymptote à $\C^n/G$ alors
 $$\H^k(X)\simeq \ima\left(\, H_c^k(X)\rightarrow H^k(X)\,\right).$$
 
\end{theo}
  
  Les travaux de Y. Ito et M. Reid, V. Batyrev, J. Denef et F. Loeser permettent
   de décrire la cohomologie usuelle de $X$ à 
  l'aide des classes de conjugaison de $G$ ((\cite{IR, Baty,DL}).
   Notons $\cC(G)$ l'ensemble des classes de conjugaison de $G$, en
  corollaire de ce théorème nous obtenons :
  \begin{cor}
  Soit $G\subset \Su (3)$ ou $G\subset \Sp(2)$ un sous groupe fini et $X\stackrel{\pi}{\longrightarrow} \C^n/G$ une
 résolution crépante de $\C^n/G$ équipée d'une métrique QALE asymptote à $\C^n/G$,
 alors :
 \begin{equation*}
 \begin{split}
 \chi_{L^2}(X)&=\sum_{k=0}^{2n} (-1)^k \dim \H^k(X)= \sum_{l=0}^{n}  \dim \H^{2l}(X)\\
 &=
 \card\big\{ [g]\in \cC(G), \ker(g-\Id)=\{0\} \,\big\}.\\
 \end{split}
 \end{equation*}
 \end{cor}
  
  En particulier, nous démontrons que l'espace des formes harmoniques $L^2$ de ${\rm Hilb}^3(\C^2)$ est de 
  dimension $1$ et qu'il est formé de formes de degré $4$.
  Notre théorème répond par l'affirmative à une question de E. Hunsicker. Lorsque $G$ agit sans point fixe
  sur $\bS^{2n-1}$ alors $X$ est une variété ALE : au dehors de $\pi^{-1}\{0\}$, $X$ est quasi isométrique
   au bout de cône $\{z\in \C^n, \|z\|\ge 1\}/G$ et dans ce cas ce théorème \ref{cohocrepante} est déjà connu (\cite{Melrose, carma}).
  
  Décrivons rapidement la géométrie à l'infini des variétés QALE considérées ici : nous supposons donc que
  $G\subset \Su (n)$ soit un sous groupe fini et que $\bS^{2n-1}/G$ soit à singularités isolées. 
  Notons $S\subset \bS^{2n-1}/G$ ce lieu singulier; c'est le quotient par $G$ d'une réunion de sphères. Si 
  $(X,g)$ est une variété QALE asymptote à $\C^n/G$ alors au dehors d'un compact $K\subset X$, 
  $$X\setminus K=\cup_{i=0}^l E_i$$
  où $E_0$ est quasi-isométrique à un bout cône sur $(\bS^{2n-1}/G)\,\setminus S^\varepsilon$ où
  $S^\varepsilon$ est un $\varepsilon$ voisinage de $S\subset \bS^{2n-1}/G$ et chaque $E_i$ a un revêtement
  fini quasi isométrique à $$\{(y,v)\in Y_i\times V_i,\ 1\le \|v\|,  \|\pi_i(y)\|\le 2\varepsilon  \}$$
  où $V_i\subset \C^n$ est un sous espace vectoriel, $A_i=\{g\in G,
  g|_{V_i}=\Id\}\not=\{\Id\}$ et $\pi_i\,:\,Y_i\rightarrow V_i^\perp/A_i$ est une
  résolution de $V_i^\perp/A_i$ équipée d'une métrique ALE asymptote à  $V_i^\perp/A_i$.
  
  L'outil important introduit ici est une suite de Mayer-Vietoris entre cohomologie $L^2$ à poids réduite.
  Pour alléger les discours on introduit la définition suivante :

  Soit $(X,g)$ une variété riemannienne (non nécessairement complète) et $\mu,w$ des fonctions 
  strictement positives
  (lisses) sur $X$, on suppose de plus que $w$ est bornée. On peut définir la cohomologie  à poids :
  
  $$\H_\mu^k(X)=\left\{\alpha\in L_\mu^2(\Lk X), d\alpha=0\right\}/
\overline{ \left\{d\alpha,\alpha\in L_\mu^2(\Lkm X), d\alpha\in L_\mu^2 \right\} }$$
  où l'adhérence est prise pour la topologie $L^2_\mu$ associée à la mesure $\mu d\vol_g$.
  On note aussi $\cD_\mu^{k-1}(d)=\left\{d\alpha,\alpha\in L_\mu^2(\Lkm X), d\alpha\in L_\mu^2 \right\}$ le
  domaine de $d$.
  
  \begin{defi}
  On dira que l'image de $d\,:\, \cD_\mu^{k-1}(d)\rightarrow L_\mu^2(\Lk X)$ est {\bf presque fermée}
 (en degré $k$ et par rapport à $w$) si lorsqu'on introduit l'espace
 $$\cC^{k-1}_{w,\mu}(X)=\{\alpha\in L_{w\mu}^2(\Lkm X), d\alpha \in L_\mu^2(\Lk X)\}$$
alors l'image de $d\,:\,\cC^{k-1}_{w,\mu}(X) \rightarrow L_\mu^2(\Lk X)$ est fermée et
son image est exactement $\overline{d\cD_\mu^{k-1}(d)}=\overline{ \left\{d\alpha,\alpha\in L_\mu^2(\Lkm X), d\alpha\in L_\mu^2 \right\}
}$.
\end{defi}
Cette notion de "presque fermeture" de l'image de $d$ est très pratique car elle permet d'obtenir des suites
de Mayer-Vietoris :

\begin{theo}Si $X=U\cup V$ ($U,V\subset X$ ouverts de $X$), on suppose que l'image
 de $d$ est presque fermée (en degré $k$ et par rapport à $w$) sur $X,U,V$ et $U\cap V$ et que 
 $$\{0\}\rightarrow\cC^{k-1}_{w,\mu}(X)
 \stackrel{\mbox{r}^*}{\longrightarrow} \cC^{k-1}_{w,\mu}(U)\oplus\cC^{k-1}_{w,\mu}(V)
\stackrel{\delta}{\longrightarrow}\cC^{k-1}_{w,\mu}(U\cap V)\rightarrow \{0\}$$ 
alors la suite exacte courte suivante est vérifiée:
 $$ \H^{k-1}_{w\mu}(U)\oplus \H^{k-1}_{w\mu}(V)
 \stackrel{\delta}{\longrightarrow} \H^{k-1}_{w\mu}(U\cap V)\stackrel{\mbox{b}}{\longrightarrow}
 \H^k_{\mu}(X) \stackrel{\mbox{r}^*}{\longrightarrow} 
 \H^k_{\mu}(U)\oplus \H^k_{\mu}(V)\stackrel{\delta}{\longrightarrow} \H^k_{\mu}(U\cap V).$$

\end{theo}
 Cette suite exacte permet en principe de déterminer $\H^k_{\mu}(X)$ lorsqu'on connaît 
 les autres espaces et les différentes flèches de cette suite exacte.
 
 Grâce à un résultat de L. Hörmander (\cite{Ho1}) et à nos précédents travaux sur la non parabolicité à l'infini, nous
 obtiendrons :
 \begin{theo}\label{Hormapara} Soit $(X,g)$ une variété riemannienne {\it complète} et 
 $w$ une fonction strictement positive lisse bornée sur $X$ (vérifiant \ref{nonpara}),  nous notons $d^*_w=w^{-1}d^*w$ l'adjoint formel de
 l'opérateur $d\,:\, L_w^2(\Lkm X)\rightarrow L_w^2(\Lk X)$ alors il existe un compact $K\subset X$ telle
 que
 $$\forall\varphi\in C^\infty_0(\Lkm (X\setminus K)),
  \|\varphi\|_{L^2_w}\le C\left[ \|d^*_w\varphi\|_{L^2}+\|d\varphi\|_{L^2}\right]$$
  si et seulement si 
  les images de $d\,:\, \cC^{k-2}_{w,w}(X)\rightarrow L^2_w(\Lkm X)$ et de 
  $d\,:\, \cC^{k-1}_{w,1}(X)\rightarrow L^2(\Lk X)$ sont fermées et si $\dim \H_w^{k-1}(X)<\infty$.
 \end{theo}
  
  Ceci est à rapprocher de la définition que nous avions introduite dans \cite{Carcre} et utilisée dans
  \cite{Cargafa} pour déterminer la cohomologie $L^2$ réduite des variétés plates au dehors d'un compact
  :
  
   \begin{defi} Soit $(X,g)$ une variété riemannienne {\it complète}, on dit que l'opérateur $d+d^*$ est
 non parabolique à l'infini s'il existe un compact $K\subset X$ tel que pour tout
 ouvert borné $U$ de $X\setminus K$, il existe une constante $C>0$ telle que
 $$\forall \varphi\in C^\infty_0(\Lambda^\bullet T^*(X\setminus K)),\ C\|\varphi\|_{L^2(U)}
 \le \|(d+d^*)\varphi\|.$$
 \end{defi}
 
 Décrivons maintenant l'organisation de ce papier : dans une première partie, on rappelle les définitions des
 espaces de cohomologie $L^2$  à poids, dans une deuxième partie nous décrivons notre critère de presque
 fermeture de l'image de $d$ et nous y donnons des conditions pour que l'image de $d$ soit presque fermée sur
 $U\cup V$ lorsqu'elle l'est sur $U,V$ et $U\cap V$ et pour qu'elle soit presque fermée sur 
 $U$ lorsqu'elle l'est sur $U\cup V$ et $U\cap V$. On compare aussi cette condition à la non parabolicité 
 à l'infini et nous montrerons le théorème \ref{Hormapara}. On calcule ensuite la cohomologie $L^2$ réduite
 des cônes et à titre d'illustration nous reprouvons un résultat de R. Melrose (\cite{Melrose}) à propos de
 la cohomologie $L^2$ des variétés à bouts coniques. On décrit ensuite la géométrie des variétés QALE.
  Enfin dans la dernière partie nous prouvons notre résultat
 principal (le théorème \ref{cohocrepante}), en fait notre méthode permet aussi de déterminer la cohomologie $L^2$ des variétés QALE
 asymptote à $\C^n/G$ (où $\bS^{2n-1}/G$ est à singularités isolées) et pas seulement des résolutions
 crépantes. La description de ces variétés faite plus haut implique qu'il y a un compact $K$ sur lequel
 la variété $X$ se rétracte. Et le bord de $K$ est une résolution de $\bS^{2n-1}/G$. Les singularités de 
 $\bS^{2n-1}/G$ sont de la forme $(\C^{m_i}/A_i\times \bS^{2n_i-1})/B_i$ où $A_i\subset \Su (m_i)$ agit sans points fixes sur
 $\C^{m_i}\setminus\{0\}$,
  $B_i$ agit sans point fixe sur $\C^{m_i}/A_i\setminus\{0\}$ et sur $\bS^{2n_i-1}$ ;
   ces singularités sont résolus par $\big(Y_i\times \bS^{2n_i-1}\big)/B_i$.
 Notons $Y=\cup_{i} Y_i/B_i$, on dispose d'une application $f\,:\, Y\rightarrow \partial K$, on peut
 considérer le sous complexe de $C^\infty(\Lk K)$ formé par les formes différentielles dont le tiré en
 arrière par $f$ est nulle et on note $H^*(K,\ker f^*)$ la cohomologie de ce complexe : un cas particulier
 de notre résultat est le suivant 
 \begin{theo} Soit $G\subset \Su (n)$ et 
 $(X,g)$ une variété QALE asymptote à $\C^n/G$ on suppose que 
 $\bS^{2n-1}/G$ est à singularités isolées et que la dimension réelle
  de ces singularités est plus grande ou égale
 à trois alors $$ \H^k(X)\simeq \left\{
 \begin{array}{ll}
 H^k(K,\partial K)\simeq H^k_c(X)& {\rm si\ } k\le 2\\

 H^k(K,\ker f^*)& {\rm si\ } k\in [3, 2n-2]\\
 H^k(K)\simeq H^k(X)& {\rm si\ } k\ge 2n-2\\
 \end{array}\right.$$
 \end{theo}
 
 \vskip0.5cm
 
\begin{merci} Je tiens à remercier C. Sorger pour avoir répondu patiemment à mes questions naïves sur la
cohomologie des résolutions crépantes. Je remercie Eugenie Hunsicker dont les commentaires sur
une version précédente de ce papier m'ont permis d'y
 repérer une erreur ; je la remercie également avec T. Hausel et R. Mazzeo et l' A.I.M.
pour avoir organisé et accueilli la conférence 
"$L^2$ Harmonic Forms in Geometry and Physics" au cours de laquelle cette question a été soulevée.
Enfin, ce travail a débuté alors que je bénéficiais d'une
délégation partielle au C.N.R.S.
\end{merci}
 
\section{$L^2$ cohomologie (à poids).}
Nous présentons ici les espaces de cohomologie $L^2$, on trouvera
d'autres présentations dans les articles de J. Cheeger, S. Zucker, J. Lott (\cite{Che,Z1,Lo} et de façon
plus abstraite dans le papier de J. Brüning et M. Lesch (\cite{BL}) .

Si $(X,g,\mu)$ est une variété riemannienne équipée d'une mesure lisse
$\mu(x)d\vol_g(x)$ (ou même
seulement dans $L^\infty_{\loc}$), on introduit l'espace $L^2_\mu(\Lk X)$ des
formes différentielles qui sont de carré sommables pour cette mesure, la norme
de cet espace de Hilbert est :
$$\|\alpha\|^2_\mu=\int_X |\alpha|^2 \mu d\vol.$$ 
On note $\cro{\ }{\  }_\mu$ le produit scalaire associé.

L'espace des $k$-formes $L^2_\mu$ fermées est l'espace $Z^k_\mu(X)$ des formes 
$\alpha \in L^2_\mu(\Lk X)$ qui sont fermées au sens des distributions, i.e.
celle qui vérifie :
$$\forall \beta \in C^\infty_0(\Lkp X),\ \cro{\alpha}{ d^*_\mu\beta}_\mu =0.$$
Où $d^*_\mu$ est l'adjoint différentiel formel de l'opérateur $$d\,:\, C^\infty_0(\Lk
X)\rightarrow L^2_\mu(\Lkp X) \ ;$$ si la variété $X$ est orientée et si $*$ est l'opérateur de
Hodge on a 
$$d^*_\mu=\pm \mu^{-1}*d*\mu$$
le signe étant fonction du degré.
Lorsque $\mu=1$, on ne mentionnera pas la référence à la fonction $\mu$ : on notera
$L^2_1=L^2$, $d^*=d_1^*$... 

Le domaine (maximal) de $d$ est 
$$\cD^k_\mu(d)=\{\alpha\in L^2_\mu(\Lk X),\  d\alpha \in L^2_\mu(\Lkp X)\}.$$
C'est à dire $\alpha\in \cD^k_\mu(d)$ si et seulement si 
$\alpha\in L^2_\mu(\Lk X)$ et s'il y a une constante $C$ tel que
$$\forall \varphi\in C^\infty_0(\Lkp X),\ |\cro{\alpha}{ d^*_\mu\varphi}_\mu|\le C\|\varphi\|_\mu.$$

Maintenant on fait l'hypothèse suivante :
\begin{equation}\begin{split}
 \text{ La variété } X \text{ est   l'intérieur d'une variété à coins } \bar X \text{ et
la métrique } g \\ \text{et le poids }\mu \text{ s'étendent à }\bar X \text{ et }\bar X \text{ est complet pour la distance
géodésique. }
\end{split}\label{HH}
\end{equation}
On note alors $C^\infty_0(\Lk X)$ l'espace des formes différentielles lisses et à support compacte
dans $X$ et $C^\infty_0(\Lk \bar X)$ l'espace des formes différentielles lisses et à support compacte
dans $\bar X$ c'est à dire celle dont le support est borné dans $(X,g)$. Avec cette hypothèse,
 $C^\infty_0(\Lk \bar X)$ est dense dans $\cD^k_\mu(d)$ lorsque
ce dernier espace est équipé de la norme
$$\alpha\mapsto \|\alpha\|_\mu+ \|d\alpha\|_\mu.$$
La $L^2_\mu$-cohomologie non réduite de $(X,g,\mu)$ est définie par
$${}^{nr} \H^k_\mu(X)=Z^k_\mu(X)/d\cD^{k-1}_\mu(d).$$
La $L^2$ cohomologie non réduite est la cohomologie d'un complexe mais lorsque l'image de $d$ n'est
pas fermée ce n'est pas un espace de Hilbert, 
c'est pourquoi on introduit la $L^2_\mu$-cohomologie réduite :
$$\H^k_{\mu}(X)=Z^k_\mu(X)/\overline{d\cD^{k-1}_\mu(d)}=Z^k_\mu(X)/\overline{dC^\infty_0(\Lkm \bar X) } .$$
On notera $B^k_\mu(X)=\overline{dC^\infty_0(\Lk \bar X) }$ où l'adhérence est pris dans
$L^2_\mu$. Évidemment lorsque l'image de $d\,:\,\cD^{k-1}_\mu(d) \rightarrow L^2_\mu(\Lk X)$ est fermée
ces deux espaces coïncident. L'objet de cet article est la cohomologie $L^2$ réduite, pour alléger 
le discours on omettra de signaler l'adjectif réduite.

La cohomologie $L^2_\mu$ se représente par des espaces de formes harmoniques.
On introduit le domaine de $d^*_\mu$, ou plus exactement de l'adjoint de 
$d\,:\, \cD^{k-1}_\mu(d)\rightarrow L^2_\mu(\Lk X)$:
lorsque $\alpha\in L^2_\mu(\Lk X)$ alors $\alpha\in \cD^k(d^*_\mu)$
si et seulement s'il y a une constante $C$ telle que
$$\forall \varphi\in C^\infty_0(\Lkm \bar X),\ |\cro{\alpha}{ d\varphi}_\mu|\le C\|\varphi\|_\mu.$$

Lorsque $\alpha\in C^\infty_0(\Lkm \bar X)$ alors $\alpha\in \cD(d^*_\mu)$ si et seulement si
le long de la partie lisse du bord de $\bar X$, $\alpha$ n'a pas de composante tangentielle.

\noindent Lorsque $\cD^k(d^*_\mu)$ est équipée de la norme du graphe
 $$\alpha\mapsto \|\alpha\|_\mu+ \|d_\mu^*\alpha\|_\mu$$
 alors $C^\infty_0(\Lk  X)$ est dense dans $\cD^k(d^*_\mu)$.
On introduit donc :

$$\cH^k_\mu(X)=\big\{\alpha\in Z^k_\mu(\Lk X),\ \alpha\in \cD(d^*_\mu)\text{ et }d^*_\mu\alpha=0\big\}.$$
La décomposition de Hodge-deRham nous apprend que :
\begin{eqnarray*}
&L^2_\mu(\Lk X)&=\cH^k_\mu(X)\oplus B^k_\mu(\Lk X)\oplus \overline{d^*_\mu\cD^{k+1}(d^*_\mu)} \\ 
& &=\cH^k_\mu(X)\oplus \overline{d C^\infty_0(\Lkm  \bar X)}\oplus \overline{d^*_\mu C^\infty_0(\Lkp  X)} \\ 
\end{eqnarray*}

\noindent l'adhérence étant bien entendu prise pour la topologie de $L^2_\mu(\Lk X)$ et on sait de plus que
$$Z^k_\mu(X)=\cH^k_\mu(X)\oplus B^k_\mu(\Lk X).$$ Donc on en déduit l'isomorphisme :
$$\H^k_{\mu}(X)\simeq \cH^k_\mu(X).$$

On peut aussi définir la $L^2_\mu$ cohomologie relative (à $\partial \bar X$ ou à 
$Y\subset \partial\bar X$ une hypersurface de $\bar X$)
 en considérant l'espace $\cD^k_\mu(d,Y)$ des formes $\alpha\in L^2_\mu(\Lk X)$
 pour lesquelles il existe une constante $C$ telle que
 $$ \forall \varphi \in C^\infty_0(\Lkp ( X\cup Y)),\ \,  
  |\cro{\alpha}{ d_\mu^*\varphi}_\mu|\le C\|\varphi\|_\mu . $$
  On définit alors 
  $$ Z^k_\mu(X,Y)=\big\{\alpha\in \cD^k_\mu(d,Y),\ d\alpha=0\,\big\}$$
  et
  $$B^k_\mu(X,Y)=\overline{ d\cD_\mu^{k-1}(d,Y)}.$$
 $$\H^k_\mu(X,Y)=Z^k_\mu(X,Y)/B^k_\mu(X,Y).$$
  Ces espaces admettent aussi une interprétation en terme
 de formes harmoniques :
 on introduit pour cela l'espace
 $\cD^k_\mu(d^*_\mu,Y)$ des formes $\alpha\in L^2_\mu(\Lk X)$ pour lesquelles il existe une constante
 $C$ telle que,
 $$ \forall \varphi \in C^\infty_0(\Lkp ( \bar X\cup Y)), \ 
  |\cro{\alpha}{ d\varphi}_\mu|\le C\|\varphi\|_\mu\  .$$
 Si on définit $$\cH^k(X,Y)=\{\alpha\in Z^k_\mu(X,Y)\cap\cD^k_\mu(d^*_\mu, Y),\  d^*_\mu\alpha=0 \}$$
  alors nous avons l'isomorphisme $$H^k_{\mu}(X,Y)\simeq \cH^k_\mu(X,Y).$$
  
\noindent Lorsque $X$ est orientée\footnote{Lorsque $X$ n'est pas orientée, un résultat analogue est valide en considérant
  les formes à valeurs dans le fibré d'orientation.}, l'opérateur $*$ de Hodge induit un isomorphisme :
   $$H^k_{\mu}(X)\simeq H^{n-k}_{\mu^{-1}}(X,\partial X).$$
\begin{rem}\label{cohol2loc}
On peut aussi introduire comme implicitement L. Hörmander \cite{Ho1, Ho2} 
et plus explicitement N. Teleman  \cite{Te}
la $L^2_{\loc}$-cohomologie. C'est la cohomologie du complexe des formes différentielles
 qui sont dans 
$L^2_\loc$ ainsi que leur différentielle. Selon la proposition 3.3 \cite{Te}, 
la cohomologie de ce complexe calcule la cohomologie de $X$
(la cohomologie de de Rham).
La dualité de Poincaré dit que la cohomologie de $X$ est isomorphe au dual de 
la cohomologie à support compact, en
particulier il n'y a pas lieu  de distingué la $L_{loc}^2$ cohomologie
 réduite ou non réduite car l'image de $d$ est fermée sur $L^2_\loc$.
  En effet, une forme fermée $L^2_{loc}$ nulle en cohomologie $L^2_\loc$ réduite est exacte
en restriction à chaque domaine compact à bord lisse de $X$, elle définit
 donc une forme linéaire nulle
 sur la cohomologie à support compact, elle est donc exacte en cohomologie $L^2_\loc$.
Aussi une forme fermée de
$L^2_\mu$ nulle en cohomologie $L^2_\mu$ réduite est nulle en cohomologie $L^2_\loc$ en particulier
elle est la différentielle d'une forme dans $L^2_\loc$.
De plus lorsque $\bar X$ est compacte alors la $L^2_\mu$ cohomologie (qui est évidemment
isomorphe à la $L^2_\loc$ cohomologie) est isomorphe à la cohomologie de $X$. 
\end{rem}

\section{la condition de presque-fermeture de l'image de $d$ et la suite de
Mayer-Vietoris}

On introduit maintenant une condition qui permet d'obtenir un petit bout de suite exacte de
Mayer-Vietoris.
\subsection{Définition}
 On considère donc $(X,g,\mu)$ une variété riemannienne équipée de la mesure $\mu(x)d\vol_g$ et 
 $w\,:\, X\rightarrow ]0,\infty[$ une fonction lisse bornée, on dira que 
 \begin{center}\it l'image de $d\,:\, \cD_\mu^{k-1}(d)\rightarrow L_\mu^2(\Lk X)$ est {\bf presque fermée}
 (en degré $k$ et par rapport à $w$) si lorsqu'on introduit l'espace
 $$\cC^{k-1}_{w,\mu}(X)=\{\alpha\in L_{w\mu}^2(\Lkm X), d\alpha \in L_\mu^2(\Lk X)\}$$
alors $B^k_\mu(X)=d\cC^{k-1}_{w,\mu}(X).$  \end{center}
Lorsque $(X,g,\mu)$ vérifie l'hypothèse (\ref{HH}) alors ceci
équivaut à ce que l'image de $d\,:\,\cC^{k-1}_{w,\mu}(X) \rightarrow L_\mu^2(\Lk X)$ soit fermée.
 Dans ces conditions une forme $L^2_\mu$ fermée $\alpha\in Z^k_\mu(X)$ est nulle en cohomologie
  $L^2_\mu$
 si et seulement si il existe $\beta\in L_{w\mu}^2(\Lkm X)$ telle que $\alpha=d\beta$
 \footnote{ au sens des distribution.}.
 Il est évident que lorsque l'image de $d$ est fermée elle est presque fermée.
 Remarquons que lorsque $w'\le w$, alors 
 $$\cC^{k-1}_{w,\mu}(X)\subset \cC^{k-1}_{w',\mu}(X)\ ;$$ ainsi si l'image de 
 $d\,:\,\cD^{k-1}_{\mu}(d) \rightarrow L_{\mu}^2(\Lk X)$ est presque fermée par rapport à $w$ alors on a 
 $$B^k_\mu(X)\subset d\cC^{k-1}_{w',\mu}(X),$$ mais on n'a pas forcément égalité.
  Dans (lemme 3.1 de \cite{Carpar}), 
 nous avons obtenu un critère qui assure l'inclusion inverse 
 (cf. aussi les articles de N. Hitchin, P. Li , J. Mc Neal, J. Jost et K. Zuo \cite{H,Li,Mc,JZ}):
 \begin{prop} 
 \label{poide}
 Si $(X,g,\mu)$ vérifie (\ref{HH}) et si $w(x)\ge\frac{1}{\psi^2(d(o,x))}$ 
 (où $o$ est un point
 fixé de $X$) où la fonction $\psi$ vérifie $$\int^\infty\frac{dt}{\psi(t)}=\infty$$ alors
 $$ d\cC^{k-1}_{w,\mu}(X)\subset B^k_\mu(X)$$
 \end{prop}
 
 \begin{defi}
 \label{nonpara}
 Lorsque $w$ vérifiera les hypothèses de cette proposition on dira que
 $w$ est à {\bf décroissance parabolique}.
 \end{defi}
 \noindent Un corollaire de ce résultat est donc le suivant :
 \begin{cor}
 \label{apoid} Si $(X,g,\mu)$ vérifie les hypothèses (\ref{HH}) et si l'image de 
 $d\,:\,\cC^{k-1}_{w,\mu}(X) \rightarrow L_{\mu}^2(\Lk X)$ est fermée alors pour tout poids $w'\le w$
 à décroissance parabolique, l'image de 
 $d\,:\,\cC^{k-1}_{w',\mu}(X) \rightarrow L_{\mu}^2(\Lk X)$ est
 encore fermée.
\end{cor}
 \subsection{Suite de Mayer-Vietoris.}
 \begin{thm}\label{MayerVietoris} Si $X=U\cup V$ ($U,V\subset X$ ouverts de $X$) et $w\,:\, X\rightarrow ]0,\infty[$ une fonction lisse,
  on suppose que l'image
 de $d$ est presque fermée (en degré $k$ et par rapport à $w$) sur $X,U,V$ et $U\cap V$ et que 
 $$\{0\}\rightarrow\cC^{k-1}_{w,\mu}(X)
 \stackrel{\mbox{r}^*}{\longrightarrow} \cC^{k-1}_{w,\mu}(U)\oplus\cC^{k-1}_{w,\mu}(V)
\stackrel{\delta}{\longrightarrow}\cC^{k-1}_{w,\mu}(U\cap V)\rightarrow \{0\}$$ 
alors la suite suivante est exacte :
 $$ \H^{k-1}_{w\mu}(U)\oplus \H^{k-1}_{w\mu}(V)
 \stackrel{\delta}{\longrightarrow} \H^{k-1}_{w\mu}(U\cap V)\stackrel{\mbox{b}}{\longrightarrow}
 \H^k_{\mu}(X) \stackrel{\mbox{r}^*}{\longrightarrow} 
 \H^k_{\mu}(U)\oplus \H^k_{\mu}(V)\stackrel{\delta}{\longrightarrow} \H^k_{\mu}(U\cap V).$$
 \end{thm}
 \proof
Si $A\subset B\subset X$ on note $\iota_{A,B}$ l'inclusion $\iota_{A,B}\,:\,
A\rightarrow B$ et on note $W=U\cap V$.

Montrons l'exactitude de la dernière flèche i.e $\ker \delta=\ima r^*$: 
soit donc $(\alpha,\beta)\in Z^k_{\mu}(U)\oplus Z^k_{\mu}(V)$ telle que 
 $$\iota_{W,U}^*\alpha-\iota_{W,V}^*\beta\in B^k_\mu(W).$$
  Il y a alors 
 $\psi\in \cC^{k-1}_{w,\mu}(W)$ tel que $\iota_{W,U}^*\alpha-\iota_{W,V}^*\beta=d\psi$ et  par
 hypothèse il y a aussi $(\psi_U,\psi_V)\in \cC^{k-1}_{w,\mu}(U)\oplus\cC^{k-1}_{w,\mu}(V)$ tel que
 $\iota_{W,U}^*\psi_U-\iota_{W,V}^*\psi_V=\psi$. Alors si $\gamma$ est la $k$-forme $L_\mu^2$ sur $X$
 telle que 
 $$\iota_{U,X}^*\gamma=\alpha-d\psi_U \et \iota_{V,X}^*\gamma=\beta-d\psi_V$$ alors $\gamma$ est
 fermée sur $X$ et on a bien $\iota_{U,X}^*\gamma-\alpha\in B^k_\mu(U)$ et
 $\iota_{V,X}^*\gamma-\beta\in B^k_\mu(V)$.
 
 On montre maintenant l'exactitude de la flèche du milieu i.e. $\ker r^*=\ima \mbox{b}$.
 Rappelons comment $\mbox{b}$ est défini, soit $$\eta\in Z^{k-1}_{w\mu}(U\cap V)\subset \cC^{k-1}_{w,\mu}(U\cap
 V),$$ 
 alors par hypothèse, on trouve $(\psi_U,\psi_V)\in \cC^{k-1}_{w,\mu}(U)\oplus\cC^{k-1}_{w,\mu}(V)$ tel que
 $$\eta=\iota^*_{W,U}\psi_U-\iota^*_{W,V}\psi_V$$
 alors la forme $\varphi$ définie sur $X$ par :
 $\iota^*_U\varphi=d\psi_U$ et $\iota^*_V\varphi=d\psi_U$ est un élément de $Z^k_\mu(X)$, sa classe
 de cohomologie $L_\mu^2$ ne dépend pas que de la classe de $L^2_{w\mu}$ cohomologie de $\eta$, et 
 alors $$\mbox{b}[\eta]=[\varphi].$$
 
 \noindent Soit donc $\gamma\in Z^k_\mu(X)$ telle que
 $$\iota_{U,X} ^*\gamma\in B^k_\mu(U)\et \iota_{V,X}^*\gamma\in B^k_\mu(V),$$
  on trouve donc 
 $(\psi_U,\psi_V)\in \cC^{k-1}_{w,\mu}(U)\oplus\cC^{k-1}_{w,\mu}(V)$ tel que :
 $\iota_{U,X}^*\gamma=d\psi_U$ et $\iota_{V,X}^*\gamma=d\psi_V$. Alors 
 $\iota_{W,U}^*\psi_U-\iota_{W,V}^*\psi_V\in Z^{k-1}_{w\mu}(W)$, et évidemment on a 
 $$\mbox{b}[\iota_{W,U}^*\psi_U-\iota_{W,V}^*\psi_V]=[\gamma].$$

 On montre finalement l'exactitude de la première flèche : c'est à dire $\ima \delta=\ker\mbox{b}$.
 Soit donc $\eta\in Z^{k-1}_{w\mu}(W)$, tel que si on choisit 
 $(\psi_U,\psi_V)\in  \cC^{k-1}_{w,\mu}(U)\oplus\cC^{k-1}_{w,\mu}(V)$ vérifiant
 $\eta=\iota_{W,U}^*\psi_U-\iota_{W,V}^*\psi_V$ alors la $k$-forme $L^2$ fermée $\gamma$
 définie par
 $$\iota_{U,X} ^*\gamma=d\psi_U \et \iota_{V,X} ^*\gamma=d\psi_V$$ est nulle en cohomologie
 $L^2_\mu$. Par hypothèse il existe $\phi\in \cC^{k-1}_{w,\mu}(X)$ telle que
 $\gamma=d\phi$. Mais alors on a
 $\alpha=\psi_U-\iota^*_{U,X}\phi\in Z^{k-1}_{w\mu}(U)$ et 
 $\beta=\psi_V-\iota^*_{V,X}\phi\in Z^{k-1}_{w\mu}(V)$
 d'où $\eta=\iota_{W,U}^*\alpha-\iota_{W,V}^*\beta$.
 
 \endproof
 \begin{rems}
 \begin{enumeroman}\label{rmMayerVietoris}
 \item On peut évidemment énoncer un théorème analogue pour la cohomologie $L^2_\mu$ relative à une
 hypersurface $Y\subset \partial X$ à condition d'introduire l'espace  
 $\cC^{k-1}_{w,\mu}(X,Y)$ des formes $\alpha\in L^2_{w\mu}(\Lkm X)$ telles qu'il existe une constante $C$ avec
 $$\forall \varphi\in
 C^\infty_0(\Lk (X\cup Y) )\ |\cro{\alpha}{d^*_\mu\varphi}_\mu|\le C \|\varphi\|_\mu.$$
 \item La preuve indique aussi que pour avoir l'exactitude au niveau des deux dernières flèches, on a uniquement
 besoin que l'image de $d$ soit presque fermée sur $U\cap V$, pour l'exactitude au niveau des deux flèches du
 milieu, il suffit que l'image de $d$ soit presque fermée sur $U$ et $V$ et au niveau des deux premières flèche
 il suffit que l'image de $d$ soit presque fermée sur $X$.
 \end{enumeroman}
 \end{rems}

 \subsection{Condition pour assurer la presque fermeture de l'image de $d$.}
Si on découpe la variété en petits bouts sur lesquels l'image de $d$ est presque fermée, il n'est pas
clair que l'image de $d$ soit presque fermée sur la variété ; voici une proposition qui assure ceci
 (c'est en quelque sorte un voyage de la cohomologie vers l'analyse fonctionnelle).
\begin{prop}\label{UVX}
Soit $(X,g)$ une variété riemannienne vérifiant (\ref{HH}) et $\mu,w,w'$ des fonctions strictement positives lisses sur $\bar X$. On suppose que
$X=U \cup V$  ($U,V\subset X$ ouverts de $X$) et que les conditions suivantes sont vérifiées:
 \begin{itemize}
\item l'image de $d$ est presque fermée sur $U$ et $V$ par rapport au poids $w$.
\item Le poids $w$ est à décroissance parabolique.
\item $d  \cC^{k-2}_{w',w\mu}(U\cap V)=B^{k-1}_{w\mu}(U\cap V)$.
\item $\cC^{k-2}_{w',w\mu}(U)\oplus\cC^{k-2}_{w',w\mu}(V)
\stackrel{\delta}{\longrightarrow}\cC^{k-2}_{w',w\mu}(U\cap V)\rightarrow \{0\}$
\item la suite suivante $\H^{k-1}_{w\mu}(U)\oplus \H^{k-1}_{w\mu}(V)\rightarrow \H^{k-1}_{w\mu}(U\cap V)
\stackrel{\mbox{b}}{\longrightarrow}
  \H^k_{\mu}(X)$ est aussi exacte,
\end{itemize}
alors l'image de $d$ est presque fermée sur $X$ par rapport au poids $w$.
\end{prop}
\proof Soit donc $\gamma\in B^k_\mu(X)$, évidemment on a aussi 
$$  \iota_{U,X} ^*\gamma\in B^k_\mu(U)\et
 \iota_{V,X}^*\gamma\in B^k_\mu(V) $$
 Il y a donc $(\psi_U,\psi_V)\in \cC^{k-1}_{w,\mu}(U)\oplus\cC^{k-1}_{w,\mu}(V)$ tels que
 $$  \iota_{U,X} ^*\gamma=d\psi_U \et
 \iota_{V,X}^*\gamma =d\psi_V.$$
On définit alors $\eta=\iota_{W,U}^*\psi_U-\iota_{W,V}^*\psi_V\in Z^{k-1}_{w\mu}(W)$ et cette forme vérifie
par construction $b[\eta]=0$. Il y a donc $\phi\in \cC^{k-2}_{w',w\mu}(W)$ et $(\alpha_U,\alpha_V)\in
Z^{k-1}_{w\mu}(U)\oplus Z^{k-1}_{w\mu}(V)$ tels que
$$\eta=d\phi+\iota^*_{W,U}\alpha_U-\iota^*_{W,V}\alpha_V.$$ 
Et on trouve aussi $(\phi_U,\phi_V)\in \cC^{k-2}_{w',w\mu}(U)\oplus\cC^{k-2}_{w',w\mu}(V)$
tels que 
$$\iota_{W,U}^*\phi_U-\iota_{W,V}^*\phi_V=\phi.$$
Si $\psi\in \cC^{k-1}_{w,\mu}(X)$ est défini par
$$\iota_{U,X} ^*\psi=\psi_U-d\phi_U-\alpha_U \et
 \iota_{V,X}^*\psi =\psi_V-d\phi_V-\alpha_V, $$ alors on a bien 
 $\gamma=d\psi$.
Ainsi, nous avons obtenu l'inclusion : $$B^k_\mu(X)\subset d\cC^{k-1}_{w,\mu}(X).$$
  L'inclusion réciproque est une conséquence du fait que $w$ est à décroissance
  parabolique.
 \endproof
  
  Nous avons bien évidemment un énoncé analogue pour la cohomologie $L^2_\mu$ relative à
  $Y\subset \partial X$.
  
  On a aussi un critère qui assure que l'image de $d$ est presque fermée sur $U$ lorsqu'elle l'est
  sur $X$ et $U\cap V$.
  \begin{prop} \label{XU} Soit $(X,g)$ une variété riemannienne telle que (\ref{HH}) soit vérifiée
à coins et $\mu,w$ des fonctions strictement positives lisses sur $\bar X$. On suppose que
$X=U \cup V$ où $U$ vérifie aussi (\ref{HH}). Les conditions suivantes
 \begin{itemize}
\item L'image de $d$ est presque fermée sur $X$ et $U\cap V$ par rapport au poids $w$.
\item Le poids $w$ est à décroissance parabolique.
\item $\{0\}\rightarrow\cC^{k-1}_{w,\mu}(X)
 \stackrel{\mbox{r}^*}{\longrightarrow} \cC^{k-1}_{w,\mu}(U)\oplus\cC^{k-1}_{w,\mu}(V)
\stackrel{\delta}{\longrightarrow}\cC^{k-1}_{w,\mu}(U\cap V)\rightarrow \{0\}$
\item $\iota_{U,X}^* \,:\, \H^k_\mu(X)\rightarrow \H^k_\mu(U)$ est injective
\end{itemize}
impliquent que l'image de $d$ est presque fermée sur $U$ par rapport au poids $w$.
  \end{prop}
  \proof Soit $\alpha\in B^k_\mu(U)$ on a donc $\iota_{W,U}^*\alpha\in B^k_\mu(W)$
  et on trouve de même   
  $(\psi_U,\psi_V)\in \cC^{k-1}_{w,\mu}(U)\oplus\cC^{k-1}_{w\mu}(V)$ tels que
  $d\iota_{W,U}^*\psi_U-d\iota_{W,V}^*\psi_V=\iota_{W,U}^*\alpha$. Maintenant
  si $\gamma\in Z^k_\mu(X)$ est défini par 
  $$\iota_{U,X} ^*\gamma=\alpha-d\psi_U \et
 \iota_{V,X}^*\gamma =-d\psi_V.$$
  Alors $\gamma$ est fermé et $\iota_{U,X}^*\gamma\in B^k_\mu(U)$ donc par hypothèse :
  $\gamma\in B^k_\mu(X)$ et on trouve $\phi\in \cC^{k-1}_{w,\mu}(X)$ tel que
  $\gamma=d\phi$, ainsi   $\alpha=d\psi_U+d\iota_{U,X}^*\phi$. D'où l'inclusion
  $$B^k_\mu(U)\subset d\cC^{k-1}_{w,\mu}(U).$$
  L'inclusion réciproque est une conséquence du fait que $w$ est à décroissance parabolique.
  \endproof
  
  \begin{rem}
  \label{L2loc}
 On remarque que l'on peut remplacer 
  l'hypothèse que l'image de $d$ est presque fermée sur $X$
  par  l'hypothèse que si $\alpha\in B^k_\mu(X)$ alors il y a une forme $\beta\in L^2_\loc$
  tel que $\alpha=d\beta$ et $\iota_{U,X}^*\beta\in L^2_{w\mu}$.
   \end{rem} 
  La preuve de la proposition précédente montre même le résultat suivant :
   \begin{cor}\label{XUiso}
   Si en plus des hypothèses précédentes, on a $\H^k_\mu(U\cap V)=\{0\}$, alors 
   $$\H^k_\mu(U)\simeq  \H^k_\mu(X).$$ 
   \end{cor}
 \subsection{Comparaison avec la non parabolicité}
 Dans \cite{Cargafa}, nous avons utilisé une condition similaire à la presque fermeture de l'image de $d$:
 \begin{defi} Soit $(M,g)$ une variété riemannienne {\it complète}, on dit que l'opérateur $d+d^*_\mu$ est
 non parabolique à l'infini s'il existe un compact $K\subset M$ tel que pour tout
 ouvert borné $U$ de $M\setminus K$, il existe une constante $C>0$ telle que
 $$\forall \varphi\in C^\infty_0(\Lambda^\bullet T^*(M\setminus K)),\ C\|\varphi\|_{L^2(U)}
 \le \|(d+d^*_\mu)\varphi\|_\mu.$$
 \end{defi}
 Soit $\tilde K$ un compact contenant $K$ dans son intérieur, on note 
 $W(\Lambda^\bullet T^*M)$ le complété de $C^\infty_0(\Lambda^\bullet T^*M)$ pour la norme :
 $$\alpha\mapsto \|\alpha\|_{L^2(\tilde K)}+\|(d+d^*_\mu)\varphi\|_\mu.$$ 
 Cet espace ne dépend pas du compact $\tilde K$ choisi, il s'injecte donc naturellement dans 
 $W^{1,2}_\loc$. De plus l'opérateur $(d+d^*_\mu)\,:\,W(\Lambda^\bullet T^*M) \longrightarrow
 L^2_\mu(\Lambda^\bullet T^*M)$ est Fredholm et on a
 $$B^k_\mu(M)=dW(\Lkm M).$$
 Si l'on fait l'hypothèse de coercivité suivante :
 $$\forall \varphi\in C^\infty_0(\Lambda^\bullet T^*(M\setminus K)),\ C\|\varphi\|_{w\mu}
 \le \|(d+d^*_\mu)\varphi\|_{\mu}\, ,$$
 alors l'opérateur $d+d^*_\mu$ est clairement non parabolique à l'infini et l'image de $d$ est
 presque fermée sur $M$ pour le poids $w$. Dans ce cas une norme sur l'espace
   $W(\Lambda^\bullet T^*M)$ est  
   $$\alpha\mapsto \|\varphi\|_{w\mu}+\|(d+d^*_\mu)\varphi\|_\mu.$$
Alors le noyau $W$ de l'opérateur $d+d_\mu^*$ est le noyau $L^2_{w\mu}$ de l'opérateur
   $d+d^*_\mu$ ; et il y a donc une constante $C>0$ telle que si 
   $\varphi\in C^\infty_0(\Lambda^\bullet T^*M)$ vérifie
   $$\cro{\varphi}{h}_{w\mu}=0,\  \forall h\in \ker_{L^2_{w\mu}}(d+d^*_\mu)$$ alors
  $$C \|\varphi\|_{w\mu}\le\|(d+d^*_\mu)\varphi\|_\mu$$
  
  L'avantage de cette définition est qu'elle ne dépend que de la géométrie à l'infini de $M$ et 
  qu'avec un peu d'analyse on  peut obtenir une suite exacte longue en cohomologie $L^2_\mu$·
  Mais ceci ne permet pas de découper la géométrie à l'infini en morceaux sur lesquels on ferait une
  analyse similaire. La condition de presque fermeture pour l'image de $d$ est assez souple mais
  nous ne savons pas si elle ne dépend que de la géométrie à l'infini :
  par exemple si $(X,g)$ est une variété riemannienne qui vérifie nos hypothèses (\ref{HH}) et $\mu$ une fonction lisse
  strictement positive sur $\bar X$. Soit $K\subset \bar X$ un compact, est-il vrai que la presque
  fermeture de l'image de $d$ sur $X$ implique la presque
  fermeture de l'image de $d$ sur $X\setminus K$ (et vice versa) ? Notons par exemple
  que l'hypothèse de presque fermeture de l'image de $d$ ne requiert pas la finitude de la
  dimension des espaces de $L^2_\mu$ cohomologie.
  
  \begin{lem} Supposons que $(X,g,\mu)$ vérifie la condition (\ref{HH})
   et que $w$ soit à décroissance parabolique, alors 
   $C^\infty_0(\Lkm \bar X)$ est dense dans $\cC^{k-1}_{w,\mu}(X)$ lorsque ce dernier espace est équipé
    de la norme du graphe $ \alpha\mapsto\|\alpha\|_{w\mu}^2+\|d\alpha\|_{\mu}.$
  \end{lem}
  \proof   C'est en fait une re-formulation du lemme (3.1) de \cite{Carpar}. Grâce aux hypothèses (\ref{HH}),
  il suffit de démontrer que si $\beta\in \cC^{k-1}_{w,\mu}(X)$, alors on peut trouver
  une suite $(\beta_l)_l$ d'élément $\cC^{k-1}_{w,\mu}(X)$ à support compact qui converge vers $\beta$.
  
  Rappelons que $w$ est à décroissance parabolique : $$w(x)\ge\frac{1}{\psi^2(d(o,x))}$$
    où $o$ est un point
 fixé de $X$ et $\psi$ vérifie $$\int^\infty\frac{dt}{\psi(t)}=\infty.$$
 On introduit alors les fonctions de troncature de 
   \cite{Carpar} :
  si $r,R$ deux nombres réels tels que $0<r<R$, on leurs associe la fonction $\phi_{r,R}$ définie
par
\begin{equation}
\label{tronc}
\phi_{r,R}(s)=\left\{\begin{array}{lll}
1& \si & s\le r\\
\int_s^R \frac{dt}{\psi(t)} \times  \left(\int_r^R \frac{dt}{\psi(t)}\right)^{-1} &\si &s\in[r,R]\\
0& \si &s\ge R\\
\end{array}
\right.\end{equation}
Il existe alors par hypothèse des suites $r_l<R_l$ telles que $\lim_{l\to\infty } r_l=\infty$ et 
 $$\left(\int_{r_l}^{R_l} \frac{dt}{\psi(t)}\right)^{-1}  \le \frac{1}{l} \ ;$$  
 
 On considère alors  $\beta_l=\phi_{r_l,R_l}(r(x)) \beta$ où $r(x)=d(o,x)$.
 Il est clair que $\lim_{l\to\infty} \|\beta_l-\beta\|_{w\mu}=0.$
  De plus puisque $d\beta_l=\phi_{r_l,R_l} d\beta+\phi_{r_l,R_l}'dr\wedge\beta$ et que
  par construction :
  $$\|\phi_{r_l,R_l}'dr\wedge\beta\|_\mu\le \frac2l \|\beta\|_{w\mu}\ , $$
  on déduit immédiatement que $\lim_{l\to\infty} \|d\beta_l-d\beta\|_{\mu}=0.$
 \endproof
 Introduisons maintenant l'opérateur
  $$T:=\left( d\,:\, \cC_{w,w\mu}^{k-2}(X)\rightarrow L^2_{w\mu}(\Lkm X)\right)$$ et
 $T^*$ son adjoint, le domaine de $T^*$
  est aussi l'espace $\cD(T^*)$ des $(k-1)$-formes $\alpha\in L^2_{w\mu}$ telles que
  $$\forall\phi\in C^\infty_0(\bar X),\ \ \phi\alpha\in  \cD^{k-1}(d^*_{w\mu})$$ et telles que
  $d^*_{w\mu}\alpha\in L^2_\mu$. En utilisant que 
 $d^*_{w\mu}\beta_l=\phi_{r_l,R_l}d^*_{w\mu}\beta-\phi_{r_l,R_l}'\inte_{\nabla r}  \beta,$
  une preuve identique  montre que
 l'on a aussi :
 \begin{lem}\label{densite}  Supposons que $(X,g,\mu)$ vérifie la condition (\ref{HH})
   et que $w$ soit à décroissance parabolique, alors 
   $C^\infty_0(\Lkm \bar X)\cap\cD^{k-1}(d^*_{w\mu})$ est dense dans $\cC^{k-1}_{w,\mu}(X)\cap \cD(T^*)$ lorsque ce dernier espace est équipé
    de la norme du graphe $ \alpha\mapsto\|\alpha\|_{w\mu}^2+\|d\alpha\|_{\mu}+\|d^*_{w\mu}\alpha\|_{\mu}.$
  \end{lem}
Grâce à ce dernier lemme, un théorème de L.
Hörmander (théorème  1.1.2 de \cite{Ho1}) implique  néanmoins le résultat suivant qui établit un premier 
lien entre les notions de non-parabolicité à l'infini et de presque-fermeture de l'image de $d$. 
  
  \begin{thm} Soit $(X,g)$ est une variété riemannienne orientée de dimension $d$ qui vérifie nos hypothèses
  (\ref{HH})et $\mu,w $ des fonctions lisses
  strictement positives sur $\bar X$ où $w$ est à décroissance parabolique. Les images des applications
  $$d\,:\, \cC^{k-1}_{w,\mu}(X)\rightarrow L^2_\mu(\Lk X)$$ 
  $$\et d\,:\, \cC^{k-2}_{w,w\mu}(X)\rightarrow
  L^2_{w\mu}(\Lambda^{k-1}T^* X) $$
  sont fermées si et seulement si il y a une constante $C>0$ telle que si $\varphi\in C^\infty_0(\Lkm \bar X)\cap
  \cD^{k-1}(d^*_{w\mu})$ vérifie 
  $$\cro{\varphi}{h}_{w\mu}=0,\ \ \forall h\in \cH^{k-1}_{w\mu}(X)$$
  alors 
  $$C \|\varphi\|^2_{w\mu}\le \|d\varphi\|^2_{\mu}+\|d^*_{w\mu}\varphi\|^2_{\mu}$$
  \end{thm}

  On peut aussi re-visiter une partie de nos travaux sur la non-parabolicité à l'infini à partir 
  de cette
  inégalité. Nous remarquons que si en plus des hypothèses de ce théorème nous supposons que
  $\cH^{k-1}_{w\mu}(X)$ est de dimension finie alors on trouve $K\subset \bar X$ un compact et
   $C$ une
  constante strictement positive telle que
  \begin{equation*}
  \begin{gathered}
  \forall \phi\in C^\infty_0(\Lkm (\bar X\setminus K))\cap \cD^{k-1}(d^*_{w\mu}),\\ 
  C\|\phi\|^2_{w\mu }\le \|d\phi\|^2_{\mu }+\|d^*_{w\mu}\phi\|^2_{\mu }.
  \end{gathered}\end{equation*}
  
  En effet si $\{h_1,...,h_l\}$ est une base orthonormée de $\cH^{k-1}_{w\mu}(X)$, nous avons
  pour $\phi\in C^\infty_0(\Lkm (\bar X\setminus K))\cap \cD^{k-1}(d^*_{w\mu})$
 $$ C\|\phi\|^2_{w\mu }\le \|d\phi\|^2_{\mu }+\|d^*_{w\mu}\phi\|^2_{\mu }+C\sum_i
 |\cro{\phi}{h_i}_{w\mu}|^2$$
 Il suffit de choisir alors $K$ afin que 
 $$\sum_i \|h_i\|^2_{L^2_{w\mu}(M\setminus K)}\le \frac{1}{10}\ .$$
 Ce résultat s'accompagne en fait d'un réciproque :
  
  \begin{thm} Supposons en plus que $w$ soit à décroissance parabolique et que
  pour un compact $K\subset \bar X$, il existe une constante $C>0$ telle que 
 \begin{equation}\label{npy}
  \begin{gathered}\forall \phi\in C^\infty_0(\Lkm (\bar X\setminus K))\cap \cD^{k-1}(d^*_{w\mu})\\
  C\|\phi\|^2_{w\mu }\le \|d\phi\|^2_{\mu }+\|d^*_{w\mu}\phi\|^2_{\mu }.
  \end{gathered}\end{equation}
  
\noindent alors l'image de $d$ est presque fermée en degré $k$ pour le poids $w$, 
   $\cH^{k-1}_{w\mu}(X)$ est de dimension finie et l'image de $d\,:\,\cC^{k-2}_{w,w\mu}(X)\rightarrow
   L^2_{w\mu}(\Lkm X)$ est aussi fermée.
  \end{thm}
 
 \proof On s'inspire largement de \cite{Carcre} : soit $\tilde K$ un compact contenant $K$ dans son
 intérieur, on lui associe la norme $N_{\tilde K}$ définie sur
  $C^\infty_0(\Lkm \bar X)\cap \cD^{k-1}(d^*_{w\mu})$ par
 $$N_{\tilde K}(\phi)=\left(\int_{\tilde K} |\phi|^2+\|d\phi\|^2_{\mu }+\|d^*_{w\mu}\phi\|^2_{\mu }\right)^{1/2}\ .$$
 Cette norme est en fait équivalente à la norme :
 $$N_{w}(\phi)= \left(\|\phi\|^2_{w\mu}+\|d\phi\|^2_{\mu }+\|d^*_{w\mu}\phi\|^2_{\mu }\right)^{1/2}.$$

 En effet la norme $N_{w}$ est à un facteur près plus grande que la norme $N_{\tilde K}$, il suffit donc
 de montrer l'estimation réciproque. Soit donc $\rho$ une fonction qui vaut $1$ dans un voisinage de
  $X\setminus\tilde K$ et à support dans $X\setminus K$. Si
  $\phi \in C^\infty_0(\Lkm \bar X)\cap \cD^{k-1}(d^*_{w\mu})$, on a 
  \begin{equation*}
  \begin{split}
  \|\phi\|^2_{w\mu}&\le 2\|(1-\rho)\phi\|^2_{w\mu}+2\|\rho\phi\|^2_{w\mu}\\
  &\le C \int_{\tilde K} |\phi|^2+C \left[\|d(\rho\phi)\|^2_{\mu }+
  \|d^*_{w\mu}(\rho\phi)\|^2_{\mu }\right]\\
   &\le C \int_{\tilde K} |\phi|^2+C\left[\|d\phi\|_{\mu }^2+\|d^*_{w\mu}\phi\|_{\mu}^2
   \right] \\
 \end{split} \end{equation*}
 Où on a utilisé d'abord l'inégalité de non parabolicité (\ref{npy}) et ensuite que 
  $$ \|d(\rho\phi)\|^2_{\mu }+\|d^*_{w\mu}(\rho\phi)\|^2_{\mu }
   \le C \int_{\tilde K} |\phi|^2+C\left[\|d\phi\|^2_{\mu }+\|d^*_{w\mu}\phi\|^2_{\mu
   }\right]$$
   Soit donc $W$ le complété de $C^\infty_0(\Lkm \bar X)\cap \cD^{k-1}(d^*_{w\mu})$
   pour la topologie induite par ces normes. Nous allons montrer que l'opérateur 
   $$(d+d^*_{w\mu})\,:\, W\longrightarrow L^2_\mu(\Lk X)\oplus L^2_\mu(\Lkm X)$$
   est semi-Fredholm : son noyau est de dimension finie et son image est fermée.

   Grâce au critère très commode de L. Hörmander (prop. 19.1.3 de \cite{Ho3}), il suffit pour cela de montrer que
   si $\{\alpha_l\}_l$ converge faiblement dans $W$ et que $\{(d+d^*_{w\mu})\alpha_l\}_l$
   converge fortement dans $L^2_\mu$ alors la suite $\{\alpha_l\}_l$ converge fortement dans $W$.
   
  Il suffit donc de démontrer qu'une telle suite converge fortement dans $L^2(\Lkm\tilde K)$.
   Mais par hypothèse une telle suite est bornée dans $W^{1,2}_\loc$ et
   converge faiblement dans $W^{1,2}_\loc$, puisque l'inclusion de $W^{1,2}_\loc$ dans $L^2(\tilde K)$
   est compacte, le résultat est immédiat.
On vient de démontrer que l'image de $d\,:\, W\longrightarrow L^2_\mu$ est fermée.
Grâce au lemme \ref{densite}, on en déduit que les opérateurs 
$d\,:\, \cC^{k-1}_{w,\mu}(X)\rightarrow L^2_\mu $ et $T^*$ ont leurs images fermées, c'est donc aussi le
cas de $T$.

 On montre maintenant que 
   $\cH^{k-1}_{w\mu}(X)$ est de dimension finie, on sait que $ker_W
   (d+d^*_{w\mu})\subset\cH^{k-1}_{w\mu}(X)$. On va vérifier que ces deux espaces coïncident.
   Soit $h\in \cH^{k-1}_{w\mu}(X)$, on veut montrer que $h$ est dans $W$, 
   il faut donc trouver une suite $(h_l)_l$  de 
   $C^\infty_0(\Lkm \bar X)\cap \cD^{k-1}(d^*_{w\mu})$ 
   telle que $\|dh_l\|_\mu\to 0$ et $\|d^*_{w\mu}h_l\|_\mu\to 0$.
  On reprend les fonction de troncature (\ref{tronc}), et 
  on considère alors  $h_l=\phi_{r_l,R_l}(r(x)) h$ où $r(x)=d(o,x)$.
 On a $dh_l=\phi_{r_l,R_l}'dr\wedge h$ et $d^*_{w\mu}h_l=-\phi_{r_l,R_l}'\inte_{\nabla r} h$
 on obtient donc $$\|dh_l\|_\mu+\|d^*_{w\mu}h_l\|_\mu\le \frac2l \|h\|_{w\mu}$$
 d'où le résultat.
 \endproof
 
 La conclusion de cette discussion est donc la suivante :
 
 \begin{thm}
 \label{conclusion}
  Lorsque $(X,g,\mu)$ vérifie (\ref{HH}) et que $w$ est à décroissante parabolique,
 la conjonction des trois conditions suivantes ne dépend que de la géométrie à l'infini de $(X,g,\mu)$ :
 \begin{enumerate}
 \item La dimension de $\H^{k-1}_{w\mu}(X)$ est finie.
 \item L'image de $d\,:\, \cC^{k-1}_{w,\mu}(X)\rightarrow L^2_\mu(\Lk X)$ est fermée.
 \item L'image de $d\,:\, \cC^{k-2}_{w,w\mu}(X,\partial X)
 \rightarrow L^2_{w\mu}(\Lambda^{k-1} T^* X)$ est fermée.
 \end{enumerate}
  Et elles équivalent à l'existence d'un compact $K\subset \bar X$ tel que :
  \begin{equation*}
  \begin{gathered}\forall \phi\in C^\infty_0(\Lkm (\bar X\setminus K))\cap \cD^{k-1}(d^*_{w\mu})\\
  C\|\phi\|^2_{w\mu }\le \|d\phi\|^2_{\mu }+\|d^*_{w\mu}\phi\|^2_{\mu }.
  \end{gathered}\end{equation*}
 \end{thm}
On peut également énoncer un résultat analogue relativement à $Y\subset \partial X $.

 \subsection{}
 On donne maintenant le critère qui plus loin nous permettra d'affirmer que la suite 
$$\{0\}\rightarrow\cC^{k-1}_{w,\mu}(X)
 \stackrel{\mbox{r}^*}{\longrightarrow} \cC^{k-1}_{w,\mu}(U)\oplus\cC^{k-1}_{w,\mu}(V)
\stackrel{\delta}{\longrightarrow}\cC^{k-1}_{w,\mu}(U\cap V)\rightarrow \{0\}$$  est exacte. 
\begin{prop}\label{condition} Si $X=U\cup V$ où $\partial V\cap U\subset \Omega$ où
$\Omega$ est quasi isométrique au bout de cône sur $]-1,1[\times \Sigma$: $\Omega\simeq C_1(]-1,1[\times
\Sigma)$. On suppose que $\partial V\cap U=C_1(\{0\}\times \Sigma)$ et $V\cap \Omega=C_1(]-1,0[\times \Sigma)$.
On considère un poids $\mu$  qui est radiale sur $\Omega$ et des poids $w, \bar w, \rho $ tels que si on note
$r$ la variable radiale sur $\Omega$ alors en restriction à $\Omega$ : $w=r^{-2} (1+\log(r))^{-2}$,
 $\rho=(1+\log(r))^{-2}$ et  $\bar w=r^{-2}$. Alors on a
$$\cC^{k-1}_{\bar w,\mu}(U)
\stackrel{\iota^*_{U\cap V,U}}{\longrightarrow}\cC^{k-1}_{\bar w,\mu}(U\cap V)\rightarrow \{0\}$$ et
également  pour tout $\psi\in \cC^{k-1}_{ w,\mu}(U\cap V)$, il existe $\phi\in \cC^{k-1}_{ w,\rho\mu}(U)$
tel que $$\iota^*_{U\cap V,U}\phi=\psi.$$
\end{prop}
\proof Les arguments de la preuve sont en fait empruntés à J. Cheeger (lemme 4.2 de \cite{Che}): on considère
$s\,:\, \Omega\rightarrow ]-1,1[$ la projection sur le deuxième facteur
$\Omega=]1,\infty[\times]-1,1[\times \Sigma$, on note $\Omega^\pm=\{\pm s >0\}$ et $\sigma$ l'involution de
$\Omega$ qui échange $\Omega^\pm$ avec $\Omega^\mp$.
Soit $\xi$ une fonction lisse sur $\Omega$ qui ne dépend que de $s$ tel que
$\xi(\Omega)\subset [0,1]$ et 
$\xi=0$ si $s\ge 1/2$, et $\xi=1$ si $s\le 1/3$. On étend cette fonction à $U\cup \Omega$
en lui imposant la valeur $0$ sur le complémentaire de $\Omega$.
Lorsque $\alpha\in \cC^{k-1}_{w,\mu}(U\cap V)$, on étend $\alpha$ à $(U\cap V) \cup \Omega$ en reflétant
$\alpha$ le long de $\partial V\cap U=C_1(\{0\}\times \Sigma)$, i.e. sur $\Omega\setminus V$, 
$\alpha=\sigma^*\alpha$ et on pose $\bar\alpha= \xi\alpha$.
On vérifie facilement que si $\nu$ est une fonction positive sur $X$ qui est radiale sur $\Omega$ alors
$$\|\bar \alpha\|_{L^2_\nu(\Omega)}\le 2 \| \alpha\|_{L^2_\nu(\Omega^+)}.$$
Les arguments de J. Cheeger montre que si $\alpha\in\cC^{k-1}_{w,\mu}(U\cap V)$   alors
$$d\bar\alpha=\overline{d\alpha}+d\xi\wedge\bar\alpha$$ 
puisque le gradient de $\xi$ est borné sur $\Omega$ par $1/r$, la proposition est bien démontrée.

\endproof
\section{$L^2$ cohomologie des cônes.}
\subsection{Inégalités de Hardy}
On rappelle ici quelques inégalités de type Hardy qui seront utiles un peu plus loin à propos de la $L^2$
cohomologie des cônes et des variétés QALE. Pour plus de détails sur ces inégalités nous renvoyons le lecteur au
traité (\cite{hardy}) et dans un cadre un peu plus général nous conseillons le survol de E.B. Davies
(\cite{Da}).

\begin{lem} \label{hardynp}
Soit $\rho\, :\, [1,\infty[\rightarrow \R_+^*$ une fonction positive telle que
$$\int^\infty \frac{dt}{\rho(t)}<\infty.$$ Si on définit  
$$G(t)=\int_t^\infty \frac{d\tau}{\rho(\tau)},$$ alors pour tout $\varphi\in W^{1,2}_{loc}([1,\infty[)$ tel que
$\varphi'\in L^2([1,\infty[, \rho dt)$ et $\lim_{t\to \infty}\varphi(t)=0$  on a
$$\frac14\int_1^\infty \left(\frac{G'(t)}{G(t)}\right)^2 |\varphi(t)|^2\rho(t)dt\le 
\int_1^\infty  |\varphi'(t)|^2\rho(t)dt.$$
\end{lem}
\proof La preuve est standard on pose $\varphi(t)=\sqrt{G(t)-G(T)\,}\, u(t)$ et on a alors
$$\int_1^T  |\varphi'(t)|^2\rho(t)dt=
\int_1^T\left[ |u'|^2 (G-G(T))+u'u G'+\frac14 u^2 \frac{G'^2}{G-G(T)}\right]\rho(t) dt$$

En intégrant par partie le terme du milieu, on obtient :

$$\int_1^T  |\varphi'(t)|^2\rho(t)dt\ge \frac{u(1)^2-u(T)^2}{2}+
\frac14 \int_1^T u^2 \frac{G'^2}{G-G(T)}\rho(t) dt$$
On obtient alors le résultat en faisant tendre $T$ vers l'infini.
\endproof
Lorsque $\int^\infty \frac{dt}{\rho(t)}=\infty$ la même preuve fournit :
 
 \begin{lem} \label{hardyp}
Soit $\rho\, :\, [1,\infty[\rightarrow \R_+^*$ une fonction positive telle que
$$\int^\infty \frac{dt}{\rho(t)}=\infty.$$ Si on définit 
$$g(t)=\int_1^t \frac{d\tau}{\rho(\tau)},$$ alors pour tout $\varphi\in W^{1,2}_{loc}([1,\infty[)$ tel que
$\varphi'\in L^2([1,\infty[, \rho dt)$ et $\varphi(1)=0$ on a 
$$\frac14\int_1^\infty \left(\frac{g'(t)}{g(t)}\right)^2 |\varphi(t)|^2\rho(t)dt\le 
\int_1^\infty  |\varphi'(t)|^2\rho(t)dt.$$
\end{lem}

De ceci on en déduit les résultats de continuité pour les opérateurs

$$M_k(v)=\int_t^\infty v(s) \frac{s^{k-1}}{t^{k-1}}ds$$
$$\et \ m_k(v)=\int_1^t v(s) \frac{s^{k-1}}{t^{k-1}}ds$$

\begin{prop}\label{CM}
\begin{enumerate}
\item si $a>2k$ et $b\in \R$, alors l'opérateur $M_k$ est borné de 
$L^2([1,\infty[,r^{a-1}(1+\log(r))^{b}dr)$ dans $L^2( [1,\infty[ ,r^{a-3}(1+\log(r))^{b}dr)$,
\item si $a<2k$ et $b\in \R$, alors l'opérateur $m_k$ est borné de \\
$L^2([1,\infty[,r^{a-1}(1+\log(r))^{b}dr)$ dans $L^2([1,\infty[,r^{a-3}(1+\log(r))^{b}dr)$.
\item Si $a=2k$ et $b<1$ alors l'opérateur $m_b$ est borné de \\
$L^2([1,\infty[ ,r^{a-1}(1+\log r)^{b}dr)$ dans $L^2([1,\infty[,r^{a-3}(1+\log r)^{b-2}dr)$.
\item Si $a=2k$ et $b>1$ alors l'opérateur $M_k$ est borné de \\
$L^2([1,\infty[ ,r^{a-1}(1+\log r)^{b}dr)$ dans $L^2([1,\infty[,r^{a-3}(1+\log
r)^{b-2}dr)$.
\end{enumerate}
\end{prop}

On considère maintenant $(V,h)$ une variété riemannienne de dimension $d-1$ et on note 
$C(V)$ le cône sur $V$, c'est la variété $]0,\infty[\times V$ équipée de la métrique
$$g=dr^2+r^2\pi^*h$$ où $r,\pi$ sont les projections sur le premier et second facteur.
Pour $R>0$ on note $C_R(V)=\{r>R\}$ et $V_R=\{R\}\times V\subset C(V)$. Pour $R>0$, les variétés $C_R(V)$
sont quasi isométriques entre elles.

\subsection{Annulation de la $L^2$ cohomologie des cônes}

\begin{prop}
Si $k\le d/2$ alors $\H^k(C_1(V))=\{0\}$
et si $k\ge d/2$ alors $\H^k(C_1(V), V_1)=\{0\}$.
\end{prop}
Ce résultat est bien connu mais la preuve nous fournira des renseignements supplémentaires sur
$B^k(C_1(V))$.
\proof On considère le champ de vecteur $X=r\frac{\partial}{\partial r}$ et $\Phi_t$ son flot :
$$\Phi_t(r,\theta)=(e^tr,\theta).$$
Si $\alpha\in Z_2^k(C_1(V))$, grâce à la formule de Cartan, nous avons pour $t>0$ :
$$ \alpha-\Phi_t^*\alpha=d\beta_t\  {\rm avec\ }\  \beta_t=-\int_0^t\Phi_\tau ^*(\inte_X\alpha)d\tau.$$
où $\inte_X$ est l'opérateur produit intérieur par $X$.
Avec l'estimation :
$$|\Phi_t^*\alpha|(r,\theta)\le e^{kt} |\alpha|(e^t r,\theta),$$
nous obtenons :
$$\int_{C_{e^t}(V)} |\Phi_t^*\alpha|^2\le e^{(2k-d)t} \int_{C_{e^t}(V)} |\alpha|.$$
Ainsi pour $k\le d/2$ on a $\alpha=L^2\lim_{t\to \infty} d\beta_t$ et puisque $\beta_t\in L^2$
on a bien $\alpha\in B_2^k(C_1(V))$.
Remarquons que puisque 
$$ |\Phi_\tau ^*(\inte_X\alpha)|(r,\theta)\le e^{k\tau} r|\alpha| (e^\tau r,\theta),$$ l'on a aussi 
$$|\beta|(r,\theta)\le \frac{1}{r^{k-1}}\int_r^{e^t r} |\alpha|(s,\theta) s^{k-1} ds .$$
Et donc on a pour $k<d/2$ :
$\alpha=dB^+\alpha$ avec 
$$|B^+\alpha|(r,\theta)|\le\frac{1}{r^{k-1}}\int_r^{\infty} |\alpha|(s,\theta) s^{k-1} ds .$$
Et avec (\ref{CM}), nous obtenons:
\begin{equation}
\label{wght}\|r^{-1} B^+\alpha\|\le C \|\alpha\|.\end{equation}

Dans le second cas, si $\alpha\in Z_2^k(C_1(V),V_1 )$ alors l'extension par zéro de
$\alpha$ à $C(V)$ est encore fermée, on notera encore  $\alpha$ cette extension et si $t<0$ on a encore
$\Phi_t^*\alpha\in Z_2^k(C_1(V),V_1 )$. Grâce à la formule de Cartan, nous obtenons 
de même pour $t<0$ :
$$ \alpha-\Phi_t^*\alpha=d\beta_t$$
On a maintenant 
$$|\beta|(r,\theta)\le \frac{1}{r^{k-1}}\int_{\inf\{e^t r,1\}}^r |\alpha|(s,\theta) s^{k-1} ds .$$
$\alpha$ et $\Phi_t^*\alpha$ définissent donc la même classe dans $\H^k(C_1(V), V_1)$, mais pour $k\ge
d/2$, $\Phi_t^*\alpha$ tend faiblement vers $0$ dans $L^2$ lorsque $t\to -\infty$ (en même fortement
si $k>d/2$). Donc cette classe de cohomologie relative est nulle. Cela prouve la deuxième annulation.
On obtient ici 
$\alpha=dB^-\alpha$ avec 
$$|B^-\alpha|(r,\theta)|\le\frac{1}{r^{k-1}}\int_1^r |\alpha|(s,\theta) s^{k-1} ds .$$
 d'où d'après (\ref{CM}):
$$\|r^{-1} B^-\alpha\|\le C \|\alpha\|,\ {\rm lorsque\ } k>d/2 $$
$$\et \|(r \log r )^{-1} B^-\alpha\|\le C \|\alpha\|,\  {\rm lorsque\ }\  k=d/2. $$

\endproof

\subsection{$L^2$ cohomologie des cônes}

\begin{prop}\label{abstorel}On suppose de plus que sur $(V,h)$ l'image de $d$ est fermée. Alors pour $k>d/2$,
 l'application de restriction $r=\iota_{C_1(V)\setminus C_2(V),C_1(V)}$
induit un isomorphisme :

$$r^*\,:\, \H^k(C_1(V))\rightarrow \H^k(C_1(V)\setminus C_2(V))\simeq \H^k(V).$$

Et si de plus $V$ est une variété à bord alors pour $k<d/2$ alors l'application d'extension par zéro est un isomorphisme :
$$\H^k(C_1(V)\setminus C_2(V), V_1\cup V_2)\simeq \H^{k-1}(V)\rightarrow \H^k(C_1(V),V_1).$$
\end{prop}

\proof
L'image de $d$ est fermée sur $C_1(V)\setminus C_2(V)$
car $C_1(V)\setminus C_2(V)$ est quasi isométrique au produit $]1,2[\times V$ et de plus
 $\H^k( C_2(V),  V_2)=\{0\}$ pour $k\ge d/2$, donc d'après (lemme 2.2 de \cite{Carpar}), on sait alors que 
 $r^*$ est injectif.
 De plus $$\pi^* \,:\, L^2(\Lk V)\rightarrow L^2(\Lk C_1(V))$$ est justement continue pour $k>d/2$
et puisque $r^*\circ \pi^*=\Id$, $r^*$ est aussi surjectif.

On va maintenant démontrer le second isomorphisme : nous commençons par démontrer que cette application est
surjective :
soit $\alpha\in Z_2^k(C_1(V),V_1)$ où $k< d/2$. Nous avons obtenu : 
$$\alpha=dB^+\alpha$$ et si
$\rho$ est une fonction lisse sur $C_1(V)$ valant $1$ sur un voisinage de $C_1(V)\setminus C_{3/2}(V)$
et nulle sur $C_{7/4}(V)$ alors 
$$\alpha=d\left(\rho B^+\alpha\right)+ d\left((1-\rho)B^+\alpha\right)$$
mais puisque grâce à l'estimée (\ref{wght}) et à (\ref{poide}), on sait que
$d[(1-\rho)B^+\alpha]$ est nul dans $\H^k(C_1(V),V_1)$
et puisque $r^* d(\rho B^+\alpha)\in Z^k_2(C_1(V)\setminus C_2(V), V_1\cup V_2)$, on a bien montré que cette
application d'extension par zéro est surjective. On montre maintenant qu'elle est injective :  soit donc 
$\alpha \in \cH^{k-1}(V)$ alors la forme $\frac{dr}{r^{d-2k+1}}\wedge\pi^*\alpha$ est 
harmonique sur $C_1(V)$, c'est un élément non nul de $\cH^k(C_1(V))$ qui représente donc un élément non nul
de $\H^k(C_1(V))$.

\endproof

\subsection{Presque fermeture de l'image de $d$}
On suppose maintenant que $V$ est une variété compacte à bord éventuellement non vide.

Pour $k\not= d/2$, , nous notons $w_k=r^{-2}$ et $w_{d/2}=(r(1+\log r))^{-2}$ alors ces poids sont
évidemment à décroissance parabolique sur $C_1(V)$ et nous avons montré que
pour $k\ge d/2$  on a
$$d\cC^{k-1}_{ w_k,1}(C_1(V),V_1)=B^k_1(C_1(V),V_1)$$
et pour $k<d/2$ on a :
$$d\cC^{k-1}_{w_k,1}(C_1(V))=B^k_1(C_1(V))$$

En fait cette dernière identité est encore valable pour $k\ge d/2$. En effet soit 
$\alpha\in B^{k}_1(C_1(V))$, il existe donc $\psi\in L^2_\loc$ tel que
$\alpha=d\psi$, on a alors pour la fonction $\rho$ précédente :
$$\alpha=d\rho\psi +d\left(B^-(\alpha-d\rho\psi)\right)$$
On a bien $B^-(\alpha-d\rho\psi), \rho\psi\in \cC^{k-1}_{w_{k},1}(C_1(V))$.

\noindent Notons que l'argument donné ici a une portée plus générale :

\begin{prop}\label{exbord1}
Si $d\,:\, \cC^{k-1}_{w,\mu}(X,Y)\rightarrow L^2_\mu(\Lk X)$ a son image fermée et si
\begin{itemize}
\item $w$ est à décroissance parabolique,
\item $Y\subset \partial\bar X$ a un voisinage $V$ dans $\bar X$ quasi isométrique à $Y\times [0,1[$ avec $\partial Y\times
 [0,1[\subset\partial\bar X$ et que $d\,:\, \cC^{k-1}_{w,\mu}(V)\rightarrow L^2_\mu(\Lk V)$ a  
 son image fermée 
\item $\H_\mu^k(X,Y)=\{0\},$
\end{itemize}
alors $d\,:\, \cC^{k-1}_{w,\mu}(X)\rightarrow L^2_\mu(\Lk X)$ a aussi son image fermée.
\end{prop}
On peut aussi montrer le résultat suivant 
\begin{prop}\label{exbord2}
Si $d\,:\, \cC^{k-1}_{w,\mu}(X)\rightarrow L^2_\mu(\Lk X)$ a son image fermée et si
\begin{itemize}
\item $Y\subset \partial\bar X$ est compact et
 a un voisinage dans $\bar X$ quasi isométrique à $Y\times [0,1[$ et $\partial Y\times
 [0,1[\subset\partial\bar X$.
  
 \item  l'application naturelle $\H^{k-1}(Y)\rightarrow \H^k(X,Y)$ est injective,
\item $w$ est à décroissante parabolique.
\end{itemize}
alors $d\,:\, \cC^{k-1}_{w,\mu}(X,Y)\rightarrow L^2_\mu(\Lk X)$ a aussi son image fermée.
\end{prop}

\proof  Soit donc $\alpha\in B^k_\mu(X,Y)$, on trouve alors $\psi\in \cC^{k-1}_{w,\mu}(X)$ tel que
$$\alpha=d\psi.$$ Soit $\rho$ une fonction lipschitzienne à support compact dans $Y\times [0,1[$ et valant
$1$ près de $Y=Y\times \{0\}$. Nous avons 
aussi 
$$\alpha=d(\rho\psi)+d[(1-\rho)\psi].$$
Puisque $(1-\rho)\psi\in \cC^{k-1}_{w,\mu}(X,Y)$,
$d(\rho\psi)$ est donc dans le noyau de l'application naturelle
$$\H^{k-1}(Y)\simeq\H^{k}(Y\times [0,1],Y\times \{0,1\} )\rightarrow \H^k(X,Y),$$
il est donc nulle par hypothèse. Il existe donc $\phi\in \cD^{k-1}(Y\times [0,1],Y\times \{0,1\})$ tel que
$d(\rho\psi)=d\phi$. Si $\bar\phi$ est l'extension par zéro de $\phi$ à $X$, on a bien
$$\alpha=d\bar\phi+d[(1-\rho)\psi]$$ et 
$\bar\phi+(1-\rho)\psi\in \cC^{k-1}_{w,\mu}(X,Y)$.
\endproof

Grâce à ces résultats et au corollaire (\ref{apoid}), on obtient donc 
\begin{prop}
 Si $w=(r(1+\log r))^{-2}$, alors
$$d\cC^{k-1}_{w,1}(C_1(V),V_1)=B^k_1(C_1(V),V_1)\et d\cC^{k-1}_{w,1}(C_1(V))=B^k_1(C_1(V))$$
\end{prop}
\subsection{$L^2$ cohomologie à poids}

Les arguments présentés auparavant se généralisent aisément à la cohomologie à poids :
soit $a,b\in \R$, on introduit $L^2_{a,b}(\Lk C_1(V))$ l'espace 
des formes différentielles de carrés sommables pour 
la mesure $r^{2a} (1+\log r)^{2b}d\vol _g$ et on notera 
$$\cC^{k-1}_{a,b}(C_1(V))=\{\alpha\in L^2_{a-1,b-1}(\Lkm C_1(V)), d\alpha\in L^2_{a,b}(\Lk C_1(V))\} $$
et $$\tilde \cC^{k-1}_{a,b}(C_1(V))=\{\alpha\in L^2_{a-1,b}(\Lkm C_1(V)), d\alpha\in L^2_{a,b}(\Lk
C_1(V))\} $$
 et on notera de la même façon :  $\cC^{k-1}_{a,b}(C_1(V), V_1)$, $\tilde\cC^{k-1}_{a,b}(C_1(V), V_1)$, $\H^k_{a,b}(C_1(V)) $,
  $\H^k_{a,b}(C_1(V), V_1)$... La même preuve montre :  
 \begin{thm}\label{cohopcone}Lorsque $V$ est une variété compacte à bord lisse et
 $(k,\frac12)\not=(\frac{d}{2}+a,b)$ alors 
 $$d\,:\, \cC^{k-1}_{a,b}(C_1(V))\rightarrow B^k_{a,b}(C_1(V) )\et 
 d\,:\, \cC^{k-1}_{a,b}(C_1(V), V_1)\rightarrow B^k_{a,b}(C_1(V), V_1) $$  ont  des images fermées.
 Même mieux si $k\not=\frac{d}{2}+a$ alors
 $$d\,:\, \tilde\cC^{k-1}_{a,b}(C_1(V))\rightarrow B^k_{a,b}(C_1(V) )\et 
 d\,:\, \tilde\cC^{k-1}_{a,b}(C_1(V), V_1)\rightarrow B^k_{a,b}(C_1(V), V_1) $$  ont  des images fermées.
 De plus \begin{itemize}
 \item Si $(k, -\frac12)<(\frac{d}{2}+a,b)$\footnote{l'ordre est ici l'ordre lexicographique}, alors
 $\H^k_{a,b}(C_1(V)) =\{0\}$.
 \item Si  $(k, -\frac12)>(\frac{d}{2}+a,b)$, alors
 $\H^k_{a,b}(C_1(V)) =\H^k(V)$.
 \item Si  $(k, \frac12)>(\frac{d}{2}+a,b)$, alors
 $\H^k_{a,b}(C_1(V), V_1) =\{0\}$.
 \item Si  $(k, \frac12)<(\frac{d}{2}+a,b)$, alors
 $\H^k_{a,b}(C_1(V), V_1) =\H^{k-1}(V)$.
\end{itemize} \end{thm}
\proof Les mêmes arguments que précédemment montrent que pour $(k,\frac12)<(\frac{d}{2}+a,b)$, alors
$d\cC^{k-1}_{a,b}(C_1(V))=Z^k_{a,b}(C_1(V))$ et pour $(k, \frac12)>(\frac{d}{2}+a,b)$, alors 
$d\cC^{k-1}_{a,b}(C_1(V), V_1)=Z^k_{a,b}(C_1(V), V_1)$. On obtient de la même façon les isomorphismes :
\begin{itemize}\item si  $(k, -\frac12)>(\frac{d}{2}+a,b)$, alors
 $\H^k_{a,b}(C_1(V)) =\H^k(V)$,
\item si  $(k, \frac12)<(\frac{d}{2}+a,b)$, alors
 $\H^k_{a,b}(C_1(V), V_1) =\H^{k-1}(V)$.
\end{itemize}
Et grâce à (\ref{exbord1},\ref{exbord2}),
 on conclut de même que pour $(k,\frac12)\not=(\frac{d}{2}+a,b)$ alors l'image de $d$ est presque fermée en tout degré
 par rapport à $w$ (ou même à $\bar w =r^{-2}$ si $k\not=\frac{d}{2}+a$).
Il reste donc à montrer que si $k=\frac{d}{2}+a$ et $b\in [-\frac12,\frac12]$, alors
$\H^k_{a,b}(C_1(V)) =\{0\}$. L'application de restriction 
$r^*\,:\, \H_{a,b}^k(C_1(V))\rightarrow \H^k(C_1(V)\setminus C_2(V))\simeq \H^k(V)$  est 
comme précédemment injective.
Il suffit donc de montrer que pour ces valeurs de $k,a,b$, cette application est nulle.
Pour cela, on montre que si $\alpha\in Z^k_{a,b}(C_1(V))$ et $\beta\in \cH^{d-1-k}(V_1,\partial V_1)$
alors
$$c=\int_{V_1}\iota_1^*\alpha\wedge \beta=0,$$
où on note $\iota_r: \{1\}\times V\rightarrow \{r\}\times V$.
Pour cela on procède comme dans la preuve de la proposition (3.7) de \cite{Carpar}:
on a bien sur
\begin{equation*}
\begin{split}
c^2&=\left(\int_{V_1}\iota_t^*\alpha\wedge \beta\right)^2\\
&\le \|\beta\|_{L^\infty}^2 t^{2k}\int_{V} |\alpha|^2(t,\theta) d\vol_{h};\\
\end{split}
\end{equation*}
D'où :
\begin{equation*}
c^2\int_1^\infty t^{2a-2k} t^{d-1}(1+\log t)^{2b}dt \le \|\beta\|_{L^\infty}^2\|\alpha\|_{L^2_{a,b}}^2.\end{equation*}
Si $k=\frac{d}{2}+a$ et $b\ge -\frac12$, on obtient bien la nullité de $c$.

\endproof
Nous devons néanmoins améliorer quelque peu ce résultat lorsque $k=\frac{d}{2} +a$ :

\begin{prop}\label{raffinement}
Soit $V$ une variété compacte à bord lisse et
 $k=\frac{d}{2}+a$ , $b\not=\frac12$ alors:
 \begin{enumeroman}
 \item Si le $(k-1)$-groupe de cohomologie de $V$ est nulle $b_{k-1}(V)=0$,
 alors $$d\,:\,\tilde \cC^{k-1}_{a,b}(C_1(V))\rightarrow B^k_{a,b}(C_1(V) )\et 
 d\,:\, \tilde\cC^{k-1}_{a,b}(C_1(V), V_1)\rightarrow B^k_{a,b}(C_1(V), V_1) $$  ont  des images fermées.
 \item Si $b_{k-1}(V)\not=0$,  notons $\cH^k(V,h)$ l'espace des formes harmoniques $L^2$ sur $V$ qui
 vérifient la condition absolue:
 $$\cH^k(V,h)=\{\alpha\in L^2(\Lk V),d\alpha=0, d^*\alpha=0, \inte_{\vec\nu}\alpha=0\}$$
 où on a noté $\vec \nu \, :\,\partial V\rightarrow TV$ le champ normal unitaire entrant de $V$.
 Alors si $\alpha\in B^k_{a,b}(C_1(V))$  , alors il existe
 $u\in \tilde \cC^{k-1}_{a,b}(C_1(V))$  et
 $v\in\cC^{k-1}_{a,b}(C_1(V))$ telle que
 $$\inte_{\frac{\partial}{\partial r}}v=0\ \et \forall r>1, \iota^*_r v\in\cH^k(V,h)$$
 tels que
 $$\alpha=du+dv.$$
 De même si $\alpha\in B^k_{a,b}(C_1(V),V_1 )$ alors il existe
  $u\in \tilde\cC^{k-1}_{a,b}(C_1(V), V_1)$ et
 $v\in\cC^{k-1}_{a,b}(C_1(V), V_1)$ tel que
 $$\inte_{\frac{\partial}{\partial r}}v=0\ \et \forall r>1, \iota^*_rv\in\cH^k(V,h)$$
avec
 $$\alpha=du+dv.$$
 \end{enumeroman}
\end{prop}
\proof On ne montre que la deuxième assertion (la première en étant un cas particulier) pour le cas absolu et
$b>1/2$
(le cas relatif à $V_1$ et $b<1/2$ est identique à une modification évidente près). Dans les autres cas s'en
déduisent grâce aux propositions (\ref{exbord1},\ref{exbord2}). 

Pour la commodité des notations, on
se sert du difféomorphisme :
$\Phi\,:\, \R_+\times V\rightarrow C_1(V)$ défini par 
$$\phi(t,\theta)=(e^t,\theta)$$ qui est conforme lorsque $\R_+\times V$ est équipé de la métrique riemannienne
produit $(dt)^2+h$. Alors $\Phi^*$ est isométrie entre $L^2_{a,b}(\Lk C_1(V))$ et 
$L^2_{\rho_{k,a,b}}(\Lk(\R_+\times V))$ où $\rho_{k,a,b}(t,\theta)=e^{(2a-2k+d)t}(1+t)^{2b}$.
Pour le degré $k=\frac{d}{2}+a$, nous avons donc $\rho_{k,a,b}:=\rho_b=(1+t)^{2b}$.
On travaille donc sur $\R_+\times V$. Et on identifie :
$$L^2_\loc(\Lk(\R_+\times V))\simeq dt\wedge L_\loc^2(\R_+, \Lkm(V))\oplus L_\loc^2(\R_+, \Lk(V)).$$

Soit donc $\alpha\in B^k_{\rho_b}(\R_+\times V)$, on peut donc écrire :
$$\alpha=dt\wedge \alpha_0+\alpha_1$$ avec 
$$\alpha_0\in L^2(\R_+, \Lkm(V))\ \et \ \alpha_1\in L^2(\R_+, \Lk(V))$$
on sait (cf. \ref{cohol2loc}) qu'il existe $\beta=dt\wedge \beta_0+\beta_1\in L^2_\loc(\Lk(\R_+\times V)) $ 
tel que $$\alpha=d\beta.$$ Si on note $d_0$ l'opérateur de différentiation extérieur sur $V$ alors
$$\left\{\begin{array}{l}
\alpha_1(t)=d_0\beta_1(t)\\
\alpha_0(t)=\beta_1'(t)-d_ 0\beta_0(t)\\
\end{array}\right.$$ 

On veut modifier $\beta$ pour qu'il satisfasse nos conditions.
Nous utilisons maintenant la décomposition de Hodge de $L^2(\Lambda^\bullet T^* V)$ ; notons 
$$\cD^k=\{u\in W^{1,2}(\Lambda^k T^* V),  \inte_{\vec\nu}u=0\ {\rm le\ long\ de}\ \partial
V\}$$
où nous avons noté $W^{1,2}(\Lambda^\bullet (T^* V))$ l'espace de Sobolev des sections de $\Lambda^\bullet( T^*
V)$ qui sont dans $L^2$ ainsi que leurs premières dérivées.
Nous savons que si $\alpha\in L^2((\Lk V)$ alors il existe un unique triplet $(h,f,g)\in \cH^k(V)\oplus
\cD^{k-1}\oplus \cD^{k+1}$ telle que
$$ \alpha=h+d_0f+d_0^*g\ \et\ f,g\perp \cH^\bullet(V),\ d_0^*f=d_0g=0,$$
de plus il existe une constante $C$ telle que
$$\|d_0f\|+\|f\|+\|d_0^*g\|+\|g\|\le C \|\alpha\|$$
où on a noté ici $\|\, .\, \|$ la norme $L^2$ sur $V$.
Nous décomposons alors pour $i=0,1$
$$\alpha_i(t)=d_0u_i(t)+d^*_0v_i(t)+w_i(t)\ \et \ \beta_i(t)=d_0f_i(t)+d^*_0g_i(t)+h_i(t)$$
Et nous obtenons les équations :
$$v_1=w_1=0$$ et 
\begin{equation}\label{eq1}
\alpha_1=d_0u_1=d_0d_0^*g_1,
\end{equation}

\begin{equation}
\label{eq2}
h_1'=w_0,
\end{equation}

\begin{equation}\label{eq3}
d^*_0g_1'=d^*_0v_ 0
\end{equation}

\begin{equation}\label{eq4}
d_0f_1'-d_0d^*_0g_0=d_0u_0
\end{equation}
Afin de résoudre (\ref{eq4}), nous supposons que
$$df_ 1=g_0'$$
et alors  (\ref{eq4}) devient :
\begin{equation}\label{eq4new}
g_0''-d_0d^*_0g_0=g_0''-(d_0d^*_0+d^*_0d_0)g_0=d_0u_0.
\end{equation}
Mais pour $\mu>0$ l'opérateur
$$L_\mu f=f''-\mu f$$
 sur $L^2(\R_+,(1+t)^{2b} dt)$ avec les conditions de Dirichlet en $t=0$ est inversible 
et qu'il existe une constante $C$ indépendante de $\mu$ tel que si $L_\mu f=g$ avec $f(0)=0$ alors
\begin{equation*}
(1+|\mu|)\|f\|_b+\|f'\|_b\le C\|g\|_b
\end{equation*}
$\|\ \|_b$ étant la norme de $L^2(\R_+,(1+t)^{2b} dt)$.

Ainsi en utilisant la décomposition spectrale de l'opérateur $(d_0d^*_0+d^*_0d_0)$ sur
$L^2(\Lambda^\bullet(TV))$ (pour les conditiosn absolu sur le bord de $V$) on trouve une solution unique $g_0$ de (\ref{eq4new}) telle que

\begin{equation}\label{est1}
\begin{split}
& g_0(0)=0\\
&\|g_0'\|_{\rho_b}+\|d_0^*g_0\|_{\rho_b}+\|g_0\|_{\rho_b}\le C\|du_0\|_{\rho_b}\\
\end{split}
\end{equation}
Nous posons alors
$$v(t)=-\int_t^\infty w_0(s)ds$$ et 
$$u=dt\wedge d_0^*g_0+g_0'+u_ 1,$$ nous obtenons alors
$$\alpha=du+dv$$
de plus $u$ et $v$ vérifient bien les propriétés requises :
$u\in L^2_{\rho_b}$ et $v\in L^2_{\rho_{b-1}}$
\endproof
\subsection{Cohomologie $L^2$ des variétés à bouts coniques :}
Nous pouvons grâce à ces calculs déterminer la cohomologie $L^2$ à poids des variétés dont un voisinage de
l'infini est un bout de cône.

\begin{thm}\label{cohosc} Soit $(X^d,g)$ une variété riemannienne dont un voisinage de l'infini est isométrique à un bout
de cône $C_1(V)$ sur une variété riemannienne compacte $V^{d-1}$ et soit $\mu,w$ des fonction lisses sur $X$
telles que sur $C_1(V)$
$$\mu(r,\theta)=r^{2a}\ \et \ w(r,\theta)=r^{-2} (1+\log r)^{-2}.$$
 Alors l'image de $d\,:\, \cD^k_\mu(d)\rightarrow L^2_\mu$ est presque fermée en tout degré par
rapport à $w$ et on a
\begin{equation*}
\H^k_\mu(X)=\left\{\begin{array}{ll}
H_c^k(X) & {\rm\ si\ } k<d/2+a\\
\ima\left(H_c^k(X)\rightarrow H^k(X)\right)  & {\rm\ si\ } k=d/2+a\\
H^k(X) & {\rm\ si\ } k>d/2+a\\
\end{array}\right. .
\end{equation*}
\end{thm}

Ce résultat est en fait bien connu, une preuve différente est faite dans l'article de T. Hausel, E. Hunsicker et
R. Mazzeo (\cite{HHM}); on peut aussi en donner une preuve avec les arguments de N. Yeganefar \cite{Y} et de
\cite{Cargafa}.
Nous montrons ici comment utiliser nos résultats précédents pour prouver ce théorème :
\proof
Nos résultats précédents impliquent que les images des applications 
$$d\,:\, \cC^{k-2}_{w,w\mu}(C_1(V))\rightarrow L^2_{w\mu}\left(\Lkm \left(C_1(V)\right)\right)$$ et de
$$d\,:\, \cC^{k-1}_{w,\mu}(C_1(V))\rightarrow L^2_{\mu}\left(\Lk \left(C_1(V)\right)\right)$$ sont fermées ;
 de plus $\H^{k-1}_{w\mu}(C_1(V))$ est bien de dimension finie.  Grâce à (\ref{conclusion}), 
 nous pouvons affirmer que l'image de $d$ est presque fermée sur $X$ en tout degré par
rapport à $w$. Si $K\subset X$ est le compact dont le complémentaire est $C_1(V)$ alors on a la suite exacte
de Mayer Vietoris :
\begin{equation}
\label{MVcone}
\begin{split}
 \H^{k-1}(K)\oplus \H^{k-1}_{w\mu}(C_1(V))
 \stackrel{\delta}{\longrightarrow} &\H^{k-1}(\partial K)\\
 \stackrel{\mbox{b}}{\longrightarrow}&
 \H^k_{\mu}(X) \stackrel{\mbox{r}^*}{\longrightarrow} 
 \H^k(K)\oplus \H^k_{\mu}(C_1(V))\stackrel{\delta}{\longrightarrow} \H^k_{\mu}(\partial K).\\
\end{split}\end{equation}

Lorsque $k<d/2+a$, on sait que $$\H^{k-1}_{w\mu}(C_1(V))=\H^k_{\mu}(C_1(V))=\{0\},$$ cette suite exacte
se réduit donc à :
\begin{equation*}
 \H^{k-1}(K)
 \stackrel{\delta}{\longrightarrow} \H^{k-1}_{}(\partial K)\stackrel{\mbox{b}}{\longrightarrow}
 \H^k_{\mu}(X) \stackrel{\mbox{r}^*}{\longrightarrow} 
 \H^k(K)\stackrel{\delta}{\longrightarrow} \H^k(\partial K).
\end{equation*}
Puisqu'on a aussi la suite exacte :
\begin{equation*}
 \H^{k-1}(K)
 \stackrel{\delta}{\longrightarrow} \H^{k-1}_{}(\partial K)\stackrel{\mbox{b}}{\longrightarrow}
 \H^k(K,\partial K) \stackrel{\mbox{r}^*}{\longrightarrow} 
 \H^k(K)\stackrel{\delta}{\longrightarrow} \H^k(\partial K),
\end{equation*}
Le lemme des cinq flèches impliquent que l'application naturelle
$$ \H^k(K,\partial K)\simeq H_c^k(X)\rightarrow \H^k_\mu(X)$$ est un isomorphisme.

Lorsque $k>d/2+a$, on sait que
$$\H^{k-1}_{w\mu}(C_1(V))=H^{k-1}(\partial K)\ \et\ \H^k_{\mu}(C_1(V))=H^{k-1}(\partial K),$$
donc les premières et dernières flèches de la suite (\ref{MVcone}) sont surjectives. On en déduit que
l'application de restriction 
$$ \H^k(X)\rightarrow  H^k(K)\simeq H^k(X)$$ est un isomorphisme.

Lorsque $k=d/2+a$, on a montré que
$$\H^{k-1}_{w\mu}(C_1(V))=H^{k-1}(\partial K)\ \et\ \H^k_{\mu}(C_1(V))=\{0\},$$
Dans la suite exacte (\ref{MVcone}) la troisième flèche est donc injective et donc
 $$ \H^k(X)\simeq \ker\left(  H^k(K)\rightarrow H^k(\partial K)\right) =\ima\left(  H^k(K,\partial
 K)\rightarrow H^k( K)\right). $$
\endproof
\begin{rem}
Nous avons en fait un peu mieux si $k\not=d/2+a$ ou si $k=d/2+a$ et $b_{k-1}(V)=0$ alors
l'image de $d\,:\, \cD^k_\mu(d)\rightarrow L^2_\mu$ est presque fermée en degré $k$ par
rapport à $$\bar w(r,\theta)=r^{-2}.$$
\end{rem}
Nous utiliserons plus loin uniquement le cas où $V$ est le quotient de la sphère $\bS^{d-1}$ par un sous
groupe fini $G\subset {\rm SO}(d)$ ($G$ agit donc sans point fixe sur $\bS^{d-1}$) on dit alors que $X$ est
une variété ALE (Asymptotiquement Localement Euclidienne) à un seul bout:
\begin{cor} 
\label{injALE}
Soit $(X^d,g)$ une variété ALE à un seul bout alors on a pour la cohomologie $L^2$ de $X$ :
\begin{equation*}
\H^k(X)=\ima\left(H_c^k(X)\rightarrow H^k(X)\right) =\left\{\begin{array}{ll}
H_c^k(X) & {\rm\ si\ } k<d\\
H^k(X) & {\rm\ si\ } k>0\\
\end{array}\right. .
\end{equation*}
En particulier, l'application $ \H^k(X)\rightarrow  H^k(X)$ est toujours injective,
et l'application $ H_c^k(X)\rightarrow  \H^k(X)$ est toujours surjective.
De plus si $d>2$, alors l'image de $d\,:\, \cD^k(d)\rightarrow L^2$ est presque fermée en tout degré  par
rapport à $\bar w$.
\end{cor}
\subsection{Une formule de Künneth}

On considère deux variétés riemanniennes $(X_1,g_1)$ et $(X_2,g_2)$ et la
 variété $X=X_1 \times X_2$ équipée de la métrique produit. 
 Un corollaire direct des arguments de S. Zucker (\cite{Z1}) est le suivant :
 \begin{prop}\label{kunneth1}
 Supposons que l'image de $d\,:\, \cD_{\mu_1}(X_1)\rightarrow L_{\mu_1}^2(\Lambda^\bullet T^* X_1)$ soit
 fermée et que l'image de $d\,:\, \cC_{w_2,\mu_2}(X_2)\rightarrow L_{\mu_2}^2(\Lambda^\bullet T^* X_2)$ soit
fermée alors si $\mu(x_1,x_2)=\mu_1(x_1)\mu_2(x_2)$ et $w(x_1,x_2)=w_2(x_2)$ l'image de
$$ d\,:\, \cC_{w,\mu}(X)\rightarrow L_{\mu}^2(\Lambda^\bullet T^* X)$$ est
fermée de plus 
$$\H_\mu^k(X)=\bigoplus_{p+q=k} \H^p_{\mu_1}(X_1)\otimes \H^p_{\mu_2}(X_2).$$
 \end{prop}

Il n'est néanmoins pas clair que si l'image de $d$ est presque fermée sur chacun des facteurs alors elle
est presque fermée sur le produit. Il est possible qu'en utilisant les techniques introduites par R. Mazzeo et A.
Vasy (\cite{MV}), on puisse dans certains cas démontrer de tels résultats.
En fait, nous aurons uniquement
besoin d'une version très affaiblie :

\begin{prop}\label{kunneth2}
On suppose que l'image de $d\,:\, \cC^\bullet_{w_1,\mu_1}(X_1)\rightarrow L_{\mu_1}^2(\Lambda^\bullet T^* X_1)$
est fermée. On considère $X=X_1\times C_1(\bS^{d-1})$ équipé de la métrique riemannienne produit et de la
mesure $\mu(x_1,x_2)=\mu_1(x_1)$, on pose $\bar w_2(x_1,x_2)=|x_ 2|^{-2}$ 
et $ w_2(x_1,x_2)=|x_ 2|^{-2}(1+\log|x_ 2|)^{-2}$ et on définit
 $$Y=X_ 2\times \partial C_1(\bS^{d-1})=X_ 2\times \bS^{d-1}\subset \partial X$$  alors
\begin{enumeroman}
\item si $d>2$ : alors pour $\alpha\in B^k_{\mu}(X)$, 
il y a $\beta_2\in \cC^{k-1}_{\bar w_2,\mu}(X) $ et
 $\beta_1\in \cC^{k-1}_{ w_1,\mu}( X)$  
tels que $$\alpha=d\beta_1+d\beta_2,$$ de plus
$$\H_\mu^k(X)=\bigoplus_{p+q=k} \H^p_{\mu_1}(X_1)\otimes
\H^p(C_1(\bS^{d-1}))=\H^{k-d+1}_{\mu_1}(X_1)$$ de même pour
$\alpha\in B^k_{\mu}(X,Y)$, 
il y a $\beta_2\in \cC^{k-1}_{\bar w_2,\mu}(X,Y) $ et
 $\beta_1\in \cC^{k-1}_{ w_1,\mu}( X,Y)$ tels que
 $$\alpha=d\beta_1+d\beta_2,$$ de plus 
$$\H_\mu^k(X,Y)=\bigoplus_{p+q=k} \H^p_{\mu_1}(X_1)\otimes
\H^p(C_1(\bS^{d-1}), \bS^{d-1})=\H^{k-1}_{\mu_1}(X_1) $$
\item Si $d=2$ : alors pour $\alpha\in B^k_{\mu}(X)$, 
il y a $\beta_2\in \cC^{k-1}_{ w_2,\mu}(X) $ et
 $\beta_1\in \cC^{k-1}_{ w_1,\mu}( X)$  
tels que $$\alpha=d\beta_1+d\beta_2,$$  de même pour
$\alpha\in B^k_{\mu}(X,Y)$, 
il y a $\beta_2\in \cC^{k-1}_{w_2,\mu}(X,Y) $ et
 $\beta_1\in \cC^{k-1}_{ w_1,\mu}( X,Y)$ tels que
 $$\alpha=d\beta_1+d\beta_2,$$ de plus 
$$\H_\mu^k(X)=\{0\}\ \et\ \H_\mu^k(X,Y)=\{0\}$$
\end{enumeroman} 
\end{prop}
On commence par démontrer le lemme suivant :
\begin{lem}
\begin{enumeroman}
\item Soit $d>2$, il existe un opérateur continu 
$$B\,:\, L^2(C_1(\bS^{d-1})) \rightarrow L^2_{\bar w_2}(C_1(\bS^{d-1}))$$ tel que
pour si $\alpha\in L^2(C_1(\bS^{d-1}))$, alors $B\alpha\in W^{1,2}_{\loc}$ et
$\inte_{\frac{\partial}{\partial r}}B\alpha=0$ sur $\bS^{d-1}=\partial C_1(\bS^{d-1})$. De plus
si $\alpha\in \cD(d)$ alors 
$$\alpha=h(\alpha)+dB\alpha+Bd\alpha,$$
où $h$ est un projecteur continu sur $\cH^\bullet(C_1(\bS^{d-1}))=\R*\frac{dr}{r^{d-1}}$.
De la même façon, il existe un opérateur continu
 $$B_{rel}\,:\, L^2(C_1(\bS^{d-1})) \rightarrow L^2_{\bar w_2}(C_1(\bS^{d-1}))$$
 tel que si $\alpha\in L^2(C_1(\bS^{d-1}))$, alors $B\alpha\in W^{1,2}_{\loc}$ et
$\iota^*_{\partial C_1(\bS^{d-1})}B\alpha=0$.
si $\alpha\in \cD(d,\bS^{d-1})$ alors 
$$\alpha=h(\alpha)+dB\alpha+Bd\alpha$$
où $h$ est le projecteur orthogonal sur $\cH^\bullet(C_1(\bS^{d-1}),\bS^{d-1})=\R\frac{dr}{r^{d-1}}$.

\item Si $d=2$, alors il existe un opérateur continu $B_{rel}\,:\, L^2\rightarrow L^2_{ w_2}$ vérifiant les
mêmes propriétés que si dessus et tel que
si $\alpha\in \cD(d,\bS^{d-1})$ alors 
$$\alpha=dB\alpha+Bd\alpha$$
\end{enumeroman} 
\end{lem}
\proof On montre d'abord le premier cas et on expliquera les aménagements nécessaires pour démonter les deux
autres cas.

Grâce à l'opérateur $$\Delta_{abs}=d_{abs}(d_{abs})^*+(d_{abs})^*d_{abs}$$ on $d_{abs}$ est la restriction de
$d$ à $\cD(C_1(\bS^{d-1}))$, on bâtit un opérateur continu
$$B_1=(d_{abs})^*\int_0^1e^{-s\Delta_{abs}}ds\,:\, L^2(\Lk (C_1(\bS^{d-1})))\rightarrow
\cD^{k-1}(d)\cap \cD^{k-1}(d_{abs}^*)$$  
De plus, nous avons pour $\alpha\in \cD(C_1(\bS^{d-1}))$
$$\alpha=S\alpha+dB\alpha+Bd\alpha$$
où $S=e^{-\Delta_{abs}}$ commute avec $d_{abs}$, de plus $S\alpha$ est un forme $C^\infty$ sur
$C_1(\bS^{d-1})$.

On peut alors comme dans \cite{carma} effectuer la moyenne de $S\alpha$ suivant ${\rm SO} (d)$, on obtient 
$$S\alpha=M\alpha+dB_2\alpha+B_2d\alpha$$
où $$B_2\,:\, L^2(C_1(\bS^{d-1}))\rightarrow L^2_{\bar w_2}(C_1(\bS^{d-1}))$$ est continue.
Et $M\alpha$ est une forme ${\rm SO} (d)-$invariante elle est donc de la forme :
$$M\alpha=u_0+u_1dr+u_{d-1} *dr+u_d *1$$
où $u_0,u_1,u_{d-1},u_d$ sont des fonctions $C^\infty$ de carré sommables sur $[1,\infty[$, de plus
$u_1(1)=u_d(1)=0$. On pose alors
$$B_3\alpha=-\left(\int_r^\infty u_1(s)ds\right)+\left(\int_1^r u_d(s)\frac{ds}{s^{d-1}}\right)\frac{*dr}{r^{d-1}}$$
Grâce à nos estimées (\ref{CM}), l'opérateur $B_3\,:\, L^2\rightarrow L^2_{\bar w_2}$ est continue de plus
il est clair que
$$M\alpha=u_{d-1}(1)\frac{*dr}{r^{d-1}}+dB_3\alpha+B_3d\alpha$$
L'opérateur $B=B_1+B_2+B_3$ convient alors.

Dans le cas relatif, la même preuve convient à condition de considérer l'opérateur $\Delta_{rel}$ et de
modifier le dernier opérateur 
$B_3$ par l'opérateur
$$B_{3,rel}\alpha=\frac{1}{r^{d-2}}\left(\int_1^\infty u_1(s)ds\right)
-\left(\int_r^\infty u_1(s)ds\right)+\left(\int_1^r u_d(s)\frac{ds}{s^{d-1}}\right)\frac{*dr}{r^{d-1}}$$

Lorsque $d=2$, il faut considérer l'opérateur 
$$\tilde B_{3,rel}\alpha=\left(\int_1^r u_1(s)ds\right)
+\left(\int_1^r u_d(s)\frac{ds}{s^{d-1}}\right)\frac{*dr}{r^{d-1}}.$$
\endproof
On peut maintenant prouver la proposition :

\proof On le fait d'abord dans le cas où $d>2$. Suivant S. Zucker (preuve du théorème 2.29 de \cite{Z1}),
 on déduit du lemme précédent un opérateur borné
$$\bar B\,:\,  L_{\mu}^2(\Lk X)\rightarrow L_{\bar w_2\mu}^2(\Lkm X),$$
tel que si $\alpha\in Z^k_\mu(X)$ alors
$$\alpha=\bar h(\alpha)+d\bar B\alpha,$$
Où $\bar h(\alpha)$ est une projection  de $\alpha$ sur les formes qui sont
dans $ L^2_{\mu_1}(\Lambda^\bullet T^*X_1)\otimes \cH^\bullet(C_1(\bS^{d-1}))$.
A proprement parlé, nous n'obtenons pour $\alpha\in Z^k_\mu(X)$ directement que l'égalité
$$\rho\alpha=\bar h(\rho\alpha)+d\bar B(\rho\alpha)+\bar B(d\rho\wedge\alpha)$$
pour toutes les fonctions $\rho$ qui ne dépendent que de la seconde variable
et dont le support est inclus un voisinage borné de $X_1\times (\partial(C_1(\bS^{d-1})))$, on conclut alors en
choisissant une suite judicieuse de telles fonctions (c'est le fait que $\bar w_2$ soit à décroissante
parabolique qui nous permet de trouver une telle suite).

On sait que $\cH^\bullet(C_1(\bS^{d-1}))$ est de dimension un, l'hypothèse que
l'image de $d$ est presque fermée sur $X_1$ implique que
l'image de $$d\,:\,  \cC^\bullet_{w_1,\mu_1}(X_1)\otimes \cH^\bullet(C_1(\bS^{d-1}))
\rightarrow  L_{\mu_1}^2(\Lambda^\bullet T^* X_1)\otimes \cH^\bullet(C_1(\bS^{d-1}))$$
est aussi fermée.
En particulier si $H$ est le projecteur orthogonal de $L^2_\mu(\Lambda T^* X)$ sur
$\cH_\mu(X)=\cH_{\mu_1}(X_1)\otimes \cH^\bullet(C_1(\bS^{d-1}))$ alors on trouve 
$$\beta_1\in \cH_{\mu_1}(X_1)\otimes \cC_{w_2,\mu_2}(X_2)\subset \cC_{w_2,\mu}(X)$$
tel que $\bar h(\alpha)=H(\alpha)+d\beta_1$, 
si alors on pose $\beta_2=\bar B\alpha$, on obtient bien le résultat annoncé.

La preuve pour le cas où $d\ge 2$ et concernant la cohomologie $L^2_\mu$ relative de $X$ est identique.
Il reste à démontrer le résultat pour le cas $d=2$ et pour la cohomologie $L^2_\mu$  absolue, et cela est une
conséquence des propositions (\ref{exbord1}) et  (\ref{kunneth1}).

\endproof
On peut évidement adapter cette preuve lorsqu'on considère d'autre poids sur $X=X_1\times X_2$, et une preuve
identique montre que 
\begin{prop}\label{kunneth3}On suppose que l'image de $d\,:\, \cC_{w_1,\mu_1}(X_1)\rightarrow L_{\mu_1}^2(\Lambda^\bullet T^* X_1)$
est fermée. On considère $X=X_1\times C_1(\bS^{d-1})$ équipé de la métrique riemannienne produit et de la
mesure $\mu(x_1,x_2)=\mu_1(x_1)$, et on pose $\bar w_2(x_1,x_2)=|x_ 2|^{-2}$ 
et $ w_2(x_1,x_2)=|x_ 2|^{-2}(1+\log|x_ 2|)^{-2}$ et on définit
 $Y=X_ 2\times \partial C_1(\bS^{d-1})=X_ 2\times \bS^{d-1}\subset \partial X$.
  Soit $Z=\emptyset$ ou $Z=Y$.
\begin{enumeroman}
\item Lorsque pour $d>2$, on considère sur $X_1\times X_2$, le poids 
$$\mu(x_1,x_2)=\mu_1(x_1)(1+\log|x_2|)^{-2},$$ 
alors pour $\alpha\in B^k_{\mu}(X,Z)$, 
il y a $\beta_2\in \cC^{k-1}_{\bar w_2,\mu}(X,Z) $ et
 $\beta_1\in \cC^{k-1}_{ w_1,\mu}( X,Z)$  
tels que $$\alpha=d\beta_1+d\beta_2,$$ de plus
$$\H_\mu^k(X)=\H^{k-d+1}_{\mu_1}(X_1)\ \et \  
\H_\mu^k(X,Y)=\H^{k-1}_{\mu_1}(X_1). $$

\item Lorsque pour $d>4$, on considère sur $X_1\times X_2$, le poids 
$$\mu(x_1,x_2)=\mu_1(x_1) |x_2|^{-2}(1+\log|x_2|)^{-2}$$ ou
$$\mu(x_1,x_2)=\mu_1(x_1) |x_2|^{-2}(1+\log|x_2|)^{-4}$$ 
alors pour $\alpha\in B^k_{\mu}(X,Z)$, 
il y a $\beta_2\in \cC^{k-1}_{\bar w_2,\mu}(X,Z) $ et
 $\beta_1\in \cC^{k-1}_{ w_1,\mu}( X,Z)$  
tels que $$\alpha=d\beta_1+d\beta_2,$$ de plus
$$\H_\mu^k(X)=\H^{k-d+1}_{\mu_1}(X_1)\ \et\ 
 \H_\mu^k(X,Y)=\H^{k-1}_{\mu_1}(X_1). $$

\item Lorsque pour $d\le 4$, on considère sur $X_1\times X_2$, le poids 
$$\mu(x_1,x_2)=\mu_1(x_1) |x_2|^{-2}(1+\log|x_2|)^{-2}$$ ou
$$\mu(x_1,x_2)=\mu_1(x_1) |x_2|^{-2}(1+\log|x_2|)^{-4}$$ 
alors 
$$\H_\mu^k(X,Y)=\{0\} $$
et pour $d>2$ alors 
$$\H_\mu^k(X)=\H^{k-d+1}_{\mu_1}(X_1) $$
par contre pour $d=2$ alors : 
$$\H_\mu^k(X)=\H^{k}_{\mu_1}(X_1)\oplus\H^{k-1}_{\mu_1}(X_1) $$
 Si de plus $d\not =4$ alors
pour  $\alpha\in B^k_{\mu}(X,Z)$, 
il y a $\beta_2\in \cC^{k-1}_{\bar w_2,\mu}(X,Z) $ et
 $\beta_1\in \cC^{k-1}_{w_1,\mu}( X,Z)$  
tels que $$\alpha=d\beta_1+d\beta_2.$$ 
Mais lorsque $d=4$, alors pour $\alpha\in B^k_{\mu}(X,Z)$, 
il y a $\beta_2\in \cC^{k-1}_{w_2,\mu}(X,Z) $ et
 $\beta_1\in \cC^{k-1}_{ w_1,\mu}( X,Z)$ 
 telle que $\alpha=d\beta_1+d\beta_2.$
 
\item Lorsque pour $d>2$, on considère sur $X_1\times X_2$, le poids 
$$\mu(x_1,x_2)=\mu_1(x_1) |x_2|^{2}(1+\log|x_2|)^{2}$$ ou
$$\mu(x_1,x_2)=\mu_1(x_1) |x_2|^{2}$$ 
alors pour $\alpha\in B^k_{\mu}(X,Y)$, 
il y a $\beta_2\in \cC^{k-1}_{\bar w_2,\mu}(X,Y) $ et
 $\beta_1\in \cC^{k-1}_{ w_1,\mu}( X,Y)$  
tels que $$\alpha=d\beta_1+d\beta_2.$$
De plus
$$\H_\mu^k(X,Y)=\H^{k-1}_{\mu_1}(X_1). $$
\item Lorsque pour $d=2$, on considère sur $X_1\times X_2$, le poids 
$$\mu(x_1,x_2)=\mu_1(x_1) |x_2|^{2}$$ 
alors pour $\alpha\in B^k_{\mu}(X,Y)$, 
il y a $\beta_2\in \cC^{k-1}_{ w_2,\mu}(X,Y) $ et
 $\beta_1\in \cC^{k-1}_{ w_1,\mu}( X,Y)$  
tels que $$\alpha=d\beta_1+d\beta_2.$$
De plus
$$\H_\mu^k(X,Y)=\H^{k-1}_{\mu_1}(X_1). $$
\item Lorsque pour $d=2$, on considère sur $X_1\times X_2$, le poids 
$$\mu(x_1,x_2)=\mu_1(x_1) |x_2|^{2}(1+\log|x_2|)^{2}$$ 
alors pour $\alpha\in B^k_{\mu}(X,Y)$, 
il y a $\beta_2\in \cC^{k-1}_{ w_2,\mu}(X,Y) $ et
 $\beta_1\in \cC^{k-1}_{ w_1,\mu}( X,Y)$  
tels que $$\alpha=d\beta_1+d\beta_2.$$
De plus
$$\H_\mu^k(X,Y)=\H^{k-1}_{\mu_1}(X_1)\oplus \H^{k-2}_{\mu_1}(X_1) . $$
\end{enumeroman}
\end{prop}

\section{La géométrie des variétés QALE.}
Dans le livre \cite{joyce}, D. Joyce construit non seulement de nouvelles variétés compactes à
 holonomie exceptionnelle, mais aussi il construit de nouveaux exemples de variétés complètes
non compactes Kähler-Einstein à courbure scalaire nulle sur certaines résolutions de $ \C^n/G$.
 Nous
intéressons ici aux espaces de formes harmoniques $L^2$ sur ces variétés, puisque la dimension 
de ces espaces
est un invariant de quasi isométrie, nous nous contentons de décrire la géométrie de ces variétés 
à quasi-isométries près et au dehors d'un compact. La construction part d'un sous-groupe fini de 
$G\subset SU(n)$
et d'une résolution localement en produit $\pi\,:\, X\rightarrow \C^n/G$ de la variété singulière. 
Rappelons d'abord ce qu'est une telle résolution :
\subsection{Résolution localement en produit}
On note $V=\C^n$ et $I$ une indexation des sous-espaces vectoriels de $\C^n$ qui sont fixés par
 un sous groupe de $G$
$$\{V_i,i\in I\}=\{V^H, H {\rm\ sous-groupe\ de\ } G\}.$$
On note aussi
$$A_i=C(V_i)=\{g\in G, gv=v\  \forall v\in V_i\}, $$
$$ N(V_i)=\{g\in G,\ g(V_i)=V_i\},\ \et B_i=N(V_i)/A_i.$$
On suppose que $0,\infty\in I$ et que
$V_0=V$ et $V_\infty = \{0\}=V^G$. On met sur $I$ l'ordre partiel suivant :
$$i\le j\Leftrightarrow V_i\supset V_j $$
Si on note $W_i=V_i^\perp$ alors $$i\le j\Leftrightarrow W_i\subset W_j .$$
On a donc $V/A_i\simeq \left(W_i/A_i\right)\times V_i$. On note aussi
$$\tilde S=\left(\cup_{i\in I\setminus\{0\}} V_i\right)/G$$ le lieu singulier de $V/G$.
Et pour $i\in I$ on définit 
$$S_i=\left(\bigcup_{j\in I\setminus\{0\}, i\not\ge j} V_j\right)/A_i
=\left(\bigcup_{j\in I\setminus\{0\}, V_i\not\subset V_j} V_j\right)/A_i.$$
et $$T_i=\{x\in V/A_i, d(A_ix,S_i)\le R\}.$$
\begin{defi} $\pi\,:\, X\rightarrow \C^n/G$ est une résolution localement produit si pour
 tout $i\in I$, il
y a une résolution $\pi_i\,:\, Y_i\rightarrow W_i/A_i$ (qui sera forcément aussi 
localement en produit) 
de $W_i/A_i$ telle que si on note $$U_i=( \pi_i\times\Id)^{-1}(T_i)\subset Y_i\times V_i$$ 
alors il
y a un biholomorphisme local $\psi_i\,:\, (Y_i\times V_i)\setminus U_i\longrightarrow X$
 qui rendent le
diagramme suivant commutatif :
 $$
\xymatrix{
  {(Y_i\times V_i)\setminus U_i} \ar[d]^{\pi_i\times\Id} \ar[r]^{\ \ \ \ \psi_i} & {X}\ar[d]^{\pi}
   \\
  {(W_i/A_i)\times V_i\setminus T_i} \ar[r]^{\ \ \ \ \ \ \ \ \ \phi_i}& {V/G}}
$$

 où $\phi_i$ est l'application naturelle
  $\phi_i\,:\,\left(W_i/A_i\right)\times V_i\longrightarrow V/G$.
\end{defi}

Le groupe $G$ agit naturellement sur $I$, D. Joyce montre alors que lorsque $(X,\pi)$
 est une résolution localement
en produit de $V/G$ alors il existe un biholomorphisme $g\,:\, Y_i\longrightarrow Y_{g.i}$
 qui rend le diagramme suivant commutatif 
$$\xymatrix{
{Y_i}\ar[d]^{\pi_i}\ar[r]^{g} &{Y_{g.i}}\ar[d]^{\pi_{g.i}}\\
{W_i/A_i}\ar[r]^{g}&{W_{g.i}/A_{g.i}}
} $$

On peut décrire ce qui se passe en dimension $n=3$. Soit donc $G\subset SU(3)$
 un sous-groupe fini, si $g\in
G\setminus\{\Id\}$ alors on a forcément $\dim\ker(g-\Id)\le 1$. En particulier,
 si $I=\{0,1,...N,\infty\}$
alors on a $V_i\cap V_j=\{0\}$ si $i\not =j$ et $i,j\in \{1,...N\}$.
 Ceci signifie exactement que $\bS^5/G$
n'a que des singularités isolées.
Notons $S=\bigcup_{i=1}^N V_i$ et $S^\eps$ le $\eps$-voisinage de $S$.
La topologie à l'infini de $X$ est le produit d'une résolution de $\bS^5/G$ avec une demi-droite.
Cette résolution de $\bS^5/G$ peut être décrite topologiquement de la façon suivante :
notons $B_i(\eps)$ la boule de rayon $\eps$ dans $W_i$ et $Y_i(\eps)=\pi_i^{-1}( B_i(\eps))$
Sur $\bS^5\setminus S^\eps$,
 l'action de $G$ est
libre (pourvu que $\eps>0$ soit choisi assez petit) et le quotient est difféomorphe à 
$\Sigma_0=(\bS^5\setminus S^\eps)/G$

Supposons que l'action de $G$ sur $\{1,...,N\}$ ait $l$ orbites $GV_1,...,GV_l$.
Alors le bord de  $(\bS^5\setminus S^\eps)/G$ a exactement $l$ composantes connexes chacune 
difféomorphe à un 
\begin{equation*}
\begin{split}
G.\big(\partial B_i(\eps)\times (V_i\cap \bS^5)\big)/G&=
\big(\partial B_i(\eps)\times (V_i\cap \bS^5)\big)/N(V_i)\\
&=\big(\left( \partial B_i(\eps)/A_i \right) \times (V_i\cap \bS^5))/B_i.\\
\end{split}\end{equation*}
 Pour $i=1,...,l$, nous notons 
 \begin{equation*}
\begin{split}
\Sigma_i&=G.\big( Y_i(\eps)\times (V_i\cap \bS^5)\big)/G\\
&= \left( Y_i(\eps)\times (V_i\cap \bS^5) \right)/B_i.\\
 \end{split}\end{equation*}
 Alors cette résolution de $\bS^5/G$  est difféomorphe à 
 $$\Sigma_0\# (\cup_{i\in\{1,...,l\}} \Sigma_i)$$ 
 Ou on a recollé chaque $\big(\partial B_i(\eps)\times (V_i\cap \bS^5))/B_i$ 
avec $$\big( \partial Y_i(\eps)\times (V_i\cap \bS^5))/B_i \subset \big(Y_i(\eps)\times (V_i\cap \bS^5))/B_i .$$

\subsection{ Métrique Quasi-Asymptotiquement-Localement-Euclidienne (QALE)}
Une métrique Quasi-Asymptotiquement-Localement-Euclidienne (QALE pour faire court) asymptote à
 $\C^n/G$ sur
une résolution localement en produit
$\pi\,:\, X\longrightarrow\C^n/G$ est une métrique $g$ sur $X$ 
telle que dans le cadre de la
définition précédente $\psi_i^*g $ soit asymptote à la métrique produit $g_i\oplus h_i$, 
où $h_i$ est une métrique 
QALE sur $Y_i$ et $h_i$ la métrique euclidienne sur $V_i$.

La vitesse à laquelle $\psi_i^*g$ tend vers la métrique produit dépend des 
distances à $S_i$ et à l'origine.
Puisque nous intéressons qu'à la classe de quasi isométrie de ces variétés,
 nous ne retiendrons que
$\psi_i^*g $ est la métrique produit  $g_i\oplus h_i$.
 
\subsection{Cas des groupes $G$ de longueur $2$}
Nous intéressons en fait à la classe la plus simple des variétés QALE (après
celle où $G$ agit sans point fixe sur $\bS^{2n-1}$, la variété étant alors ALE).
\begin{defi}Soit $G$ un sous groupe de $SU(n)$, on appelle longueur de $G$,
$$L(\C^n,G)=\max\{k, {\rm\ il\ y\ a\ } i_0<i_1...<i_k, i_j\in I\}.$$
\end{defi}

Ainsi la longueur de $G$ vaut $1$ si et seulement si $S=\{0\}$ i.e. si et seulement si 
$G$ agit sans point fixe sur $\bS^{2n-1}$.
Et la longueur de $G$ vaut $2$ si et seulement si les singularités de 
$\bS^{2n-1}/G$ sont isolées ; par exemple lorsque
$G\subset \Su(3)$ ou $G\subset \Sp(2)$, on a forcément $L(\C^3,G)\le 2$.

Lorsque la longueur de $G$ vaut $1$, une métrique ALE sur une résolution de $\C^n/G$ 
est au dehors d'un compact quasi isométrique à une métrique plate sur $(\C^n\setminus B(R))/G$.

On peut aussi décrire la géométrie de variétés QALE asymptote à $\C^n/G$ 
lorsque la longueur de $G$ est $2$:

On numérote $G.V_1,...,G.V_l$, les différentes orbites de $I\setminus \{0,\infty\}$
 sous l'action de $G$. Le sous groupe
$A_i$ agit sans point fixe sur $W_i\setminus \{0\}$ et notons $\pi_i\,:\,Y_i\rightarrow W_i/A_i$ une
résolution localement en produit de $W_i/A_i$.
Si $y\in Y_i$ on note $|y|=\|\pi_i(y)\|$ et $\bB$ la boule unité de $\C^n$.
Chaque $Y_i$ est donc muni une métrique $g_i$ plate au dehors d'un compact. 
On équipe maintenant $Y_i\times (V_i\setminus \bB)$ de la métrique produit et 
on définit 
 $$C_i=\{(y,v)\in Y_i\times (V_i\setminus \bB), |y|\le 2\eps |v|\}$$
Lorsque $g\in G$ on a donc une application $g\,:\, C_i\rightarrow C_{g.i}$, notons
$E_i$ le quotient $G.C_i/G$, on a aussi $E_i=C_i/B_i$.
Soit maintenant 
$$E_0=\Big( \C^n\setminus 
\big( \bB\cup \left\{v\in \C^n, \exists i\   d(v,V_i)< \eps d(v, W_i)\,\right\}\big) \Big)/G$$

Si on note $S_i^\eps =\{v\in \C^n, d(v,V_i)\le \eps\}$ alors
$E_0$ est un bout de cône sur $\big(\bS^{2n-1}\setminus (\bigcup_{i>1} S_i^\eps) )/G$.
Pour $\eps$ assez petit les $S_i^{2\eps}$ sont deux à deux disjoints et
 à quasi isométrie près nous identifions :
$$X\setminus \pi^{-1}(\bB)=\cup_{i=0}^l E_i$$ où
$$E_i\cap E_j=\emptyset,\  {\rm lorsque } \ i\not=j \et i,j\not=0$$
$$\et  E_i\cap E_0=C_1\big( \left((S_i^{2\eps}\cap \bS^{2n-1})\setminus S_i^{\eps} \right)/G)$$

\section{$L^2$ cohomologie des variétés QALE de rang 2.}
Nous allons maintenant déterminer la cohomologie $L^2$ des variétés QALE asymptotes à $\C^n/G$ où
$G\subset\Su (n)$ est tel que le quotient $\bS^{2n-1}/G$ n'a que des singularités isolées. Soit $X$ une
telle variété, rappelons qu'au dehors d'un compact $K\subset X$, $\cO=X\setminus K$ est la réunion de
$E_0,E_1,...,E_l$. Où pour $\eps>0$ assez petit :

$E_0$ est un bout de cône sur $(\bS^{2n-1}\setminus S^\eps)/G $ avec
$$S=\{v\in \bS^{2n-1},\ \exists g\in G\setminus\{\Id\}, gv=v\}$$ et
 $$S^\eps=\{v\in \bS^{2n-1}, d(v,S)<\eps\}.$$

L'ensemble $S$ est en fait une réunion de sous espace vectoriel de $\C^n$ ne s'intersectant qu'en $\{0\}$ et
$G$ agit sur cette famille de sous-espaces vectoriels. On suppose
$$S=\cup_{i=1}^l G.V_i$$ avec 
$$G.V_i\cap G.V_j=\{0\}\ \  {\rm si\  }\ i\not=j.$$
Notons $W_i=V_i^{\perp}$, $A_i=\{g\in G, \forall v\in V_i, gv=v\}$, $
N(V_i)=\{g\in G, g(V_i)=V_i\}$ et $B_i=N(V_i)/A_i$.
Soit $\pi_i\,:\,Y_i\rightarrow W_i/A_i$ une résolution de  $W_i/A_i$
 équipée d'une métrique ALE asymptote à $W_i/A_i$ alors on suppose que $B_i$ agit sur
$Y_i\times V_i$ et $E_i=\widehat E_i/B_i$ où
$$\widehat E_i=\{(y,v)\in Y_i\times V_i,\ \eps<|v|,\  |\pi_i(y)|<2\eps |v|\}.$$
alors $E_i\cap E_0$ est quasi isométrique à un bout de cône sur 
$\left(\bS_{W_i}/A_i\times \bS_{V_i}\times ]\eps,2\eps[\right)/B_i$
où on note $\bS_V=\bS^{2n-1}\cap V$ la sphère unité d'un sous espace vectoriel $V\subset \C^n$.

\subsection{La cohomologie de $X\setminus K$} On commence par calculer la cohomologie de $\cO=X\setminus
K$. On sait que $\cO$ se rétracte sur $\partial K$ et que
$\partial K=\cup_{i=0}^l \Sigma_i$ où
$\Sigma_0=(\bS^{2n-1}\setminus S^\eps)/G$ et 
$\Sigma_i=\{(y,v)\in Y_i\times\bS_{V_i}, |y|<2\eps\}/B_i$.
Compte tenu du fait que la cohomologie de la sphère $\bS^{2n-1}$ est nulle en degré $k\not=0,2n-1$, on
a pour les degrés $k\in ]0,2n-2[$ :
\begin{equation}
\label{coho1}
H^k(\Sigma_0)\oplus H^k(S)^{G}\simeq H^k(\partial \Sigma_0) .
\end{equation}
Et compte tenu du fait que $H^{2n-2}(S)=\{0\}$, nous obtenons aussi :
\begin{equation}
\label{coho1p}
\{0\}\rightarrow H^{2n-2}(\Sigma_0)\rightarrow  H^{2n-2}(\partial \Sigma_0)\simeq \R^{l}\rightarrow
H^{2n-1}(\bS^{2n-1})^G\simeq \R\rightarrow \{0\}.
\end{equation}
où si un groupe $G$ agit linéairement sur un espace vectoriel $H$ on a noté $H^G$, 
le sous espaces des vecteurs $G$-invariants. On utilisera la suite exacte de Mayer-Vietoris :
 \begin{equation}
\label{coho2}
\bigoplus_{i=0}^l H^{k-1}(\Sigma_i)\rightarrow H^{k-1}(\partial \Sigma_0)
\rightarrow H^k(\partial K)\rightarrow
\bigoplus_{i=0}^l H^{k}(\Sigma_i)\rightarrow H^{k}(\partial \Sigma_0).\end{equation}
Notons $n_i=\dim V_i$. Remarquons que $B_i$ agit trivialement sur la cohomologie de $\bS_{V_i}$ ainsi 
$$H^k(S)^{G}=\bigoplus_{i=1}^l H^k(\bS_{V_i})$$ et 
 $$H^k(\Sigma_i)=H^k(Y_i)^{B_i}\oplus H^{k-2n_i+1}(Y_i)^{B_i} .$$
 En particulier, l'application $H^k(\Sigma_i)\rightarrow H^k(\bS_{V_i})$ est surjective. Avec (\ref{coho1}),
 on en déduit que pour $k\not=0,2n-2,2n-1$, l'application
 $$\bigoplus_{i=0}^l H^{k}(\Sigma_i)\rightarrow H^{k}(\partial \Sigma_0)$$ est surjective.
Et donc pour $k\not= 0,1,2n-2,2n-1$, la suite exacte (\ref{coho2}) implique que 
  $$\{0\} \rightarrow  H^k(\partial K)\rightarrow
\bigoplus_{i=0}^l H^{k}(\Sigma_i)\rightarrow H^{k}(\partial \Sigma_0)\simeq H^k(\Sigma_0)\bigoplus
H^k(S)^{G}\rightarrow \{0\} ; $$
ainsi pour ces degrés $H^k(\partial K)$ est isomorphe à
\begin{equation*}
\bigoplus_{i=1}^l \ker \big(H^{k}(\Sigma_i)\rightarrow H^k(\bS_{V_i})^{B_i})
=\bigoplus_{i, k>2n_i-1} H^{k-2n_i+1}(Y_i)^{B_i}\oplus \bigoplus_{i} H^{k}(Y_i)^{B_i}.\end{equation*}
Pour les degrés $k=1, 2n-2$ : on obtient notamment grâce à (\ref{coho2}) :
$$H^1(\partial K)\simeq \bigoplus_{i\ge 1} H^1(Y_i)^{B_i}$$
Mais, nous avons toujours $2n_i-1\ge 1$, en conséquence il est aussi vrai que :
$$H^1(\partial K)\simeq \bigoplus_{i, 1>2n_i-1} H^{1-2n_i+1}(Y_i)^{B_i}\oplus \bigoplus_{i}
H^{1}(Y_i)^{B_i}.$$
Nous obtenons aussi avec (\ref{coho1p}) :
$$H^{2n-2}(\partial K)\simeq \bigoplus H^{2\dim Y_i-1}(Y_i)^{B_i}.$$
Mais puisqu'on a toujours $4\le \dim Y_i=2n-2n_i\le 2n-2$, il est aussi vrai que
$$H^{2n-2}(\partial K)\simeq \bigoplus_{i, 2n-2>2n_i-1} H^{2n-2-2n_i+1}(Y_i)^{B_i}\oplus \bigoplus_{i}
H^{2n-2}(Y_i)^{B_i}.$$
 On obtient donc 
 \begin{prop}\label{cohopK}
 Pour $k\not\in \{0,2n-1\}$
 $$H^k(\cO)=H^k(\partial K)=
 \bigoplus_{i, k>2n_i-1} H^{k-2n_i+1}(Y_i)^{B_i}\oplus \bigoplus_{i} H^{k}(Y_i)^{B_i}.$$
 \end{prop}

 \subsection{Cohomologie $L^2$ de $\cO$.}Nous considérons sur $\cO$ les poids
 $\bar w(x)=|x|^{-2}$ et $w(x)=|x|^{-2}(1+|\log x|)^{-2}$
 où si $\pi \,:\, \cO\rightarrow \C^n/G$, on note $|x|=\|\pi(x)\|$.
 On sait que $$\cO=\bigcup_{i\ge 0} E_i$$ et que
 $E_0$, $E_0\cap E_i$ sont des bouts de cônes on sait que l'image de $d$ y 
 est presque fermée en tout degré et par rapport au poids $w$, 
 de plus on connaît leurs cohomologies $L^2$.
  On doit donc maintenant calculer la cohomologie $L^2$ de $E_i$.
 \subsubsection{Cohomologie $L^2$ de $E_i$}
   On sait que $E_i=\widehat E_i/B_i$ et donc :
  $$\H^k(E_i)= \H^k(\widehat E_i)^{B_i}.$$

  $\widehat E_i$ est une partie de $Y_i\times (V_i\setminus \bB)$ (on note $\bB$ la boule unité de 
  $\C^n$). Sur $\widehat E_i$,  $\bar w$ est comparable 
  $\bar w'(y,v)=\|v\|^{-2}$ et $w$ à $w'(y,v)=\|v\|^{-2}(1+|\log \|v\| )^{-2}$. 
  Notons également $\bar w_1(y,v)=|y|^{-2}$.
  
  Nous commençons par remarquer que (\ref{kunneth2}) :
  $$\H^k(Y_i\times(V_i\setminus \bB) )\simeq \bigoplus_{p+q=k}\H^p(Y_i)\otimes \H^q(V_i\setminus \bB)$$
  De plus si $Y_i^{\eps}= \{y\in Y_i, |y|<\eps\} $ alors puisque  l'application 
  $ \H^p(Y_i)\rightarrow \H^p(Y_i^{\eps})=H^p(Y_i^{\eps})$ est injective (cf. \ref{injALE}), il en
  est de même de l'application
  $$\H^k\big(Y_i\times(V_i\setminus \bB) \big)
  \rightarrow\H^k\left(Y_i^{\eps}\times(V_i\setminus \bB)\right).$$
  Ainsi puisque $Y_i^{\eps}\times(V_i\setminus \bB)\subset \widehat E_i\subset Y_i\times(V_i\setminus \bB)$,
  l'application de restriction :
  $$\iota^*_{\widehat E_i}\,:\,\H^k(Y_i\times(V_i\setminus \bB) )\rightarrow\H^k(\widehat E_i )$$
  est aussi injective.
  
  On suppose maintenant que $k<n$. Soit $\alpha\in B^k(\widehat E_i)$. Posons
  $$\Omega_i=C_1(]\eps/2,2\eps[\times \bS^{2m_i-1}\times \bS^{2n_i-1})\subset
  \widehat E_i$$
  où $n_i=\dim_\C V_i$ et $m_i=n-n_i$. Puisque 
  $Y_i\times(V_i\setminus \bB)=\widehat E_i\cup V_i$ avec $V_i\cap \widehat
  E_i=\Omega_i$, nous avons 
  $$\iota^*_{\Omega_i}\alpha\in B^k(\Omega_i)=d\cC^{k-1}_{\bar w,1}(\Omega_i),$$
  il y a donc $\psi\in \cC^{k-1}_{\bar w,1}(\Omega_i)$ tel que 
  $$\iota^*_{\Omega_i}\alpha=d\psi$$
  Maintenant grâce à (\ref{condition}), nous savons qu'il existe $\bar \psi\in \cC^{k-1}_{\bar w,1}(\widehat E_i)$
  tel que $\iota^*_{\Omega_i}\bar \psi=\psi$. En particulier, 
  l'extension par zéro de $\alpha-d\bar \psi$ à $Y_i\times(V_i\setminus \bB)$ est fermée :
  $$\alpha-d\bar \psi\in Z^k(Y_i\times(V_i\setminus \bB))$$
  Mais l'image de la classe de cohomologie $L^2$ de cette forme
   par $\iota^*_{\widehat E_i}$ est la classe de cohomologie de $\alpha\in B^k(\widehat E_i)$, puisque cette
   application est injective, la proposition (\ref{kunneth2}) nous permet de trouver 
   $u\in \cC^{k-1}_{\bar w_1,1}(Y_i\times(V_i\setminus \bB))$ et
    $v\in \cC^{k-1}_{w',1}(Y_i\times(V_i\setminus \bB))$
  tels que
  $$\alpha-d\bar\psi=du+dv$$
  Compte tenu du fait que sur $\widehat E_i$ : $\bar w_1\ge C w'$ et que $w'$ et comparable à $w$, alors
  $$\phi= \bar\psi+\iota^*_{\widehat E_i}(u+v)\in \cC^{k-1}_{w,1}(\widehat E_i)\ \et \ \alpha=d\phi.$$
  Si nous partons d'un $\alpha\in Z^k(\widehat E_i)$ alors puisque $Z^k(\Omega_i)=B^k(\Omega_i)$ lorsque 
  $k<n$, alors nous obtenons qu'en fait l'application :
  $$\iota^*_{\widehat E_i}\,:\,\H^k(Y_i\times(V_i\setminus \bB) )\rightarrow\H^k(\widehat E_i )$$
  est un isomorphisme. Remarquons également que puisque 
  l'application d'extension par zéro 
  $$ H_c^p(Y_i)\simeq \H^p(Y_i^\eps, \partial Y_i^\eps) \rightarrow \H^p(Y_i)$$ est surjective 
  (cf. \ref{injALE}), il en est de même de l'application :
  $$\H^k(Y_i^\eps\times (V_i\setminus \bB), \partial Y_i^\eps\times (V_i\setminus \bB))\rightarrow 
  \H^k(\widehat E_i ).$$
  
  Nous avons donc démontré :
  \begin{prop}
  Si $k<n$, alors $d\cC_{w,1}^{k-1}(\widehat E_i)=B^k(\widehat E_i)=$ de plus 
  $$\H^k(\widehat E_i )\simeq \oplus_{p+q=k} \H^p(Y_i)\otimes \H^q(V_i\setminus \bB)=
  \left\{\begin{array}{ll} \H^{k-2n_i+1}(Y_i)&{\rm\ si\ }n_i>1\\
  \{0\}&{\rm\ si\ }n_i=1.\\ \end{array}\right.$$
  Et l'application d'extension par zéro 
  $$\H^k(Y_i^\eps\times (V_i\setminus \bB), \partial Y_i^\eps\times (V_i\setminus \bB))\rightarrow 
  \H^k(\widehat E_i )$$ est surjective.
  \end{prop} En considérant la cohomologie relative à $Y_i\times \bS^{2n_i-1}$, nous obtenons également :
  \begin{prop}\label{relatifL2}
  Si $k>n$, alors $d\cC_{w,1}^{k-1}(\widehat E_i,Y_i\times \bS^{2n_i-1})=B^k(\widehat E_i,Y_i\times \bS^{2n_i-1})$ de plus 
  $$\H^k(\widehat E_i,Y_i\times \bS^{2n_i-1} )\simeq 
  \oplus_{p+q=k} \H^p(Y_i)\otimes \H^q(V_i\setminus \bB,\bS^{2n_i-1})=
  \left\{\begin{array}{ll} \H^{k-1}(Y_i)&{\rm\ si\ }n_i>1\\
  \{0\}&{\rm\ si\ }n_i=1\\ \end{array}\right.$$
  Et l'application d'extension par zéro 
  $$\H^k(Y_i^\eps\times (V_i\setminus \bB), \partial Y_i^\eps\times (V_i\setminus \bB)\cup Y_i\times \bS^{2n_i-1} )\rightarrow 
  \H^k(\widehat E_i,Y_i\times \bS^{2n_i-1}) )$$ est surjective.
  \end{prop}
  
 Il reste maintenant à déterminer la cohomologie $L^2$ de $\widehat E_i$ en degré $n$. Le raffinement
 présenté en (\ref{raffinement}) montre que cette même preuve est valide lorsque $b_{n-1}(\bS^{2m_i-1}\times
 \bS^{2n_i-1})=0$, c'est à dire lorsque $n_i\not=m_i$. Et même mieux si de plus $n_i>1$ alors, on trouve que
 $d\cC_{\bar w,1}^{n-1}(\widehat E_i,Y_i\times \bS^{2n_i-1})=B^n(\widehat E_i,Y_i\times \bS^{2n_i-1}).$
  
   Supposons donc que $n_i=m_i=n/2$. Alors nous avons forcément $n_i\ge 2$. Notons :
   $\rho=(\log(\|v\|+1)^{-2}$. On reprend nos arguments si 
   $\alpha\in Z^n(\widehat E_i, Y_i\times \bS^{2n_i-1})$, alors nous trouvons comme précédemment
   $\psi\in \cC^{n-1}_{ w,1}(\Omega_i, Y_i\times \bS^{2n_i-1})$ tel que 
  $$\iota^*_{\Omega_i}\alpha=d\psi$$
  Maintenant grâce à (\ref{condition}), nous savons qu'il existe 
  $\bar \psi\in \cC^{n-1}_{ w',\rho}(\widehat E_i, Y_i\times \bS^{2n_i-1})$
  tel que $\iota^*_{\Omega_i}\bar \psi=\psi$. En particulier, 
  l'extension  par zéro de $\alpha-d\bar \psi$ à $Y_i\times(V_i\setminus \bB)$est fermée :
  $$\alpha-d\bar \psi\in Z_\rho^n(Y_i\times(V_i\setminus \bB, Y_i\times \bS^{2n_i-1}).$$
  Les mêmes causes produisant les mêmes effets, il est également vrai que l'application :
  $$\iota^*_{\widehat E_i}\,:\,\H_\rho^n(Y_i\times(V_i\setminus \bB), Y_i\times \bS^{2n_i-1} )
  \rightarrow\H_\rho^n(\widehat E_i , Y_i\times \bS^{2n_i-1})$$
  est injective.
  Et puisque $n_i>1$, d'après (\ref{kunneth3}), l'application naturelle :
  $$\H^n(Y_i\times(V_i\setminus \bB), Y_i\times \bS^{2n_i-1} )\rightarrow 
  \H^n_\rho(Y_i\times(V_i\setminus \bB), Y_i\times \bS^{2n_i-1} )$$ est un isomorphisme.
  
  Mais l'image de la classe de cohomologie $L_\rho^2$ de $\alpha-d\bar \psi$
   par $\iota^*_{\widehat E_i}$ est la classe de cohomologie de $\alpha\in Z_{\rho}^k(\widehat E_i)$, 
   la proposition (\ref{kunneth3}) nous permet de trouver 
   $u\in \cC^{k-1}_{\bar w_1,\rho}(Y_i\times(V_i\setminus \bB), Y_i\times \bS^{2n_i-1})$, 
    $v\in \cC^{k-1}_{\bar w,\rho}(Y_i\times(V_i\setminus \bB), Y_i\times \bS^{2n_i-1})$ et 
    $h\in \cH^n(Y_i\times(V_i\setminus \bB), Y_i\times \bS^{2n_i-1} )$
  tels que
  $$\alpha-d\bar\psi=h+du+dv$$
  Si on pose $\phi=\psi+\iota^*_{\widehat E_i}(u+v)$ alors
  $\alpha=\iota^*_{\widehat E_i} h+d\phi$ et $\phi\in L^2_{w}(\Lambda^{n-1}T^*(E_i))$.
  Nous avons ainsi montrer que la proposition (\ref{relatifL2}) est encore valide pour $k=n$.
  
  En revenant à $E_i$, nous obtenons
  \begin{prop}\label{cohoL2Ei}\begin{enumeroman}
  \item Si $k<n$, alors $d\cC^{k-1}_{w,1}( E_i)=B^k( E_i)$ de plus 
  $$\H^k( E_i )\simeq \oplus_{p+q=k} \H^p(Y_i)^{B_i}\otimes \H^q(V_i\setminus \bB)^{B_i}=
  \left\{\begin{array}{ll} \H^{k-2n_i+1}(Y_i)^{B_i}&{\rm\ si\ }n_i>1\\
  \{0\}&{\rm\ si\ }n_i=1\\ \end{array}\right.$$
  Et l'application d'extension par zéro 
  $$\H^k((Y_i^\eps\times (V_i\setminus \bB))/B_i, (\partial Y_i^\eps\times (V_i\setminus \bB))/B_i)\rightarrow 
  \H^k( E_i )$$ est surjective.
  \item Si $k\ge n$, alors $d\cC^{k-1}_{w,1}( E_i,\partial\cO)=B^k( E_i,\partial\cO)$ de plus 
  $$\H^k( E_i,\partial\cO )\simeq \oplus_{p+q=k} \H^p(Y_i)^{B_i}\otimes \H^q(V_i\setminus \bB,\bS^{2n_i-1})^{B_i}=
  \left\{\begin{array}{ll} \H^{k-1}(Y_i)^{B_i}&{\rm\ si\ }n_i>1\\
  \{0\}&{\rm\ si\ }n_i=1\\ \end{array}\right.$$
  \noindent Et l'application d'extension par zéro 
  $$\H^k( (Y_i^\eps\times (V_i\setminus \bB))/B_i
  ,\partial\cO\cup (\partial Y_i^\eps\times (V_i\setminus \bB))/B_i)
  \rightarrow \H^k( E_i )$$ est surjective.
  \end{enumeroman}\end{prop}
  
    \subsubsection{Cohomologie $L^2$ de $\cO$} Nous pouvons maintenant déterminer la cohomologie $L^2$ de
  $\cO$. Nous commençons en degré $k<n$. Remarquons que
  puisque $$(Y_i^\eps\times (V_i\setminus \bB))/B_i \subset \cO$$ et que 
  $$\H^k((Y_i^\eps\times (V_i\setminus \bB))/B_i, (\partial Y_i^\eps\times (V_i\setminus \bB))/B_i)
  \rightarrow \H^k( E_i )\rightarrow\{0\}$$  surjective alors l'application 
  $$\H^k(\cO)\rightarrow \bigoplus_{i \ge 1} \H^k(E_i)$$ est surjective.
  Nous allons montrer que cette application est injective : soit donc $\alpha\in Z^k(\cO)$ tel que
  pour tout $i\ge 1$ : $$\iota^*_{E_i}\alpha\in B^k(E_i)=d\cC^{k-1}_{w,1}( E_i).$$
  Remarquons que puisque $k<n$, et que $E_0$ est un bout de cône, on sait que 
  $$\iota^*_{E_0}\alpha\in B^k(E_0)=d\cC^{k-1}_{w,1}( E_0).$$
  On trouve donc pour chaque $i\ge 0$, $\psi_i\in \cC^{k-1}_{w,1}( E_i)$ tel que
  $$\iota^*_{E_i}\alpha=d\psi_i$$
  Pour $i\ge 1$, nous posons alors $$\eta_i=\iota^*_{E_i\cap E_0}\psi_i -\iota^*_{E_i\cap E_0}\psi_0\in Z^{k-1}_w(E_i\cap E_0).$$
  Mais $E_i\cap E_0$ est un bout de cône et $k-1<n-1$, donc d'après (\ref{cohopcone}) nous avons 
  $$Z^{k-1}_w(E_i\cap E_0)=B^{k-1}_w(E_i\cap E_0)=d\cC^{k-2}_{\bar w,w}(E_i\cap E_0)$$
  Alors grâce à (\ref{condition}), nous trouvons $u_i\in \cC^{k-2}_{\bar w,w}(E_i)$
  tel que 
  $$d(\iota^*_{E_i\cap E_0}u_i)=\eta_i.$$
  La forme $\psi$ définie sur $\cO$ par $\iota^*_{E_0}\psi=\psi_0$ et
  $\iota^*_{E_i}\psi=\psi_i-du_i$ est bien un élément de $\cC_{w,1}^{k-1}(E_i)$ et
  $$d\psi=\alpha.$$
\noindent Puisque $w$ est à décroissance parabolique, nous avons montré que $\alpha\in B^k(\cO)$ et aussi que
l'image de $d$ est presque fermée sur $\cO$ en degré $k<n$ par rapport au poids $w$.
  
\noindent Évidemment une preuve identique permet de conclure en degré $k>n$ et pour la cohomologie $L^2$ de 
$\cO$ relative à $\partial \cO$.

Nous traitons maintenant le cas de degré $k=n$ pour la cohomologie $L^2$ relative :

Il nous faut donc démontrer que si $\alpha\in Z^n(\cO,\partial \cO)$ est tel que 
$$\forall i\ge 1, \ \iota^*_{E_i}\alpha\in B^k(E_i,\partial \cO)$$ alors
$$\alpha\in d\cC^{n-1}_{w,1}(\cO,\partial\cO).$$

\noindent Remarquons que si pour tout $i\ge 1$, nous avons $b_{n-2}(E_i\cap E_0)=0$ (par exemple si $n$ est pair)
alors la preuve précédente est valide.
Nous supposons donc que $n$ est impair. Alors on remarque que 
$b_{n-1}(\Sigma_0)=0$ et que pour tout $i\ge 1$, $b_{n-1}(\bS^{2m_i-1}\times\bS^{2n_i-1})=0$. Ainsi
grâce à (\ref{raffinement}), nous pouvons trouver
$\psi_0\in \cC^{n-1}_{\bar w,1}(E_0,\partial\cO)$ telle que
$$\iota^*_{E_0}\alpha=d\psi_0.$$
Nous savons également qu'il existe 
$\psi_i\in  \cC^{n-1}_{w,1}( E_i,\partial \cO)$ tel que
  $$\iota^*_{E_i}\alpha=d\psi_i$$
Pour $i\ge 1$, nous posons alors 
$$\eta_i=\iota^*_{E_i\cap E_0}\psi_i -\iota^*_{E_i\cap E_0}\psi_0\in Z^{n-1}_w(E_i\cap E_0,\partial \cO).$$
\begin{enumerate}
\item Si $b_{n-2}(E_i\cap E_0)=0$, alors comme précédemment nous trouvons
$u_i\in \cC^{k-2}_{\bar w,w}(E_i,\partial \cO)$
  tel que 
  $$d(\iota^*_{E_i\cap E_0}u_i)=\eta_i.$$
  
\item Si maintenant $b_{n-2}(E_i\cap E_0)\not=0$ mais que $n_i>1$, alors la preuve précédente montre que
nous pouvons en fait supposer que $\psi_i\in \cC^{n-1}_{\bar w,1}( E_i,\partial \cO)$ et donc que
$$\eta_i\in Z^{n-1}_{\bar w}(E_i\cap E_0,\partial \cO)=B^{n-1}_{\bar w}(E_i\cap E_0,\partial \cO)=
d\cC^{n-2}_{ w,\bar w}(E_i\cap E_0,\partial \cO)$$
On trouve donc $ u_i\in \cC^{n-2}_{ w,w}(E_i,\partial \cO)$
  tel que 
  $$d(\iota^*_{E_i\cap E_0}u_i)=\eta_i.$$
\item Il reste le cas où $n_i=1$ et $b_{n-2}(E_i\cap E_0)\not=0$ c'est à dire le cas où $n=3$.
Dans ce cas grâce au raffinement proposé en (\ref{raffinement},ii) nous trouvons
$u\in \cC^{1}_{\bar w, w}(E_i\cap E_0,\partial \cO)$ et 
$$v=f(|x|)\sigma\in \cC^{1}_{ w, w}(E_i\cap E_0,\partial \cO)$$ où
$\sigma$ est le tiré en arrière de la forme "volume" de $\bS^{1}$ sur
$$E_ i\cap E_0=\simeq C_1(]\eps/2,2\eps[\times \bS^3\times \bS^1)/B_i$$
tel que $\eta=du+dv$
Nous pouvons évidemment trouver $\bar u\in \cC^{1}_{\bar w, w}(E_i,\partial \cO)$
tel que
$\iota^*_{E_i\cap E_0}\bar u=u$ et On peut aisément étendre $v$ à $E_i$ grâce à
\begin{equation}\label{astuce}
\bar v=f(s(x))\sigma\end{equation}
où on définit d'abord $s$ sur $\widehat E_i$ par
$$s(y,v)=\left\{\begin{array}{ll}
\sqrt{|y|^2+\|v\|^2} &{\rm\  si\ } |y|\ge \frac{\eps}{2} \|v\|\\
\sqrt{1+(\eps/2)^2}\|v\|&{\rm\  si\ } |y|\le \frac{\eps}{2} \|v\|\\
\end{array}\right. ,$$
cet fonction est clairement $B_i$ invariante et $\bar v$ est bien définie et il est clair que
$\bar v \in\cC^{1}_{ w, w}(E_i,\partial \cO)$. On pose alors
$u_i=\bar u+\bar v\in \cC^{1}_{ w, w}(E_i,\partial \cO)$ et nous avons également
$$d(\iota^*_{E_i\cap E_0}u_i)=\eta_i.$$
\end{enumerate}
La forme $\psi$ définie sur $\cO$ par $\iota^*_{E_0}\psi=\psi_0$ et
  $\iota^*_{E_i}\psi=\psi_i-du_i$ est bien un élément de $\cC_{w,1}^{k-1}(E_i,\partial \cO)$ et
  $$d\psi=\alpha.$$
  
  Nous avons donc démontré que :
  \begin{prop}\label{coho2O}
  \begin{enumeroman}
  \item Si $k<n$,alors $d\cC^{k-1}_{w,1}(\cO)=B^k(\cO)$ et 
  $$\H^k(\cO)=\bigoplus_{i,n_i>1} \H^{k-2n_i+1}(Y_i)^{B_i}.$$
  \item Si $k\ge n$,alors $d\cC^{k-1}_{w,1}(\cO,\partial \cO)=B^k(\cO,\partial \cO)$ et 
  $$\H^k(\cO,\partial \cO)=\bigoplus_{i,n_i>1} \H^{k-1}(Y_i)^{B_i}.$$
  \end{enumeroman}
  \end{prop}

  \subsection{Cohomologie $L_w^2$ de $\cO$}
Les mêmes arguments que précédemment nous permettent de démontrer que
\begin{enumeroman}
\item Si $k<n$, alors $B^{k-1}_w(\cO)=d\cC^{k-2}_{w,w}(\cO)$ et 
$$\H^{k-1}_w(\cO)=\bigoplus_{p+q=k-1} \H^p(Y_i)^{B_i}\otimes 
\H_w^q(\C^{n_i}\setminus \bB)=\bigoplus_i \H^{k-2n_i}(Y_i)^{B_i}\oplus \bigoplus_{i,n_i=1}
\H^{k-1}(Y_i)^{B_i}$$

\item Si $k>n$, alors $B^{k-1}_w(\cO,\cO)=d\cC^{k-2}_{w,w}(\cO,\cO)$ et 
$$\H^{k-1}_w(\cO,\partial\cO)=\bigoplus_{p+q=k-1} \H^p(Y_i)^{B_i}\otimes 
\H_w^q(\C^{n_i}\setminus \bB,\bS^{2n_i-1})=\bigoplus_{i, n_i>2} \H^{k-2}(Y_i)^{B_i}.$$
\end{enumeroman}
 Comme précédemment le point délicat est le cas de degré $k-1=n-1$. Nous reprenons nos arguments en commençant
 par
 calculer la cohomologie $L_w^2$ de $\widehat E_i$ en degré $n-1$, celui ci se traite de façon similaire au
 cas de la cohomologie $L^2$ :
 \begin{enumerate}
 \item Si $b_{n-2}(\bS^{2m_i-1}\times \bS^{2n_i-1})=0$, alors nous montrons comme précédemment que
 $$B^{n-1}_w(\widehat E_i,Y_i\times\bS^{2n_i-1} )=d\cC^{n-2}_{w,w}(\widehat E_i,Y_i\times\bS^{2n_i-1}) ;$$
 et que 
 $$\H^{n-1}_w(\widehat E_i,Y_i\times\bS^{2n_i-1} )
 =\H^{n-1}_w(Y_i\times (\C^{n_i}\setminus \bB),Y_i\times\bS^{2n_i-1} )$$
 avec de plus si $n_i\not =2$ : 
 $$B^{n-1}_w(\widehat E_i,Y_i\times\bS^{2n_i-1} )=d\cC^{n-2}_{\bar w,w}(\widehat E_i,Y_i\times\bS^{2n_i-1})$$
 
 \item Supposons maintenant que $b_{n-2}(\bS^{2m_i-1}\times \bS^{2n_i-1})\not=0$, c'est à dire que
 $|n_i-m_i|=1$, alors nous avons deux cas:
 \begin{enumeroman}
 \item le premier est celui où
 $2n_i-1=n-2$ : soit $\alpha\in B^{n-1}_w(\widehat E_i,Y_i\times\bS^{2n_i-1})$ alors le raffinement
  présenté en (\ref{raffinement}) nous permet de trouver :
  $u\in \cC^{n-2}_{\bar w,w}(\Omega_i,Y_i\times\bS^{2n_i-1})$ et 
  $v\in \cC^{n-2}_{ w,w}(\Omega_i,Y_i\times\bS^{2n_i-1})$ tels que
  $v=f(r)\sigma$, $\sigma$ étant le tiré en arrière de la forme volume de $\bS^{2n_i-1}$ par
  la projection sur le dernier facteur de $\Omega_i=C_1(C_1(]\eps/2,2\eps[\times \bS^{2m_i-1}\times
  \bS^{2n_i-1})$ et
  $$\iota^*_{\Omega_i}\alpha=du+dv$$
  On peut alors comme (\ref{astuce}) étendre $u,v$ en 
  $\bar v,\bar u\in \cC^{n-2}_{ w,w}(\Omega_i,Y_i\times\bS^{2n_i-1})$ et en posant $\bar\psi=\bar u+\bar v$,
  les arguments précédents mènent de la même façon au résultat voulu :
  l'extension  par zéro de $\alpha-d\bar \psi$ à $Y_i\times(V_i\setminus \bB)$ est fermée :
  $$\alpha-d\bar \psi\in Z^{n-1}_{w'}(Y_i\times(V_i\setminus \bB),Y_i\times\bS^{2n_i-1})$$
  Mais l'image de la classe de cohomologie $L_{w'}^2$ de cette forme
   par $\iota^*_{\widehat E_i}$ est la classe de cohomologie de $\alpha\in B_{w'}^{n-1}(\widehat E_i,Y_i\times\bS^{2n_i-1})$, puisque
   comme précédemment cette
   application est injective, la proposition (\ref{kunneth3}) nous permet de trouver 
   $f\in \cC^{n-2}_{\bar w_1,w'}(Y_i\times(V_i\setminus \bB),Y_i\times\bS^{2n_i-1})$ et
    $g\in \cC^{n-2}_{w',w'}(Y_i\times(V_i\setminus \bB),Y_i\times\bS^{2n_i-1})$
  tels que
  $$\alpha-d\bar\psi=df+dg$$
  Compte tenu du fait que sur $\widehat E_i$ : $\bar w_1\ge C w'$ et que $w'$ et comparable à $w$, alors
  $$\psi= \bar\psi+\iota^*_{\widehat E_i}(f+g)\in \cC^{n-2}_{w,w}(\widehat E_i,Y_i\times\bS^{2n_i-1}).$$
  
  Ainsi 
  $$\cC^{n-2}_{w,w}(\widehat E_i,Y_i\times\bS^{2n_i-1})=B^{n-1}_w(\widehat E_i,Y_i\times\bS^{2n_i-1})$$ et de
  la même façon l'application :
  $$\iota^*_{\widehat E_i}\,:\, \H_{w'}^{n-1}(Y_i^\eps\times (V_i\setminus \bB),
 Y_i\times\bS^{2n_i-1})
 \rightarrow \H^{n-1}_{w'}( \widehat E_i,Y_i\times\bS^{2n_i-1} )$$ est un isomorphisme
  et l'application d'extension par zéro 
 $$\H_{w'}^{n-1}(Y_i^\eps\times (V_i\setminus \bB),
  \partial (Y_i^\eps\times (V_i\setminus \bB))\cup Y_i\times\bS^{2n_i-1})
 \rightarrow \H^{n-1}_{w'}( \widehat E_i,Y_i\times\bS^{2n_i-1} )$$ est également surjective.

 \item Supposons maintenant que $2m_i-1=n-2$, remarquons qu'alors $m_i=(n-1)/2$ et $n_i=m_i+1$ puisque 
 $m_i\ge 2$, ainsi $n_i\ge 3$.
 Soit donc  $\alpha\in Z^{n-1}_{w'}(\widehat E_i,Y_i\times\bS^{2n_i-1})$,
 alors nous trouvons comme précédemment
   $\psi\in \cC^{n-2}_{ w',w'}(\Omega_i, Y_i\times \bS^{2n_i-1})$ tel que 
  $$\iota^*_{\Omega_i}\alpha=d\psi$$
  Maintenant grâce à (\ref{condition}), nous savons qu'il existe 
  $\bar \psi\in \cC^{n-1}_{ w',\rho w'}(\widehat E_i, Y_i\times \bS^{2n_i-1})$
  tel que $\iota^*_{\Omega_i}\bar \psi=\psi$. En particulier, 
  l'extension à $Y_i\times(V_i\setminus \bB)$ par zéro de $\alpha-d\bar\psi$ est fermée :
  $$\alpha-d\bar \psi\in Z_{\rho w'}^{n-1}(Y_i\times(V_i\setminus \bB, Y_i\times \bS^{2n_i-1})).$$
  Les mêmes causes produisant les mêmes effets, il est également vrai que l'application :
  $$\iota^*_{\widehat E_i}\,:\,\H_{\rho w'}^{n-1}(Y_i\times(V_i\setminus \bB), Y_i\times \bS^{2n_i-1} )
  \rightarrow\H_{\rho w'}^{n-1}(\widehat E_i , Y_i\times \bS^{2n_i-1})$$
  est injective.
  Et également d'après (\ref{kunneth3}), l'application naturelle :
  $$\H_{ w'}^{n-1}(Y_i\times(V_i\setminus \bB), Y_i\times \bS^{2n_i-1} )\rightarrow 
  \H_{\rho w'}^{n-1}(Y_i\times(V_i\setminus \bB), Y_i\times \bS^{2n_i-1} )$$ est un isomorphisme
  
  Mais l'image de la classe de cohomologie $L_{\rho w'}^2$ de $\alpha-d\bar \psi$
   par $\iota^*_{\widehat E_i}$ est la classe de cohomologie de $\alpha\in Z_{\rho w'}^{n-1}(\widehat E_i, Y_i\times \bS^{2n_i-1})$, 
   la proposition (\ref{kunneth3}) nous permet de trouver 
   $u\in \cC^{n-2}_{\bar w_1,\rho w'}(Y_i\times(V_i\setminus \bB), Y_i\times \bS^{2n_i-1})$,
    $v\in \cC^{n-2}_{\bar w',\rho w'}(Y_i\times(V_i\setminus \bB), Y_i\times \bS^{2n_i-1})$ et 
    $h\in \cH^{n-1}_{w'}(Y_i\times(V_i\setminus \bB), Y_i\times \bS^{2n_i-1} )$
  tels que
  $$\alpha-d\bar\psi=h+du+dv$$
  Si on pose $\phi=\psi+\iota^*_{\widehat E_i}(u+v)$ alors
  $\alpha=\iota^*_{\widehat E_i} h+d\phi$ et $\phi\in L^2_{(w')^2}(\Lambda^{n-2}T^*(E_i))$.
  Nous avons donc obtenu le même résultat qu'en (i).
 \end{enumeroman}
 \end{enumerate}
  En revenant à $E_i=\widehat E_i/B_i$, nous avons démontré que  
  $d\cC^{n-2}_{w,w}( E_i,\partial\cO)=B^{n-1}_w( E_i,\partial\cO)$
  même mieux si $|n_i-m_i|\not= 1$ et $n_i\not=2$ alors
  $$ d\cC^{n-2}_{\bar w,w}( E_i,\partial\cO)=B^{n-1}_w( E_i,\partial\cO).$$ De plus 
  $$\H^{n-1}_w( E_i,\partial\cO )\simeq \bigoplus_{p+q=n-1} \H^p(Y_i)^{B_i}\otimes 
  \H_{w'}^q(V_i\setminus \bB,\bS^{2n_i-1})^{B_i}=
  \left\{\begin{array}{ll} \H^{k-1}(Y_i)^{B_i}&{\rm\ si\ }n_i>2\\
  \{0\}&{\rm\ si\ }n_i\le 2\\ \end{array}\right.$$
  \noindent Et l'application d'extension par zéro 
  $$\H^{n-1}_w( (Y_i^\eps\times (V_i\setminus \bB))/B_i
  ,\partial\cO\cup (\partial Y_i^\eps\times (V_i\setminus \bB))/B_i)
  \rightarrow \H^{n-1}_w( E_i )$$ est surjective.
  
  On peut maintenant calculer la cohomologie $L^2_w$ de $\cO$ relative à $\partial \cO$ en degré $n-1$:
  pour les mêmes raisons l'application 
  $$\H^{n-1}_w(\cO,\partial\cO)\rightarrow \bigoplus_{i \ge 1} \H^{n-1}_w(E_i,\partial\cO)$$ est surjective.
  
Il nous faut donc démontrer que si $\alpha\in Z^{n-1}_w(\cO,\partial \cO)$ est tel que 
$$\forall i\ge 1, \ \iota^*_{E_i}\alpha\in B^{n-1}_w(E_i,\partial \cO)$$ alors
$$\alpha\in d\cC^{n-2}_{w,w}(\cO,\partial\cO).$$
\begin{description}
\item[Premier cas] supposons que $n$ est impair et $n>3$,
 nous pouvons trouver
$\psi_0\in \cC^{n-2}_{ w,w}(E_0,\partial\cO)$ telle que
$$\iota^*_{E_0}\alpha=d\psi_0.$$
Nous savons également qu'il existe 
$\psi_i\in  \cC^{n-2}_{w,w}( E_i,\partial \cO)$ tel que
  $$\iota^*_{E_i}\alpha=d\psi_i.$$
Pour $i\ge 1$, nous posons alors 
$$\eta_i=\iota^*_{E_i\cap E_0}\psi_i-\iota^*_{E_i\cap E_0}\psi_0
 \in Z^{n-2}_{ w^2}(E_i\cap E_0,\partial \cO).$$ 
 mais puisqu'alors nous avons $b_{n-3}(\bS^{2m_i-1}\times \bS^{2n_i-1})=0$, la proposition
 (\ref{raffinement}), nous permet de trouver $u_i\in \cC^{n-3}_{\bar w, w^2}(E_i\cap E_0,\partial \cO)$
 tel que $\eta_i=du_i$, on peut grâce à \ref{condition} étendre $u_i$ en $\bar u_i\in 
 \cC^{n-3}_{\bar w, w^2}(E_i,\partial \cO)$ et alors la forme $\psi$ définie sur $\cO$ par $\iota^*_{E_0}\psi=\psi_0$ et
  $\iota^*_{E_i}\psi=\psi_i-d\bar u_i$ est bien un élément de $\cC^{n-2}_{w,w}(\cO,\partial \cO)$ et
  $$d\psi=\alpha.$$
  \item[Deuxième cas] si $n=3$, alors nous n'avons clairement pas $b_{n-3}(\bS^{3}\times \bS^{1})=0$.
  Cependant grâce à (\ref{raffinement}), nous pouvons 
  en fait trouver pour chaque $i>0$, $f_i\in \cC^{0}_{\bar w, w^2}(E_i\cap E_0,\partial \cO)$ et 
  $g_i\in \cC^{0}_{ w, w^2}(E_i\cap E_0,\partial \cO)$ une fonction radiale tel que
  $$\eta_i=df_i+dg_i,$$ on étend alors aisément ces fonctions à $E_i$ pour obtenir de même une fonction
  $\bar u_i\in \cC^{0}_{ w, w^2}(E_i,\partial \cO)$ tel que 
  $$\iota^*_{E_i\cap E_0}d\bar u_i =\eta_i,$$
  Ceci permet alors de conclure comme précédemment.
  
\item[Troisième cas]  supposons que $n$ est pair alors nous avons forcément pour tout $i$ : $|n_i-m_i|\not =1$ 
et également $b_{n-2}(\Sigma_0)=0$, alors nous pouvons en fait trouver  
$\psi_0\in \cC^{n-2}_{ \bar w,w}(E_0,\partial\cO)$ tel que
$$\iota^*_{E_0}\alpha=d\psi_0.$$
Nous savons également qu'il existe 
$\psi_i\in  \cC^{n-2}_{w,w}( E_i,\partial \cO)$ tel que
  $$\iota^*_{E_i}\alpha=d\psi_i.$$
Pour $i\ge 1$, nous posons alors 
$$\eta_i=\iota^*_{E_i\cap E_0}\psi_i-\iota^*_{E_i\cap E_0}\psi_0
 \in Z^{n-2}_{w^2}(E_i\cap E_0,\partial \cO).$$ 
\begin{enumerate}
\item Si $b_{n-3}(E_i\cap E_0)=0$, alors comme précédemment nous trouvons
$\bar u_i\in \cC^{n-3}_{\bar w,w^2}(E_i,\partial \cO)$
  tel que 
  $$d(\iota^*_{E_i\cap E_0}\bar u_i)=\eta_i.$$
  
\item Si maintenant $b_{n-3}(E_i\cap E_0)\not=0$ et
$n-3=2n_i-1$ alors  grâce à (\ref{raffinement}) et en 
reprenant l'argument utilisé dans (\ref{astuce})
on trouve $ \bar u_i\in \cC^{n-3}_{ w,w^2}(\cO,\partial \cO)$
  tel que 
  $$d(\iota^*_{E_i\cap E_0}u_i)=\eta_i.$$
\item Il reste le cas où $b_{n-3}(E_i\cap E_0)\not=0$ mais que $2m_i-1=n-3$ et $2n_i-1=n+1$, remarquons que
puisque $m_i\ge 2$, nous avons forcément $n\ge 6$ et $n_i\ge 4>2$ en particulier, nous pouvions en fait
 choisir $\psi_i\in  \cC^{n-2}_{\bar w,w}( E_i,\partial \cO)$ et alors
 $\eta_i\in Z^{n-2}_{\bar ww}( E_i\cap E_0,\partial \cO)$ et on trouve alors
 $u_i\in \cC^{n-2}_{w,\bar ww}( E_i\cap E_0,\partial \cO)$ tel que $du_i=\eta_i$, cette forme 
 peut être étendue en $\bar u_i\in \cC^{n-2}_{w,ww}( E_i,\partial \cO)$.
\end{enumerate}
La forme $\psi$ définie sur $\cO$ par $\iota^*_{E_0}\psi=\psi_0$ et
  $\iota^*_{E_i}\psi=\psi_i-d\bar u_i$ est bien un élément de $\cC^{n-2}_{w,w}(\cO,\partial \cO)$ et
  $$d\psi=\alpha.$$
  
 \end{description}
  Nous avons donc obtenu :
\begin{prop}\label{coho2pO}
 \begin{enumeroman}
  \item Si $k<n$,alors $d\cC^{k-2}_{w,w}(\cO)=B_w^{k-1}(\cO)$ et 
  $$\H^{k-1}(\cO)=\bigoplus_{i} \H^{k-2n_i}(Y_i)^{B_i}\oplus \bigoplus_{i, n_i=1}\H^{k-1}(Y_i)^{B_i} .$$
  \item Si $k\ge n$,alors $d\cC^{k-2}_{w,w}(\cO,\partial \cO)=B^{k-1}(\cO,\partial \cO)$ et 
  $$\H^{k-1}(\cO,\partial \cO)=\bigoplus_{i,n_i>2} \H^{k-2}(Y_i)^{B_i}.$$
  \end{enumeroman}
\end{prop}

En conséquence pour $k<n$, les images de 
$d\,:\,\cC^{k-2}_{w,w}(\cO)\rightarrow L_w^2(\Lkm \cO)$  et de 
$d\,:\,\cC^{k-1}_{w,1}(\cO)\rightarrow L^2(\Lk \cO)$ sont fermées et l'espace
$\H^{k-1}_w(\cO)$ est de dimension finie,
grâce au théorème (\ref{conclusion}), il en est de même sur toute variété isométrique à $\cO$ au dehors
d'un compact ; par exemple pour les degrés $k<n$, les images de
 $d\,:\,\cC^{k-1}_{w,1}(\cO,\partial\cO)\rightarrow L^2(\Lk \cO)$
et de $d\,:\,\cC^{k-1}_{w,1}(X)\rightarrow L^2(\Lk X)$ sont fermées.
Grâce à la version relative du même théorème, on peut affirmer que pour $k\ge n$, les images de
$d\,:\,\cC^{k-1}_{w,1}(\cO)\rightarrow L^2(\Lk \cO)$
et de $d\,:\,\cC^{k-1}_{w,1}(X)\rightarrow L^2(\Lk X)$ sont fermées.

 Nous devons également connaître la $L^2_{1/w}$ cohomologie de $\cO$ relative à $\partial \cO$ en degré 
 $(k+1)$ pour $k\ge n$.  Lorsque $k>n$, nous obtenons
 de la même façon que 
 $d\cC^{k}_{w,1/w}(\cO,\partial \cO)=B^{k+1}_{1/w}(\cO,\partial \cO)$ et 
  $$\H^{k+1}_{1/w}(\cO,\partial \cO)=
  \bigoplus_{i} \H^{k}(Y_i)^{B_i}\oplus \bigoplus_{i,n_i=1} \H^{k-1}(Y_i)^{B_i}.$$
  
  Concernant le degré $k=n+1$, nous pouvons reprendre les arguments précédents mais nous ne pouvons conclure
  que lorsque $n\not =3$. En effet, ces arguments ne permettent pas dans ce cas de déterminer
  la $L^2_{1/w}$ cohomologie de $\widehat E_i$ relative à $Y_i\times (V_i\setminus \bB)$. Nous ne pouvons
 en effet conclure que lorsque $b_n(\bS^{2m_i-1}\times \bS^{2n_i-1})=0$, ou lorsque 
  $b_n(\bS^{2m_i-1}\times \bS^{2n_i-1})\not=0$ mais que $n_i>1$ ou encore lorsque
  $n_i=1$ mais $H^n(\bS^{2m_i-1}\times \bS^{2n_i-1})\simeq H^n( \bS^{2n_i-1})$ ;
   mais nous ne pouvons conclure dans le dernier cas  où $b_n(\bS^{2m_i-1}\times \bS^{2n_i-1})\not=0$, $n_i=1$ et $n=2m_i-1$, ce qui se
produit exactement lorsque $n=3$. Nous obtenons ainsi :
  
 \begin{prop}\label{coho2mO}
 Si $k\ge \max{n,4} $,alors $d\cC^{k}_{w,1/w}(\cO,\partial \cO)=B^{k+1}_{1/w}(\cO,\partial \cO)$ et 
  $$\H^{k+1}_{1/w}(\cO,\partial \cO)=\bigoplus_{i} \H^{k}(Y_i)^{B_i}
  \oplus \bigoplus_{i,n_i=1} \H^{k-1}(Y_i)^{B_i}.$$
\end{prop}

Nous pouvons maintenant déterminer la cohomologie $L^2$ de $\cO$ en degré $k\ge \max\{n,4\}$ . En effet, on dispose 
maintenant de la suite exacte :
$$ H^{k-1}(\partial \cO)
 \longrightarrow \H^{k}_{}(\cO,\partial \cO)\longrightarrow
\H^k(\cO)\longrightarrow H^k(\partial \cO)\longrightarrow\H^{k+1}_{1/w}(\cO,\partial \cO).$$
Nos calculs de la cohomologie de $\partial K$ (\ref{cohopK}) et de la cohomologie $L^2$ relative de 
$\cO$ (\ref{coho2O}) montre que la première flèche de cette suite exacte est surjective, nous obtenons donc  
$$ \H^k( \cO)\simeq \ker\left(H^k(\partial \cO)\rightarrow\H^{k+1}_{1/w}(\cO,\partial \cO) \right).$$
Puisque pour $k\in [1, 2n-2]$
$$H^k(\partial \cO)
=\bigoplus_{i, k>2n_i-1} H^{k-2n_i+1}(Y_i)^{B_i}\oplus \bigoplus_{i} H^{k}(Y_i)^{B_i}$$
et 
$$\H^{k+1}_{1/w}(\cO,\partial \cO)
=\bigoplus_{i,} \H^{k}(Y_i)^{B_i}\oplus \bigoplus_{i,n_i=1} \H^{k-1}(Y_i)^{B_i},$$
nous obtenons également :
$$\H^k(\cO)\simeq \bigoplus_{i,n_i>1} \H^{k-2n_i+1}(Y_i)^{B_i}.$$
Pour $k=2n-1$, nous avons alors $\H^{2n}_{1/w}(\cO,\partial \cO)=\{0\}$ et donc 
$$\H^{2n-1}( \cO)\simeq H^{2n-1}(\cO)\simeq \R.$$

Remarquons que nous pouvons également déterminer la cohomologie $L^2$ de $\cO$ en degré $n=3$, car nous avons
$$\H^3(\cO,\partial \cO)=\{0\}$$ et grâce à la dualité induite par l'opérateur de Hodge nous obtenons
également $$\H^3(\cO)=\{0\}= \bigoplus_{i,n_i>1} \H^{3-2n_i+1}(Y_i)^{B_i}.$$

On peut de même déterminer la cohomologie $L^2_w$ de $\cO$ pour les degrés $k-1\ge n$ grâce à la suite exacte :
$$ H^{k-2}(\partial \cO)
 \longrightarrow \H^{k-1}_{w}(\cO,\partial \cO)\longrightarrow
\H^{k-1}_{w}(\cO)\longrightarrow H^{k-1}(\partial \cO)\longrightarrow\H^{k}_{}(\cO,\partial \cO).$$
Comme précédemment la première flèche de cette suite est surjective et 
grâce à (\ref{coho2pO}), on obtient également :
$$\H^{k-1}_w(\cO)\simeq \bigoplus_{i\ge 1} \H^{k-2n_i}(Y_i)^{B_i}\oplus\bigoplus_{i,n_i=1}\H^{k-1}(Y_i)^{B_i}.$$

Nous avons donc obtenu :
\begin{prop}\label{cohoinfty} Pour $k\in [1,2n-2]$,
la cohomologie $L^2$ de $\cO$ est donnée par
\begin{equation}\label{cohoi}\H^k(\cO)\simeq \bigoplus_{i,n_i>1} \H^{k-2n_i+1}(Y_i)^{B_i}
\simeq\bigoplus_{i,k>2n_i-1,n_i>1} H^{k-2n_i+1}(Y_i)^{B_i}.\end{equation}
Et en degré $(2n-1)$ : $$\H^{2n-1}( \cO)\simeq H^{2n-1}(\cO)\simeq \R.$$
En degré $k-1\in [0,2n-2]$, la cohomologie $L^2_w$ de $\cO$ est donnée par :
\begin{equation}\label{cohowi}\begin{split}
\H^{k-1}_w(\cO)
&\simeq \bigoplus_{i\ge 1} \H^{k-2n_i}(Y_i)^{B_i}\oplus\bigoplus_{i,n_i=1}\H^{k-1}(Y_i)^{B_i}\\
&\simeq \bigoplus_{i\ge 1,k>2n_i} H^{k-2n_i}(Y_i)^{B_i}\oplus\bigoplus_{i,n_i=1, k>1}H^{k-1}(Y_i)^{B_i}.\\
\end{split}
\end{equation}
\end{prop}
\subsection{Cas des résolutions crépantes}\ 
\subsubsection{Cas général} Supposons que $X\stackrel{\pi}{\longrightarrow} \C^n/G$ soit une résolution
crépante. Nous pouvons dans ce cas 
donner une description explicite de la cohomologie $L^2$ de $X$. En effet, on connaît la cohomologie de $X$
en fonction des classes de conjugaison de $G$ :
si $g\in G$ on note $\age(g)=\theta_1+\theta_2+...+\theta_n$ où les valeurs propres de $g$ sont 
$(e^{i2\pi \theta_1},e^{i2\pi \theta_2},..., e^{i2\pi \theta_n})$ avec $\theta_i\in [0,1[$, puisque
$\det g=1$ on a forcément $\age(g)\in\N\cap [0,n[$, il est aussi clair que $\age(g)$ ne dépend que de la
classe de la conjugaison de $g$ $[g]\in \cC(G)$.
Les résultats de Y. Ito et M. Reid ($n=3$) et V. Batyrev et indépendamment J. Denef et F. Loeser sont les
suivants (\cite{IR, Baty,DL}) :
la cohomologie de $X$ est nulle en degré impair et en degré pair 
$$\dim  H^{2k}(X)=\card\{ [g]\in \cC(G),\ \age(g)=k\}.$$
Notre résultat est le suivant :
  
\begin{thm} Soit $G\subset \Su (n)$ un sous groupe fini, on suppose que $\bS^{2n-1}/G
$ le quotient de la sphère par $G$ soit à singularités isolées. Si $X\stackrel{\pi}{\longrightarrow} \C^n/G$ une
 résolution crépante de $\C^n/G$ équipée d'une métrique QALE asymptote à $\C^n/G$ alors
 $$\H^k(X)\simeq \ima\left( H_c^k(X)\rightarrow H^k(X)\right).$$
 
\end{thm}

\proof Nous avons obtenu la suite exacte suivante :
\begin{equation}
\label{se1}
 H^{k-1}(K)\oplus \H^{k-1}_{w}(\cO)
 \stackrel{\delta}{\longrightarrow} \H^{k-1}_{}(\partial K)\stackrel{\mbox{b}}{\longrightarrow}
 \H^k(X) \stackrel{\mbox{r}^*}{\longrightarrow} 
 H^k(K)\oplus \H^k(\cO)\stackrel{\delta}{\longrightarrow} H^k(\partial K).
\end{equation}
{\bf Si $k$ est pair :} alors puisque les $Y_i$ sont des résolutions crépantes de $W_i/A_i$, elles n'ont de
cohomologie qu'en degré pair ainsi d'après (\ref{cohoinjty}) et (\ref{injALE}) on a forcément 
$\H^k(\cO)=\{0\}$ ; de plus d'après (\ref{cohoinjty}), on sait que
$$ \H_w^{k-1}(\cO)=\bigoplus_{i\ge 1} \H^{k-2n_i}(Y_i)^{B_i}$$
et aussi (\ref{cohopK}): $$ H^{k-1}(\partial K)\simeq \bigoplus_{i\ge 1, k>2n_i} H^{k-2n_i}(Y_i)^{B_i}.$$
Or lorsque $l>0$,  $\H^l(Y_i)=H^l(Y_i)$, on obtient donc l'isomorphisme : 
$$\H_w^{k-1}(\cO)\simeq H^{k-1}(\partial K).$$ La suite exacte (\ref{se1}) se réduit à :
$$\{0\}\rightarrow  \H^k(X) \rightarrow
 H^k(K)\rightarrow H^k(\partial K),$$
 Donc $\H^k(X)$ est isomorphe au noyau de $ H^k(K)\rightarrow H^k(\partial K)$ qui est exactement l'image de 
  $H^k(K,\partial K)\rightarrow H^k(K)$. Puisque $X$ se rétracte sur $K$, on a bien démontré le résultat.
  
  {\bf Si $k$ est impair :}
  dans ce cas on sait $H^k(K)=\{0\}$ et $H^k(K,\partial K)=\{0\}$ donc
  $$H^{k-1}(K)\rightarrow H^{k-1}(\partial K)\rightarrow \{0\}\ .$$
 De la suite exacte (\ref{se1}), il ne reste que
  $$ \{0\}\rightarrow  \H^k(X) \rightarrow
 \H^k(\cO)\rightarrow H^k(\partial K).$$
Or $$\H^k(\cO)=\bigoplus_{i\ge 1, n_i>1} \H^{k-2n_i-1}(Y_i)^{B_i}$$ et
  $$H^k(\partial K)=\bigoplus_{i\ge 1,k>2n_i-1 } H^{k-2n_i-1}(Y_i)^{B_i},$$
 sauf pour $k=2n-1$ où $\H^{2n-1}(\cO)\simeq \H^{2n-1}(\partial K)$.
  La dernière flèche de cette suite est donc injective; la cohomologie $L^2$ de $X$ est donc nulle.
  \endproof
  
  On peut étudier maintenant étudier deux cas particuliers :
 \subsubsection{le cas particulier où $G\subset \Su (3)$}
On suppose ici que $G$ est un sous groupe fini de $\Su(3)$, alors $\C^3/G$ admet toujours des résolutions
crépantes et soit $X\stackrel{\pi}{\longrightarrow} \C^3/G$ l'une d'entre elles. 
Lorsque $g\in G$, on a forcément $$\dim\ker(g-\Id)\in\{0,1,3\}$$ et donc le lieu singulier de
$\bS^5/G$ est une réunion disjointes de cercles. En particulier, si $K=\pi^{-1}(\bB)$ et
$\cO=X\setminus K$ alors  pour $k\le 3$ : 
$\H^k(\cO)=\{0\}$  et  pour $k\ge 3$ $\H^k(\cO,\partial \cO)=\{0\}$.
Puisque $d\cC_{w,1}^{k-1}(\cO)=Z^k(\cO)$ pour $k\le 3$, on en déduit facilement que
l'application naturelle $H^k_c(X)\rightarrow \H^k(X)$ est surjective pour $k\le 3$.
Puisque $\cO$ est connexe et que $H^1(\cO)=\{0\}$, le lemme (2.1) de \cite{Carpar} montre que cette
application est aussi injective pour $k\le 2$. De plus, puisque pour $k\ge 3$, $\H^k(\cO,\partial \cO)=\{0\}$, 
le lemme (2.2) de \cite{Carpar}) implique que l'application $\H^k(X)\rightarrow H^k(X)$ est injective pour
$k\ge 3$ ; 
on en déduit immédiatement :
\begin{thm}Soit $G\subset SU(3)$ un sous-groupe fini et $\pi\,:\, X\longrightarrow \C^3/G$
 une résolution
crépante, on équipe $X$ d'une métrique QALE alors 
$$\H^k(X)=\left\{\begin{array}{lll}
\{0\}& {\rm si\ }& k\not=2,4\\
H^2_c(X)\simeq \ima (H^2_c(X)\rightarrow H^2(X))& {\rm si\ }& k=2\\
H^4(X)\simeq \ima (H^4_c(X)\rightarrow H^4(X))& {\rm si\ }& k=4\\
\end{array}\right. .$$
En général, pour toute variété $(M,g)$ isométrique au dehors d'un compact avec $X$ alors
$$\H^k(M)=\left\{\begin{array}{lll}

H^k_c(M)& {\rm si\ }& k\le 2\\
 \ima (H^3_c(M)\rightarrow H^3(M))& {\rm si\ }& k=3 \\
H^k(M)& {\rm si\ }& k\ge 4\\
\end{array}\right. .$$
 
\end{thm}

On revient au cas de $X$ une résolution crépante de $\C^3/G$ :
$$\chi_{L^2}(X)=\sum_{l=0}^6 (-1)^l \dim\H^l(X)=2\dim \H^4(X)=2\card\{ [g]\in \cC(G),\ \age(g)=2\}.$$
Mais on a toujours 
$$\age(g)+\age(g^{-1})=3-\dim\ker(g-\Id),$$
ainsi $$\ker(g-\Id)\not=\{0\}\Leftrightarrow \age(g)=\age(g^{-1})=1$$
et $$\ker(g-\Id)=\{0\}\Leftrightarrow \{\age(g),\age(g^{-1})\}=\{1,2\}.$$
On a donc montré le résultat suivant
\begin{cor}
$$\chi_{L^2}(X)=\dim \H^2(X)+\dim \H^4(X)=\card\{ [g]\in \cC(G),\ \ker(g-\Id)=\{0\}\}.$$
\end{cor}
Considérons l'exemple 9.3.5 du livre de D. Joyce : $G\simeq \Z_4$ est le groupe engrendé par
$(z_1,z_2,z_3)\mapsto (-z_1,iz_2,iz_3)$. Si $X$ est la résolution crépante de $\C^3/G$ équipée d'une métrique 
QALE alors :
$$\dim \H^2(X)=\dim \H^4(X)=1$$

Si on considère l'exemple 9.3.6 de ce même livre : $G=\Z_2^2$ est le sous groupe de $\Su(3)$ formé par les
matrices diagonales à entrée $\pm 1$· Dans ce cas tous les éléments de $G$ fixe au moins une droite de $\C^3$
et donc les résolutions crépantes de $\C^3/G$ ne portent pas de formes harmoniques $L^2$.

 \subsubsection{le cas particulier où $G\subset \Sp (2)$}
 Soit donc $G$ un sous groupe fini de $\Sp(2)$ et $X\stackrel{\pi}{\longrightarrow} \C^n/G$
  une résolution crépante\footnote{Cela n'existe pas forcément.}.
Mais la dimension du noyau d'un élément de $G$ est  $0,2$ ou $4$.
Puisque que les variétés ALE hyperkählerienne de dimension complexe $2$
ne portent des formes harmoniques $L^2$ qu'en degré $2$ et qu'en ce degré la cohomologie $L^2$ est la
cohomologie ordinaire, on obtient que
pour $$k\not=5,7 \ \H^k(\cO)=\H^k_w(\cO)=\{0\}
\ \et\ \H^5(\cO)=\H^5_w(\cO)=\bigoplus_{i\ge 1} H^2(Y_i)^{B_i}.$$

Pour la cohomologie ordinaire de $\cO$ on sait que :
$H^k(\cO)=\{0\}$ pour $k=3,4,6,8$ et de plus
$H^5(\cO)=\bigoplus_{i\ge 1} H^2(Y_i)^{B_i}\simeq \H^5(\cO)$ et aussi
$H^7(\cO)\simeq \H^7(\cO)$.
Grâce à nos suites exactes de Mayer-Vietoris, on obtient alors facilement :
$\H^k(X)=H^k(X)$ dès que $k\ge 4$. En fait cette analyse est valide pour toute variété isométrique au
dehors d'un compact avec $X$, on obtient donc :

\begin{thm}Soit $G\subset \Sp(2)$ un sous-groupe fini et $\pi\,:\, X\longrightarrow \C^4/G$
 une résolution
crépante, on équipe $X$ d'une métrique QALE alors 
$$\H^k(X)=\left\{\begin{array}{lll}
\{0\}& {\rm si\ }& k\not=2,4,6\\
H^2_c(X)\simeq \ima (H^2_c(X)\rightarrow H^2(X))& {\rm si\ }& k=2\\
H^4(X)\simeq H^4_c(X)\simeq H^4(X) & {\rm si\ }& k=4\\
H^6(X)\simeq \ima (H^6_c(X)\rightarrow H^6(X)) & {\rm si\ }& k=6\\\end{array}\right. .$$
En général, pour toute variété $(M,g)$ isométrique au dehors d'un compact avec $X$ alors
$$\H^k(M)=\left\{\begin{array}{lll}

H^k_c(M)& {\rm si\ }& k\le 4\\
 
H^k(M)& {\rm si\ }& k\ge 4\\
\end{array}\right. .$$
 
\end{thm}

Remarquons que si une résolution crépante de $\C^4/G$ n'existe pas toujours, il existe toujours une
résolution crépante de $(\C^4\setminus \bB)/G$.

Revenons maintenant au cas d'une résolution crépante de $\C^4/G$.
Lorsque $g\in\Sp(2)$ on a de même :
$$\age(g)+\age(g^{-1})=4-\dim\ker(g-\Id).$$
En conséquence on a aussi $$\dim\ker(g-\Id)=2\Leftrightarrow \age(g)=\age(g^{-1})=1.$$
Et puisque 
$\dim\H^2(X)=\dim\H^6(X)=\card\{[g]\in\cC(D), \age(g)=3\}$ et 
$\dim\H^4(X)=\dim H^4(X)=\card\{[g]\in\cC(D), \age(g)=2\}$, on obtient encore :
\begin{cor}
\begin{equation*}
\begin{split}
\chi_{L^2}(X)&=\dim \H^2(X)+\dim \H^4(X)+\dim \H^6(X)\\
&=\card\{ [g]\in \cC(G),\ \ker(g-\Id)=\{0\}\}.\\
\end{split}\end{equation*}
\end{cor}

Étudions maintenant les exemples du livre de D. Joyce :
le premier (exemple 9.3.9) est le schéma de Hilbert de $3$ points sur $\C^2$.
Dans ce cas $G$ est le groupe $S_3$ des permutations de $\{1,2,3\}$ agissant sur $\{(x,y,z)\in (\C^2)^3,
x+y+z=0\}\simeq \C^4$. Ce schéma de Hilbert de $3$ points sur $\C^2$ peut être muni d'une métrique
hyperkählerienne QALE et alors la variété obtenue ne porte de formes harmoniques $L^2$ qu'en degré $4$ où
cet espace est de dimension $1$.

Le deuxième exemple est en fait le schéma de Hilbert de $2$ points sur une surface hyperkählerienne ALE, le
groupe $G$ en question est le produit semi-direct de $\Z_2$ avec le produit $H\times H$ où
$H$ est un sous-groupe fini de $\Su(2)=\Sp(1)$. Dans ce cas si $Y$ est la résolution crépante de 
 $\C^2/H$, ce Schéma de Hilbert ${\rm Hilb}^2(Y)$ ne porte de formes harmoniques $L^2$ qu'en degré $4$ et 
 $$\dim \H^4({\rm Hilb}^2(Y))=\frac{b_2(Y)(b_2(Y)+1)}{2}.$$
 \subsection{Le cas général}Nous allons maintenant donner une interprétation topologique des espaces de
 formes harmoniques $L^2$ des variétés QALE. La réponse est ici moins lisible que pour les résolutions
 crépantes car on ne connaît pas la cohomologie des résolutions localement en produit de $\C^n/G$. 
 Soit donc $X$ une variété QALE asymptote à $\C^n/G$ où $G$ est un sous groupe de $\Su(n)$
 tel que $\bS^{2n-1}/G$ soit à singularités isolés. Nos résultats précédents impliquent que la suite
 exacte courte suivante est valide :
 $$ H^{k-1}(K)\oplus \H^{k-1}_{w}(\cO)
 \stackrel{\delta}{\longrightarrow} \H^{k-1}_{}(\partial K)\stackrel{\mbox{b}}{\longrightarrow}
 \H^k(X) \stackrel{\mbox{r}^*}{\longrightarrow} 
 H^k(K)\oplus \H^k(\cO)\stackrel{\delta}{\longrightarrow} H^k(\partial K).
$$

On peut déjà déterminer la cohomologie $L^2$ de $X$ en degré $1$ :
nous avons en effet $\H^0_w(\cO)=\{0\}$ et 
$$\H^1(\cO)=\bigoplus_{1>2n_i-1, n_i>1} \H^{2-2n_i}(Y_i)=\{0\}\ ;$$
On conclut comme dans la preuve du théorème (\ref{cohosc}) que 
$$\H^1(X)\simeq H^1_c(X).$$
De la même façon, nous obtenons :
$$\H^{2n-1}(X)\simeq H^{2n-1}(X).$$
On va grâce à cette suite interpréter 
$\H^k(X)$ à l'aide de la cohomologie d'un sous complexe de
$C^\infty(\Lk K)$. Notons $$Y=\bigcup_{i,n_i>1} Y_i/B_i$$ et $$Y_P=\bigcup_{i}Y_i/B_i\cup \bigcup_{i,n_i=1}
(Y_i\times \bS_{V_i})/B_i$$ on dispose d'applications 
$$f\,:\, Y\rightarrow \partial K\subset K\ \et f_P\,:\, Y_P\rightarrow \partial K\subset K.$$
On note $H^\bullet(K,\ker f^*)$ (resp. $H^\bullet(K,\ker f_P^*)$) la cohomologie du sous complexe de
 $C^\infty(\Lambda^\bullet T^* K)$ formé par les formes nulles lorsqu'elle sont tirées en arrière par $f$
 (resp. $f_P$). Notre résultat est le suivant :
 \begin{prop} Pour $k\in [3, 2n-2]$, nous avons l'isomorphisme 
 $$\H^k(X)\simeq \ima\left( H^k(K,\ker f_P^*)\rightarrow H^k(K,\ker f^*)\right).$$
 \end{prop}
 \proof
 Soit $U\simeq [0,1[\times \partial K$ un voisinage tubulaire de $\partial K$ dans $\bar\cO$, le bord de
 $U$ est constitué de deux copies de $K$, on suppose que le bord de $D:=K \cup U$ est $\{1\}\times\partial
 K$, on peut aussi considérer les applications $f\,:\, Y\rightarrow \partial D$ et 
 $ f_P\,:\, Y_P\rightarrow \partial D$
 et on a évidemment :
  $$H^k(D,\ker f^*)\simeq H^k(K,\ker f^*) \ \et\  H^k(D,\ker f_P^*)\simeq H^k(K,\ker f_P^*).$$
  On peut aussi calculer la cohomologie du complexe des formes différentielles sur $U$ dont
le tiré en arrière par $f$ (ou par $f_P$) de leur restriction à $\{1\}\times\partial K$ est nulle.
 En effet, on dispose des suites exactes :
\begin{equation}\label{coho3}
H^{k-1}(U)\rightarrow H^{k-1}(Y)\rightarrow H^k(U,\ker f^*)\rightarrow H^k(U)\rightarrow H^k(Y)
\end{equation}
\begin{equation}\label{coho4}H^{k-1}(U)\rightarrow H^{k-1}(Y_P)\rightarrow H^k(U,\ker f_P^*)
\rightarrow H^k(U)\rightarrow H^k(Y_P). \end{equation}
Notre calcul de la cohomologie de $U$  (\ref{cohopK}) pour les degrés $k\in [1, 2n-2]$
$$H^k(U)\simeq \bigoplus_{i, k>2n_i-1} H^{k-2n_i+1}(Y_i)^{B_i}\oplus \bigoplus_{i} H^{k}(Y_i)^{B_i}$$
et les isomorphismes 
$$H^k(Y)=\bigoplus_{i,n_i>1}  H^{k}(Y_i)^{B_i}\ \et 
\ H^k(Y_P)=\bigoplus_{i, n_i=1} H^{k-1}(Y_i)^{B_i}\oplus \bigoplus_{i} H^{k}(Y_i)^{B_i}$$
impliquent que la première flèche de la suite (\ref{coho3}) est surjective dès que $k\ge 2$ ainsi dans ces
degrés :
\begin{equation}\label{coho5}
\begin{split}
H^k(U,\ker f^*)&\simeq \ker\big(H^k(U)\rightarrow H^k(Y))\\
  &\simeq \bigoplus_{i,k>2n_i-1} H^{k-2n_i+1}(Y_i)^{B_i}\oplus \bigoplus_{i,n_i=1} H^{k}(Y_i)^{B_i}
\end{split}\end{equation}
De même pour $k\ge 3$, la première flèche de la suite (\ref{coho4}) est surjective, on obtient donc pour
$k\in [3,2n-2]$ :
\begin{equation}\label{coho6}
\begin{split}
H^k(U,\ker f_P^*)&\simeq \ker\big(H^k(U)\rightarrow H^k(Y_P))\\
  &\simeq \bigoplus_{i,k>2n_i-1,n_i>1} H^{k-2n_i+1}(Y_i)^{B_i}
\end{split}\end{equation}
Mais on sait que pour $k\not=0$ alors $\H^k(Y_i)=H^k(Y_i)$.
Ainsi si on compare (\ref{coho6}) avec (\ref{cohoinfty}, i)) on obtient pour $k\in [3, 2n-2]$ :
\begin{equation}
\label{iso1}
\H^k(\cO)\simeq H^k(U,\ker f_P^*)
\end{equation} et en comparant (\ref{coho5}) avec (\ref{cohoinfty}, ii)),
on obtient pour $k\in [3, 2n-2]$ :\begin{equation}
\label{iso2}
\H^{k-1}_w(\cO)\simeq H^{k-1}(U,\ker f^*).
\end{equation}

Puisque $Y\subset Y_P$ nos obtenons le diagramme commutatif
 suivant où les flèches horizontales sont exactes :

\begin{equation*}\begin{split}
\xymatrix{
H^{k-1}(K)\oplus H^{k-1}(U,\ker f_P^*)\ar[d]\ar[r]&
H^{k-1}(U) \ar[d]\ar[r]&{}\\
H^{k-1}(K)\oplus H^{k-1}(U,\ker f^*)\ar[r]&
H^{k-1}(U) \ar[r]&{}\\
 }\\
 \xymatrix{
\ar[r]& H^k(D,\ker f_P^*)\ar[d]\ar[r]&H^{k}(K)\oplus H^{k}(U,\ker f_P^*)
\ar[d]\ar[r]&H^{k}(U)\ar[d]\\
 \ar[r]& H^k(D,\ker f^*)\ar[r]&H^{k}(K)\oplus H^{k}(U,\ker f^*)
\ar[r]&H^{k}(U)\\
 }\\
 \end{split}\end{equation*}
 De plus nos calculs montrent que la première et la quatrième flèches verticales sont injectives et
 les deuxième et dernière flèches verticales sont des isomorphismes.
 Alors on se sert du lemme suivant qu'on montrera un peu plus loin :
 \begin{lem}\label{diagramme}
Dans le diagramme commutatif suivant :
 $$
\xymatrix{
A\ar[d]\ar[r]&B\ar[d]\ar[r]&C\ar[d]\ar[r]&D\ar[d]\ar[r]&E\ar[d]\\
 A'\ar[r]&B'\ar[r]&C'\ar[r]&D'\ar[r]&E'\\}$$
on suppose que les flèches horizontales sont exactes,
 que les quatrième  et dernière flèches verticales sont injectives
 et que
 la deuxième flèche verticale est surjective alors la suite suivante est exacte :
 $$A'\rightarrow B'\simeq\ima(B\rightarrow B') \rightarrow \ima(C\rightarrow C')\rightarrow
  D\simeq\ima(D\rightarrow D')
 \rightarrow E\simeq\ima(E\rightarrow E').$$
 \end{lem}
 Pour en déduire que la suite suivante est exacte :
\begin{equation*} \begin{split}
 H^{k-1}(K)\oplus H^{k-1}(U,\ker f^*)&\rightarrow
H^{k-1}(U)\rightarrow\\
\rightarrow \ima\big[H^k(D,\ker f_P^*)\rightarrow H^k(D,\ker f^*)]&\rightarrow H^{k}(K)\oplus H^{k}(U,\ker f_P^*)
\rightarrow H^{k}(U) \\
\end{split}\end{equation*}

En comparant cette suite exacte à notre suite exacte :
$$ 
 H^{k-1}(D)\oplus H_w^{k-1}(\cO )\rightarrow
H^{k-1}(U)\rightarrow\H^k(X)\rightarrow H^{k}(D)\oplus \H^{k}(\cO)
\rightarrow H^{k}(U) ,$$ à l'aide des isomorphismes (\ref{iso1},\ref{iso2}), on en déduit
 l'isomorphisme voulu. \endproof
 Prouvons maintenant le lemme (\ref{diagramme}).
 \proof Donnons des noms aux différentes flèches de ce diagramme commutatif :
  $$
\xymatrix{
A\ar[d]^{f_0}\ar[r]^{\varphi_1}&B\ar[d]^{f_1}\ar[r]^{\varphi_2}&
C\ar[d]^{f_2}\ar[r]^{\varphi_3}&D\ar[d]^{f_3}\ar[r]^{\varphi_4}&E\ar[d]^{f_4}\\
 A'\ar[r]^{\varphi'_1}&B'\ar[r]^{\varphi'_2}&C'\ar[r]^{\varphi'_3}&D'\ar[r]^{\varphi'_4}&E'\\}$$
 Puisque $f_4$ est injectif, nous avons :
 
 $$\ker f_4\circ\varphi_4=\ker \varphi_4=\ker \varphi'_4\circ f_3$$
 et donc $$\ker \varphi'_4\cap \ima f_3=f_3(\ker \varphi_4)=f_3(\ima \varphi_3)=\ima\varphi_3'\circ f_2.$$
 Ceci montre l'exactitude de la dernière flèche.
 
 Remarquons que puisque $f_1$ est surjectif, l'image de $\varphi_2'$ est forcement incluse
 dans l'image de $f_2$, elle est en fait égale à l'image de $f_2\circ\varphi_2.$
 Et donc de la suite exacte
 $A'\rightarrow B'\rightarrow C'$, on en extrait la suite exacte 
 $A'\rightarrow B'\rightarrow \ima f_2$.
 
 Nous obtenons également :
 $$\ker \varphi'_3\cap \ima f_2=f_2(\ker \varphi_3)
 =f_2(\ima \varphi_2)=\ima\varphi_2'\circ f_1=\ima\varphi_2'.$$ 
 \endproof
Il ne reste que le cas de degré $2$. Dans ce cas, on sait aussi que
$$\H^1_w(\cO)=\bigoplus_{i, n_i=1} \H^1(Y_i)^{B_i}\ 
\et \ \H^2(\cO)=\bigoplus_{i, 3-2n_i>1,n_i>1} \H^{3-2n_i}(Y_i)^{B_i}=\{0\}.$$
  Si pour chaque $i$, on considère un point $p_i\in \bS_{V_i}$, alors 
  un petit calcul montre que la cohomologie de $U$ relative à 
$Z=\partial K\setminus ( \bigcup_{i,n_i=1}( Y_i\times\{B_i. p_i\})/ B_i )$ vaut en degré $1$ :
$$ H^1(U,Z)=\bigoplus_{i, n_i=1} H^1(Y_i)^{B_i}.$$
Les mêmes arguments montrent à l'aide du lemme (\ref{diagramme}) que 
$$\H^2(X)\simeq \ima\left[ H^2(K,\partial K)\rightarrow  H^2(K,Z)\right].$$
Ainsi nous avons obtenu :

\begin{thm} Soit $(X,g)$ une variété QALE asymptote à $\C^n/G$ où $G\subset \Su(n)$ telle que les
singularités de $\bS^{2n-1}/G$ soient isolées alors 
 $$ \H^k(X)\simeq \left\{
 \begin{array}{ll}
 H^k(K,\partial K)\simeq H^k_c(X)& {\rm si\ } k\le 1\\
\ima\left( H^2(K,\partial K)\rightarrow  H^2(K,Z)\right) & {\rm si\ } k=2\\
\ima\left( H^k(K,\ker f_P^*)\rightarrow H^k(K,\ker f^*)\right)& {\rm si\ } k\in [3, 2n-2]\\
 H^k(K)\simeq H^k(X)& {\rm si\ } k\ge 2n-1\\
 \end{array}\right.$$
 
 Si de plus la dimension des singularités de $\bS^{2n-1}/G$ est supérieure ou égale à $3$ alors :
 $$ \H^k(X)\simeq \left\{
 \begin{array}{ll}
 H^k(K,\partial K)\simeq H^k_c(X)& {\rm si\ } k\le 2\\

 H^k(K,\ker f^*)& {\rm si\ } k\in [3, 2n-2]\\
 H^k(K)\simeq H^k(X)& {\rm si\ } k\ge 2n-2\\
 \end{array}\right.$$
\end{thm}

 \subsection{Cohomologie $L^2$ à poids} Si on considère le poids $\mu$ tel que $\mu =\|x\|^{2a}$
  en dehors d'un compact, alors nous pouvons également  déterminer la cohomologie $L^2$ à
 poids des variétés QALE asymptotes à $\C^n/G$ où les singularités de $\bS^{2n-1}/G$ sont isolés,
  lorsque $a>n$, on trouve  $$\H^k_\mu(X)\simeq H_c^k(X)$$ et pour $a<-n$, nous obtenons
 $$\H_\mu^k(X)\simeq H^k(X).$$

\end{document}